\documentclass[11pt]{article}\input{amssym.def}\input{amssym}

\newenvironment{annotacia}{\centerline{\sc Abstract}\vspace{2mm}\narrower\narrower\sf}

\def\theequation{\thesubsection.\arabic{equation}}
\makeatletter\@addtoreset{equation}{subsection}\@addtoreset{equation}{section}\makeatother

\def\nn{\nonumber}\def\lb{\label}\def\nin{\noindent}\def\Mt{M \raisebox{1mm}{$\intercal$}}
\def\be{\begin{equation}}\def\ee{\end{equation}}\def\ba{\begin{eqnarray}}\def\ea{\end{eqnarray}}
\def\tr{{\rm Tr}\,}\def\Tr#1{{\rm Tr}_{\! R^{\mbox{\scriptsize$(#1)$}}}}
\def\TR#1#2{{\rm Tr}_{\! #2^{\mbox{\,\scriptsize$(#1)$}}}}\def\str#1{\rule[#1mm]{0pt}{1mm}}

\def\cR{{\cal R}}\def\cF{{\cal F}}\def\ot{\otimes}

\def\a{\alpha}\def\b{\beta}\def\g{\gamma}\def\d{\delta}\def\D{\Delta}\def\e{\epsilon}
\def\id{{\mathrm{id}}}\def\vf{v_{_\cF}}\def\ti{\tilde}

\newcounter{theorem}\makeatletter
\@addtoreset{theorem}{section}\makeatother

\newtheorem{prop}[theorem]{Proposition}\newtheorem{rem}[theorem]{Remark}\newtheorem{lem}[theorem]{Lemma}
\newtheorem{def-lem}[theorem]{Definition-Lemma}\newtheorem{def-prop}[theorem]{Definition-Proposition}
\newtheorem{defin}[theorem]{Definition}\newtheorem{theor}[theorem]{Theorem}
\newtheorem{cor}[theorem]{Corollary}

\begin{document}

\title{ Orthogonal and Symplectic Quantum Matrix Algebras\\
and Cayley-Hamilton Theorem for them}
\author{ \rule{0pt}{7mm} Oleg Ogievetsky\thanks{oleg@cpt.univ-mrs.fr}\\
{\small\it Center of Theoretical Physics, Luminy, 13288 Marseille, France}\\[-3pt]
{\small \&}\\[-5pt]
{\small\it P.N. Lebedev Physical Institute, Theoretical Department}\\[-2pt]
{\small\it Leninsky prospekt 53, 117924 Moscow, Russia}\\
\rule{0pt}{7mm} Pavel Pyatov\thanks{pyatov@thsun1.jinr.ru}\\
{\small\it Max Planck Institute for Mathematics}\\[-2pt]
{\small\it Vivatsgasse 7, D-53111 Bonn, Germany}\\[-2pt]
{\small \&}\\[-5pt]
{\small\it Bogoliubov Laboratory of Theoretical Physics}\\[-2pt]
{\small\it JINR, 141980 Dubna, Moscow region, Russia}}
\date{}
\maketitle

\begin{annotacia} For families of orthogonal and symplectic types quantum matrix (QM-) algebras, we
derive corresponding versions of the Cayley-Hamilton theorem. For a wider family of
Birman-Murakami-Wenzl type QM-algebras, we investigate a structure of its characteristic subalgebra
(the subalgebra in which the coefficients of characteristic polynomials take values). We define 3 sets
of generating elements of the characteristic subalgebra and derive recursive Newton and Wronski
relations between them. For the family of the orthogonal type QM-algebras, additional reciprocal relations
for the generators of the characteristic subalgebra are obtained. \end{annotacia}

\newpage

\tableofcontents

\bigskip

\section{Introduction}\lb{sec1}

A notion of a quantum matrix group, also called the RTT-algebra, is implicit in structures of
the quantum inverse scattering method. It has been given a formal definition in the works of
V. Drinfel'd, L. Faddeev, N. Reshetikhin and L. Takhtajan \cite{D1,FRT}. Since then, various aspects
of the quantum matrix group theory have been elaborated, especially in attempts to
define differential geometric structures on non-commutative spaces (see, e.g. \cite{Man,SchWZ}\footnote{
Our bibliography on the subject is not exhaustive here and below.
We give only those references, which are close in ideas and technique to our present work
and may be considered as our sources of motivation.}). In particular, a different family of algebras
generated by matrix components -- so-called reflection equation (RE-) algebras \cite{C,KS},
has been brought into consideration. Soon it was realized that, for both the RTT- and the RE-algebras,
some of the basic concepts of the classical matrix algebra, like the notion of the spectral invariants
and the characteristic identity (the Cayley-Hamilton theorem) can be properly generalized (see
\cite{EOW,NT,PS,Zh}). So, it comes out that the matrix notation used for the definition of
the RTT- and the RE-algebras is not only technically convenient, but it dictates certain structure
properties for the algebras themselves. It is then natural to search for a possibly most general
algebraic setting for the matrix-type objects. Such family of algebras was introduced in refs.\cite{Hl}
and \cite{IOP1}, and in the latter case the definition was dictated by a condition that
the standard matrix theory statements should have their appropriate generalizations. These algebras were
called quantum matrix (QM-) algebras although one should have in mind that the QM-algebras
are generated by the matrix components rather than by the matrix itself.

\smallskip
The RTT- and the RE-algebras are probably the most important subfamilies from a variety of QM-algebras.
They are distinguished both from the algebraic point of view (the presence of additional
bi-algebra and bi-comodule structures) and from the geometric point of view (their interpretation
as, respectively, the algebras of quantized functions and of quantized invariant differential
operators on a group). However, for the generalization of the basic matrix algebra statements,
it is not only possible but often more clarifying to use a weaker structure settings of the QM-algebras.

\smallskip
So far, the program of generalizing the Cayley-Hamilton theorem was fully accomplished for the linear
type QM-algebras. For the $GL(m)$-type algebras, the results were described in \cite{GPS1,IOPS,IOP1}
and for the $GL(m|n)$-type algebras in \cite{GPS2,GPS3}. These works generalize earlier results
on characteristic identities by A.J. Bracken, M.D. Gould, H.S. Green, P.D. Jarvis and R.B. Zhang
in the Lie (super)algebra case \cite{BGr,Gr,JGr} and in the quantized universal enveloping algebra case
\cite{GZB}, and by I. Kantor and I. Trishin in the matrix superalgebra case \cite{KT1,KT2}.

\smallskip
In the present work we are investigating the QM-algebras of the orthogonal and symplectic types.
Our main goal is to derive, for the quantum matrices of these types,
the characteristic identities, to use them for a
description of the spectra of the quantum matrices and to establish analogs of the Newton and Wronski
relations. To be able to do that, we have to carry out in the sections 2 and 3 some preparatory work.

\smallskip
In the section \ref{sec2} we collect few basically known results about Birman-Murakami-Wenzl (BMW)
algebras. The role that the BMW algebras play in the orthogonal/symplectic cases is similar to that
of the Hecke algebras in the general linear case. Namely, given  an R-matrix representation of the
BMW-algebra, one can define the corresponding QM-algebra and then specific properties of the R-matrix
representation define the type of the QM-algebra and hence control the form of the corresponding
characteristic identity. In the subsection \ref{subsec2.2} we introduce three sets of the idempotents
in the BMW algebra --- called antisymmetrizers, symmetrizers and contractors --- which are
used later for an identification of a specific type of the R-matrix representation.

\smallskip
In the section \ref{sec3} we consider the R-matrix representations of the BMW algebras. We define
standard notions of an R-trace~\footnote{This operation is also called a quantum trace or, shortly,
a $q$-trace in the literature.}, of a skew-invertibility, of a compatible pair of R-matrices (the
subsection \ref{subsec3.1}) and of a twisted R-matrix (the
subsection \ref{subsec3.2}). To investigate the skew-invertibility of the
R-matrix after a twist, we derive an expression for the twisted R-matrix, which is different from the standard one.

\smallskip
Next we describe the
BMW type R-matrices (the subsection \ref{subsec3.3}) and then the R-matrices of
the orthogonal and symplectic subtypes
(the subsection \ref{subsec3.4}).
Most part of a technical preparative work is done in
the subsections \ref{subsec3.2}, \ref{subsec3.3} and \ref{subsec3.5}. Here we develop the R-matrix
technique, which is later used in the main sections \ref{sec4}--\ref{sec6}.

\smallskip
A very common (although not the only) source of the skew-invertible R-matrices and of
their compatible pairs is provided by
the theory of the quasitriangular Hopf algebras.
In the appendix, we comment on the universal (i.e., quasitriangular Hopf algebraic)
counterparts of the matrix relations derived in section \ref{sec3}.

\smallskip
In the beginning of the section \ref{sec4} we introduce the QM-algebras of different (generic, BMW,
orthogonal, symplectic) types. We then define the characteristic subalgebra of the QM-algebra. It is
the subalgebra to which the coefficients of the Cayley-Hamilton identity belong. As it was shown
in \cite{IOP1}, the characteristic subalgebra is abelian. In the subsection \ref{subsec4.2} we
describe three generating sets for the characteristic subalgebra of the BMW-type QM-algebra. In the
subsection \ref{subsec4.3} we derive reciprocal relations, which are specific to generators of the
characteristic subalgebra in the orthogonal case.
Particular examples (and for particular, so-called, "standard"  R-matrices)
of the reciprocal equalities
relating the quantum determinants and the
2-contractions were obtained for the RTT-algebras by T. Hayashi \cite{Hay}.

\smallskip
Next, in the subsection \ref{subsec4.4}, we construct another important requisite of the Cayley-Hamilton
theorem --- the quantum matrix product "$\star $". The definition is given for the generic QM-algebra.
In general, the $\star \, $-product is different from the usual matrix product (however, for the subfamily
of RE-algebras, the $\star \, $-product coincides with the matrix product) but it obeys the same
algebraic properties. The $\star \, $-product is proven to be associative and hence the $\star \,$-powers
of the same quantum matrix $M$ commute. We determine then the set of all different descendants of  the
quantum matrix $M$ in the BMW case and prove that this set is $\star \, $-commutative. It turns out that,
unlike the general linear (Hecke) case, it is not possible to express all these descendants in terms of the
$\star \, $-powers of $M$ only. The expressions include also a new operation "\raisebox{1mm}{$\intercal$}",
which can be treated as a "matrix multiplication with a transposition". This fact causes some technical
complications in the proof of the Cayley-Hamilton theorem, but on the other hand, it gives rise to a
diversity of the characteristic identities in the orthogonal and symplectic cases.

\smallskip
A construction of the inverse matrix of generators and a definition of two components
of the orthogonal QM-algebra close the section \ref{sec4}.

\smallskip
Given all the definitions of the section \ref{sec4}, the derivation of the main theorems in the sections
\ref{sec5} and \ref{sec6} becomes a rather technical task, with the tools developed in
the sections \ref{sec2} and \ref{sec3}. We prove the Cayley-Hamilton theorems separately
for the symplectic case and for three different orthogonal cases in theorems \ref{theorem5.4},
\ref{theorem5.6} and \ref{theorem5.8}. We also define in each case a homomorphism from the
characteristic subalgebra to the algebra of symmetric polynomials in some set of commuting variables
(see the corollaries \ref{corollary5.5}, \ref{corollary5.7} and \ref{corollary5.9}). The Cayley-Hamilton
identities under the action of these homomorphisms are completely factorized  and hence the variables
of the symmetric polynomials can be treated as quantum matrix eigenvalues.

\smallskip
The Newton and Wronski recursions relating three different generating sets of the characteristic
subalgebra of the BMW-type QM-algebra are proven in the theorem \ref{theorem6.1}. For the orthogonal and
symplectic cases considered in the section \ref{sec5}, we use the homomorphic maps  defined in the corollaries
\ref{corollary5.5}, \ref{corollary5.7} and \ref{corollary5.9} to calculate expressions for all three
generating sets of the characteristic subalgebra in terms of quantum matrix eigenvalues (see the corollary
\ref{corollary6.4}). The latter calculation was initially aimed at comparing our results with
expressions given for the power sums for the RE-algebras in a recent paper by A. Mudrov \cite{Mudr} (see
the subsection 8.3 there). Although the  derivation methods are very different the results agree up to some
obvious changes in a notation. Notice that compared to \cite{Mudr} we are working in a more general
setting. The generalization goes in several directions. First, we do not assume a ``standard''
Drinfel'd-Jimbo's form for the R-matrices defining the algebra and, moreover, we do not use any
deformation assumptions in our constructions. Next, we are working with a wider family of QM-algebras.
And, finally, we are working directly in the algebra without passing to representations\footnote{
Passing to the representations level is hardly possible except in the RE-algebra case. The reason is
that the characteristic subalgebra belongs to the center of the RE-algebra, which is not true for the
general QM-algebra.}.

\medskip
Concluding the introduction let us comment on yet another generalization of the Cayley-Hamilton theorem,
which was suggested by I.M. Gelfand, D. Krob, A. Lascoux, B. Leclerc, V.S. Retakh and J. Thibon in
\cite{G-T}. In the latter approach the authors used a notion of quasideterminants, which allowed them to
write the Cayley-Hamilton identity in a very general setting of a free associative algebra generated by
indeterminates. The matrix powers are usual (there is no R-matrix, which affects the definition of the matrix
powers in our approach). The coefficients of the characteristic identity however become diagonal matrices
rather than scalars and, moreover, the components of these diagonal matrices are rational functions rather
than polynomials in generators. This is the price to be paid for the generality. The first examples
\cite{EOW,Zh} of the characteristic identities for the RTT-algebras constructed for the standard
Drinfeld-Jimbo R-matrix were obtained namely in the form of \cite{G-T}: the matrix powers were usual, the
coefficients of the characteristic identity were diagonal matrices (however, the entries of these matrices
are polynomials in the algebra generators, which shows a specific character of the RTT-algebras). In
\cite{OV} it was shown that the characteristic identities of \cite{EOW,Zh} contain the same information as
the characteristic identities obtained in our approach \cite{IOPS}. Notice, that even for the subfamily of
the RTT-algebras constructed for a generic Hecke-type R-matrix (not by the standard one), there is no hope
to get the characteristic identity with the polynomial coefficients in the approach of \cite{G-T}.

\section*{Acknowledgments}

The authors express their gratitude to Robert Coquereaux, Nikolai Iorgov, Alexei Isaev, Andrei Mudrov,
Arun Ram and Pavel Saponov for numerous fruitful discussions and valuable remarks. The work of the first
author (O. O.) was supported by the ANR project GIMP. The work of the second author (P. P.) was
partially supported by the grant of the Heisenberg--Landau program and by the RFBR grant No. 05-01-01086.

\section{Some facts about Birman-Murakami-Wenzl algebras}\lb{sec2}

In this preparatory section we collect definitions and derive few results on the Birman-Murakami-Wenzl
algebras. We give a minimal information, which is required for the main part of the paper.
The reader will find a more detailed presentation of the Birman-Murakami-Wenzl algebras in, e.g.,
papers \cite{W} and \cite{LR}.

\subsection{Definition}\lb{subsec2.1}

The braid group ${\cal B}_{n}$, in Artin presentation, is defined by generators
$\{\sigma_i\}_{i=1}^{n-1}$ and relations
\ba\lb{braid}\sigma_i \sigma_{i+1} \sigma_i = \sigma_{i+1} \sigma_i \sigma_{i+1}\,
&& \forall\; i=1,2,\dots ,n-1\ea
and
\ba\lb{braid2}\sigma_i \sigma_j = \sigma_j \sigma_i\  && \forall\; i,j:\; |i-j|>1\ .\ea
The {\em Birman-Murakami-Wenzl (BMW) algebra} ${\cal W}_{n}(q,\mu)$ \cite{BW,M1} is a finite
dimensional quotient algebra of the group algebra ${\Bbb C}{\cal B}_{n}$. It depends on
two complex parameters $q$ and $\mu$. Let
\be\lb{kappa}\kappa_{i} := {(q1-\sigma_i)(q^{-1} 1 +\sigma_i)\over \mu (q-q^{-1})} ,
\qquad i=1,2,\dots ,n-1 \ .\ee
The quotient algebra ${\cal W}_{n}(q,\mu)$ is specified by conditions
\be\lb{bmw2}\sigma_i \kappa_i = \kappa_i \sigma_i = \mu \kappa_i\, , \qquad
\kappa_i \sigma_{i+1}^{\pm 1} \kappa_i =  \mu^{\mp 1} \kappa_i\ .\ee
The relations (\ref{bmw2}) also imply
\be\lb{cubic} (q1-\sigma_i)(q^{-1}1+\sigma_i)(\mu 1-\sigma_i) = 0\, ,\ee
\be\lb{bmw7} e(\sigma_i-\raisebox{1pt}{$(q-q^{-1})$}1) \kappa_{i+1}
(\sigma_i-\raisebox{1pt}{$(q-q^{-1})$}1) =(\sigma_{i+1}-\raisebox{1pt}{$(q-q^{-1})$}1) \kappa_i
(\sigma_{i+1}-\raisebox{1pt}{$(q-q^{-1})$}1)\, ,\hspace{15mm}\ee
\be\lb{bmw3} \kappa_i \sigma_{i+1}^{\pm 1}  =\kappa_i \kappa_{i+1}\sigma_i^{\mp 1}\,
,\qquad\qquad\;\, \sigma_{i}^{\pm 1} \kappa_{i+1}  =
\sigma_{i+1}^{\mp 1} \kappa_{i} \kappa_{i+1}\, ,\ee
\be\lb{bmw5} \kappa_i \kappa_{i+1} \kappa_i = \kappa_i\, ,\qquad\qquad\qquad\quad
\kappa_i^2 =\eta\, \kappa_i\, ,\quad\; \mbox{where}\quad\;
\eta:= {\displaystyle {(q-\mu)(q^{-1}+\mu)\over\mu (q-q^{-1})}}\, .\ee
The parameters $q$ and $\mu$ of the BMW algebra are taken in domains
$q\in{\Bbb C}\backslash \{0,\pm 1\}$\footnote{For particular values  $\mu=\pm q^{i}$, $i\in {\Bbb Z}$,
the limiting cases $q\rightarrow \pm 1$ to the Brauer algebra
\cite{Br}
can be consistently defined.} and
$\mu\in{\Bbb C}\backslash\{0,q,-q^{-1}\}$, so that the elements $\kappa_i$ are well defined and
non-nilpotent. Further restrictions on $q$ and $\mu$ will be imposed in the subsections \ref{subsec2.2} and
\ref{subsec3.4}.

\medskip
The braid groups and their quotient BMW algebras admit a chain of monomophisms
\be\begin{array}{l}{\cal B}_2\hookrightarrow\dots\hookrightarrow {\cal B}_n\hookrightarrow
{\cal B}_{n+1}\hookrightarrow\dots \ \ ,\\[1em]{\cal W}_2\hookrightarrow\dots\hookrightarrow
{\cal W}_n\hookrightarrow {\cal W}_{n+1}\hookrightarrow\dots\end{array}\lb{h-emb}\ee
defined on the generators as
\be\lb{h-emb2}{\cal B}_{n}\ ({\mathrm or}\ \ {\cal W}_{n})\ni \sigma_i \mapsto
\sigma_{i+1}\in {\cal B}_{n+1}\ ({\mathrm or}\ \ {\cal W}_{n+1})\;\; \forall\; i=1,\dots ,n-1.\ee
We denote by $\alpha^{(n)\uparrow i}\in {\cal B}_{n+i}\ ({\mathrm{or}}\ \ {\cal W}_{n+i})$
an image of an element $\alpha^{(n)}\in {\cal B}_n\ ({\mathrm{or}}\ \ {\cal W}_{n})$
under a composition of the mappings (\ref{h-emb})--(\ref{h-emb2}). Conversely, if for some
$j<(n-1)$, an element $\alpha^{(n)}$ belongs to the image of ${\cal B}_{n-j}\ ({\mathrm{or}}\
\ {\cal W}_{n-j})$ in ${\cal B}_n\ ({\mathrm{or}}\ \ {\cal W}_{n})$ then by $\alpha^{(n)\downarrow j}$
we denote the preimage of $\alpha^{(n)}$ in ${\cal B}_{n-j}\ ({\mathrm{or}}\ \ {\cal W}_{n-j})$.

\medskip
This notation will be helpful in the next subsection where we discuss
three distinguished sequences of idempotents in the BMW algebras.

\subsection{Baxterized elements, (anti)symmetrizers and contractors}\lb{subsec2.2}

A set of elements $\sigma_i(x)$, $i=1,2,\dots, n-1,$ depending on a complex parameter $x$, in a
quotient of the group algebra ${\Bbb C}{\cal B}_{n}$ is called a set of {\em baxterized elements} if
\ba\lb{bYBE}\sigma_i(x)\, \sigma_{i+1}(xy)\, \sigma_i(y)\, =\,
\sigma_{i+1}(y)\, \sigma_i(xy)\, \sigma_{i+1}(x)\, \ea
for $i=1,2,\dots, n-1$ and
\ba\lb{bcomm}
\sigma_i(x)\, \sigma_j(y)\, =\,\sigma_j(y)\, \sigma_i(x)\,  \ea
if $|i-j|>1$.

\begin{lem}\lb{lemma2.1}
{\rm\bf \cite{J,I}} For the algebra ${\cal W}_{n}(q,\mu)$, the baxterized elements exist. There are two
sets of the baxterized elements $\{\sigma^{\varepsilon}_i\}$, $\varepsilon =\pm 1$, given by
\be\lb{ansatz}\sigma_i^{\varepsilon}(x)\, :=\,  1\, +\, {x-1\over q-q^{-1}}\, \sigma_i\, +\,
{x-1\over \alpha_{\varepsilon} x+1}\, \kappa_i\, ,
\ee
where $\alpha_{\varepsilon}\, :=\, -\varepsilon q^{-\varepsilon} \mu^{-1}$.
\end{lem}

The complex argument $x$, traditionally called {\em the spectral parameter},
is chosen in a domain $\ $ ${\Bbb C}\setminus\{-\alpha_{\varepsilon}^{-1}\}$.

\medskip
There exists an algebra isomorphism
\be\lb{homS-A}\iota :\ {\cal W}_n (q,\mu) \rightarrow {\cal W}_n (-q^{-1},\mu)\ ,\
\sigma_i\mapsto\sigma_i\ \ .\ee
The map $\iota$ interchanges the
two sets in (\ref{ansatz}), $\iota (\sigma_i^{\varepsilon}(x))=\sigma_i^{-\varepsilon}(x)$.

\medskip
In terms of the baxterized generators we construct two series of elements $a^{(i)}$ and $s^{(i)}$,
$i=1,2,\dots ,n$, in the algebra ${\cal W}_n(q,\mu)$. They are defined iteratively in two ways:
\ba\lb{ind1}a^{(1)} := 1\ \ \ {\mathrm{and}}\ \ \ s^{(1)}\, :=\, 1\, ,\ \ \ \ \ \ \ \ \ \ \, &&
\\[1em]\lb{a^k}a^{(i+1)} :={q^i\over (i+1)_q} a^{(i)}\, \sigma^{-}_i(q^{-2i})\,
a^{(i)}\ \ &{\mathrm{or}}&\ \  a^{(i+1)}:=\, {q^i\over (i+1)_q} a^{(i)\uparrow 1}\, \sigma^{-}_1(q^{-2i})\,
a^{(i)\uparrow 1}\, ,
\\[1em]\lb{s^k}s^{(i+1)} :={q^{-i}\over (i+1)_q}\, s^{(i)}\, \sigma^{+}_i(q^{2i})\,
s^{(i)}\ \ \ \  &{\mathrm{or}}&\ \  s^{(i+1)}:=\, {q^{-i}\over (i+1)_q}\, s^{(i)\uparrow 1}\,
\sigma^{+}_1(q^{2i})\, s^{(i)\uparrow 1}\, ,\ea
where $i_q$ are usual $q$-numbers, $i_q\, :=\, (q^i - q^{-i})/(q-q^{-1})$. Below we show that these two
definitions coincide.

\smallskip
To avoid singularities in the definition of $a^{(i)}$ (respectively, $s^{(i)}$), $i=1,2,\dots ,n$,
we impose further restrictions on the parameters of ${\cal W}_n(q,\mu)$:
\be\lb{mu}j_q\,\neq\, 0\, , \quad \mu \neq -q^{-2j+3} \;\; (\mbox{respectively,~}\
\mu \neq q^{2j-3})\, \quad \forall\;j = 2,3,\dots , n\ .\ee

\smallskip
The elements $a^{(i)}$ and $s^{(i)}$ are called an {\em $i$-th order antisymmetrizer} and
an {\em $i$-th order symmetrizer}, respectively.

\smallskip
The second order antisymmetrizer and symmetrizer
\be\lb{as-2}\hspace{-1mm}
a^{(2)}={q\over 2_q}\sigma_1^{-}(q^{-2})={(q1-\sigma_1)(\mu 1-\sigma_1)\over 2_q (\mu+q^{-1})} , \quad
s^{(2)}={q^{-1}\over 2_q}\sigma_1^{+}(q^2)={(q^{-1}1+\sigma_1)(\mu 1-\sigma_1)\over 2_q (\mu-q)}\,\ee
are the idempotents participating in a resolution of unity in ${\cal W}_2(q,\mu)$
(c.f. with eq.(\ref{cubic})$\,$),
\be\lb{resolution}1 \, =\, a^{(2)} + s^{(2)} + \eta^{-1}\kappa_1\, ;\ee
the spectral decomposition of the generator $\sigma_1$ of ${\cal W}_2(q,\mu)$ is
\be\lb{specdec}\sigma_1 \, =\, -q^{-1} a^{(2)} +q s^{(2)} +\mu \eta^{-1}\kappa_1\, ;\ee
Under the map (\ref{homS-A}) one has $\iota(a^{(2)})=s^{(2)}$,
$\iota(s^{(2)})=a^{(2)}$ and $\iota(\kappa_1)=\kappa_1$.

\smallskip
Likewise for $a^{(2)}$ and $s^{(2)}$, one can introduce higher order analogues for the third idempotent
entering the resolution. Namely, define iteratively
\be\lb{kappa-i}c^{(2)}\, :=\, \eta^{-1}\kappa_1\, , \qquad c^{(2i+2)}\, :=\,
c^{(2i)\uparrow 1}\, \kappa_1 \kappa_{2i+1}\, c^{(2i)\uparrow 1}\, .\ee
The element $c^{(2i)}$ is called an {\em $(2i)$-th order contractor}. Main properties of
the (anti)symmetrizers and contractors are summarized below.

\begin{prop}\lb{proposition2.2}
Two expressions given for the antisymmetrizers and symmetrizers in eqs.(\ref{a^k}) and (\ref{s^k})
are identical. The elements $a^{(n)}$ and $s^{(n)}$ are central primitive idempotents
in the algebra ${\cal W}_n(q,\mu)$. One has
\ba\lb{idemp-1}a^{(n)}\sigma_i  =\sigma_i a^{(n)} = -q^{-1} a^{(n)},&&
s^{(n)}\sigma_i = \sigma_i s^{(n)} = q s^{(n)}\quad\quad\;\;\forall\; i=1,2,\dots ,n-1\,
\ea
and
\ba\lb{idemp-2}a^{(n)} a^{(m)\uparrow i} = a^{(m)\uparrow i} a^{(n)} = a^{(n)} ,&&
s^{(n)} s^{(m)\uparrow i} = s^{(m)\uparrow i} s^{(n)} = s^{(n)} \qquad\;\;\, \mbox{if~~}\; m+i\leq n\, .
\hspace{12mm}\ea
The antisymmetrizers $a^{(n)}$, for all $n=2,3,\dots$, are orthogonal to the symmetrizers
$s^{(m)}$, for all $m=2,3,\dots$ .

\medskip
The element $c^{(2n)}$ is a primitive  idempotent in the algebra ${\cal W}_{2n}(q,\mu)$
and in the algebra ${\cal W}_{2n+1}(q,\mu)$. One has
\ba\lb{idemp-c1}c^{(2n)} c^{(2i)\uparrow n-i} =c^{(2i)\uparrow n-i} c^{(2n)} = c^{(2n)}\, &&
\forall\; i=1,2,\dots , n\, ;\\[1em]\lb{idemp-c2}c^{(2n)} \sigma_i = c^{(2n)} \sigma_{2n-i}\, , \qquad
\sigma_i c^{(2n)} = \sigma_{2n-i}\, c^{(2n)} \, &&\forall\; i=1,2,\dots ,n-1\, \ea
and
\ba\lb{idemp-c3}c^{(2n)}\sigma_n = \sigma_n\, c^{(2n)} =\mu c^{(2n)}\, . &&\ea
The contractors $c^{(2n)}$ are orthogonal to the antisymmetrizers $a^{(m)}$ and to the symmetrizers
$s^{(m)}$ for all $m>n$.
\end{prop}

\nin{\bf Proof.~} The explicit formula (\ref{a^k}) for idempotents, which we call  antisymmetrizers here, appears
in \cite{TW}, although without referring to the baxterized elements (see the proof of the lemma
7.6 in \cite{TW}).\footnote{Different expressions for the antisymmetrizers and symmetrizers, which are less suitable for
our applications, were derived in \cite{HSch}.} Our proof of the formulas (\ref{idemp-1}) and (\ref{idemp-2}) relies
on the relations (\ref{bYBE}) for the baxterized generators.

\medskip
We first check that the elements $a^{(i)}$ defined iteratively by the first formula in (\ref{a^k}) satisfy
eqs.(\ref{idemp-1}) and (\ref{idemp-2}). The equalities (\ref{idemp-1}) for the antisymmetrizers are
equivalent to
$$a^{(n)}s^{(2)\uparrow i-1}=s^{(2)\uparrow i-1}a^{(n)}=
a^{(n)}c^{(2)\uparrow i-1}=c^{(2)\uparrow i-1}a^{(n)}=0\ ,$$
which, in turn, are equivalent to
\be\lb{divide}a^{(n)}\,\sigma_i^{-}(q^2) = \sigma_i^{-}(q^2)\, a^{(n)} = 0\,
\quad \forall\; i=1,2,\dots n-1\, .\ee
Indeed, the spectral decomposition of $\sigma^{-}_i(q^2)$ contains (with nonzero coefficients) only two
idempotents, $s^{(2)\uparrow i-1}$ and $c^{(2)\uparrow i-1}$: $$\sigma^{-}_i(q^2)=\, q\, 2_q\,
(s^{(2)\uparrow i-1}+\frac{1+q\mu}{q^3+\mu}\, c^{(2)\uparrow i-1})\ .$$
To avoid a singularity in the expression for $\sigma_i^-(q^2)$, we have to  assume additionally
$\mu \neq -q^3$ for the rest of the proof. However, the expressions entering the relations
(\ref{idemp-1}) and (\ref{idemp-2}) are well defined and continuous at the point $\mu=-q^3$
(unless $-q^3$ coincides with one of the forbidden by eq.(\ref{mu}) values of $\mu$), so
the validity of the relations (\ref{idemp-1}) and (\ref{idemp-2})
at the point $\mu=-q^3$ follows by the continuity.

\medskip
Notice that the equalities $a^{(n)} \sigma_i=-q^{-1} a^{(n)}$ are equivalent to the equalities
$\sigma_i a^{(n)} =-q^{-1} a^{(n)}$ due to an antiautomorphism
$$\varsigma :\ {\cal W}_n (q,\mu) \rightarrow {\cal W}_n (q,\mu)\ ,\
\sigma_i\mapsto\sigma_i\ ,\ \varsigma (xy)=\varsigma (y)\varsigma (x)\ ,$$
since $\varsigma(a^{(n)})=a^{(n)}$ by construction.

\medskip
We now prove the equalities (\ref{idemp-1}) and (\ref{idemp-2}) by induction on $n$.

\smallskip
{}For $n=2$,~ $a^{(2)} \sigma_1  = -q^{-1} a^{(2)}$,~ by (\ref{as-2}) and (\ref{cubic}).

\smallskip
Let us check the equalities for some fixed $n>2$ assuming that they are valid for all smaller values of $n$.
Notice that as a byproduct of the definition (\ref{a^k}) (the first equality)
and the induction assumption, eqs.(\ref{divide}) and (\ref{idemp-2}) are satisfied, respectively,
for all $i=1,2,\dots ,n-2$ and for all $m,i:\; m+i\leq n-1$. It remains to check eq.(\ref{divide}) for
$i=n-1$ and eq.(\ref{idemp-2}) for $m=n-i$. Respectively, we calculate
\ba\nonumber a^{(n)}\,\sigma^{-}_{n-1}(q^2) &\sim&a^{(n-1)}\sigma^{-}_{n-1}(q^{-2n+2})
a^{(n-1)}\,\sigma^{-}_{n-1}(q^2)\\[1em]\nonumber&\sim&(a^{(n-1)}a^{(n-2)})\,\sigma^{-}_{n-1}(q^{-2n+2})
\sigma^{-}_{n-2}(q^{-2n+4})\sigma^{-}_{n-1}(q^2)\, a^{(n-2)} \hspace{10mm}\\[1em]\nonumber&=&
(a^{(n-1)}\sigma^{-}_{n-2}(q^2))\,\sigma^{-}_{n-1}(q^{-2n+4})
\sigma^{-}_{n-2}(q^{-2n+2})\, a^{(n-2)} = 0 \ ,\ea
("$\sim$" means "proportional") and
\ba\nonumber a^{(n)}\, a^{(n-i)\uparrow i}&=& {q^{n-i-1}\over (n-i)_q}\,
(a^{(n)} a^{(n-i-1)\uparrow i})\,\sigma^{-}_{n-1}(q^{-2(n-i-1)})\,a^{(n-i-1)\uparrow i}\\[1em]\nonumber
&=&{q^{n-i-1}\over (n-i)_q}\,a^{(n)}\, (1 + q^{i-n}(n-i-1)_q )\, a^{(n-i-1)\uparrow i} = a^{(n)}\ .\ea
Here in both cases, the definition of antisymmetrizers (\ref{a^k}) (the first equality), the induction
assumption and eq.(\ref{bYBE}) were used.  The centrality and primitivity of the idempotents
$a^{(n)}\in {\cal W}_n(q,\mu)$ follow then from eqs.(\ref{idemp-1}).

\medskip
To prove equivalence of the two expressions for antisymmetrizers given in eq.(\ref{a^k}), consider an
inner algebra automorphism\footnote{$\tau^{(n)}$ is the lift of the longest element of the symmetric
group. The equality $\tau^{(n)}\, \sigma_i\, (\tau^{(n)})^{-1}=\sigma_{n-i}$ holds already in the group
algebra ${\Bbb C}{\cal B}_{n}$. It follows that $(\tau^{(n)})^2$ is central in ${\Bbb C}{\cal B}_{n}$.}
\be\lb{innalis}\begin{array}{rcl} {\cal W}_n (q,\mu) \rightarrow {\cal W}_n (q,\mu)\, :&&
\sigma_i\mapsto \tau^{(n)}\, \sigma_i\, (\tau^{(n)})^{-1}=\sigma_{n-i}\ ,\\[1em]
\mbox{where}\hspace{2mm} && \tau^{(1)}=1, \;\;\; \tau^{(j+1)}=
\tau^{(j)}\,\sigma_j\sigma_{j-1}\dots\sigma_1\, .\end{array}\ee
Under this automorphism the first expression in eq.(\ref{a^k}) transforms into the second one.
However, the elements $a^{(n)}$ are central, so they do not change under
conjugation with $\tau^{(n)}$, which proves consistency of the equalities (\ref{a^k}).

\medskip
All the assertions concerning the symmetrizers follow from the relations for the antisymmetrizers by
an application of the map (\ref{homS-A}) because $\iota(a^{(n)})=s^{(n)}$ and $\iota(s^{(n)})=a^{(n)}$.

\medskip
The orthogonality of the antisymmetrizers and the symmetrizers is a byproduct of the relations (\ref{idemp-1}):
$$-q^{-1}a^{(n)}s^{(m)}= (a^{(n)}\sigma_1) s^{(m)}=a^{(n)}(\sigma_1 s^{(m)})=qa^{(n)}s^{(m)}\ .$$

\medskip
The equalities (\ref{idemp-c1}) can be proved by induction on $n$. They are obvious in the case $n=1$. Let us
check them for some fixed $n\geq 2$, assuming they are valid for all smaller values of $n$. Notice that
the iterative definition (\ref{kappa-i}) together with the induction assumption approve
eqs.(\ref{idemp-c1}) for all values of index $i$, except $i=n$. Checking the case $i=n$ splits in two
subcases: $n=2$ and $n>2$. In the subcase $i=n=2$, we have $c^{(4)}=\eta^{-2}\kappa_2\kappa_3\kappa_1
\kappa_2$ and
$$\left(c^{(4)}\right)^2 =\eta^{-4} \kappa_2\kappa_3\kappa_1\kappa_2^2\kappa_3\kappa_1\kappa_2 =
\eta^{-3}\kappa_2\kappa_3(\kappa_1\kappa_2\kappa_1)\kappa_3\kappa_2 =
\eta^{-3}\kappa_2\kappa_3\kappa_1\kappa_3\kappa_2 =\eta^{-2}\kappa_2\kappa_3\kappa_1\kappa_2 = c^{(4)}\, ,
$$
while in the subcase $i=n>2$, the calculation is carried out as follows
\ba\nonumber\left(c^{(2n)}\right)^2 &=&c^{(2n-2)\uparrow 1}\kappa_1\kappa_{2n-1}c^{(2n-2)\uparrow 1}
\kappa_1\kappa_{2n-1}c^{(2n-2)\uparrow 1}\\[1em]\nonumber&=&
\left(c^{(2n-2)\uparrow 1}c^{(2n-4)\uparrow 2}\right) (\kappa_1\kappa_2\kappa_1)
(\kappa_{2n-1}\kappa_{2n-2}\kappa_{2n-1}) \left(c^{(2n-4)\uparrow 2}c^{(2n-2)\uparrow 1}\right)\\[1em]
\nonumber &=&c^{(2n-2)\uparrow 1}\kappa_1\kappa_{2n-1}c^{(2n-2)\uparrow 1}= c^{(2n)}\, .\ea
Here in both calculations we used the definition (\ref{kappa-i}), the induction assumption and the
relations (\ref{bmw5}).

\medskip
Taking into account the relations (\ref{idemp-c1}), one can derive an alternative expression for the contractors
\be\begin{array}{ccl} c^{(2i)}&=&c^{(2i-2)\uparrow 1}\kappa_1\kappa_{2i-1}c^{(2i-2)\uparrow 1}
=c^{(2i-2)\uparrow 1}\kappa_1\kappa_{2i-1}c^{(2i-4)\uparrow 2}\kappa_2\kappa_{2i-2}c^{(2i-4)\uparrow 2}
\\[1em] &=&
(c^{(2i-2)\uparrow 1}c^{(2i-4)\uparrow 2})\kappa_1\kappa_{2i-1}\kappa_2\kappa_{2i-2}c^{(2i-4)\uparrow 2}
=c^{(2i-2)\uparrow 1}\kappa_{2i-1}\kappa_{2i-2}\kappa_1\kappa_2c^{(2i-4)\uparrow 2}\\[1em]
&=&\dots =c^{(2i-2)\uparrow 1}\left(\kappa_{2i-1}\kappa_{2i-2}\dots\kappa_{i+1}\right)
\left(\kappa_1\kappa_2\dots\kappa_{i-1}\right)c^{(2)\uparrow i-1}\\[1em] &=&
\eta^{-1} c^{(2i-2)\uparrow 1}\left(\kappa_{2i-1}\kappa_{2i-2}\dots\kappa_{i+1}\right)
\left(\kappa_1\kappa_2\dots\kappa_i\right) .\end{array}\lb{kappa-i2}\ee
Now, using this expression and noticing that, by eqs.(\ref{bmw3}),
$$\kappa_{i+1} \kappa_{i-1} \kappa_i \sigma_{i-1} =\kappa_{i+1} \kappa_{i-1} \sigma^{-1}_i =
\kappa_{i-1} \kappa_{i+1} \sigma^{-1}_i =\kappa_{i+1} \kappa_{i-1} \kappa_i \sigma_{i+1}\ ,$$
we conclude that the equality (\ref{idemp-c2}) is satisfied for $i=n-1$. In particular, (\ref{idemp-c2})
holds for $n=2$ and $i=1$. It is enough (by induction on $n$) to prove (\ref{idemp-c2}) for $i=1$.
Then observe, again by eq.(\ref{bmw3}), that
\be\lb{secoide}\kappa_i\kappa_{i\pm 1}\kappa_{i\pm 2}\,\sigma_i=\kappa_i\kappa_{i\pm 1}\sigma_i\,
\kappa_{i\pm 2}=\kappa_i\sigma^{-1}_{i\pm 1}\kappa_{i\pm 2}=\sigma_{i\pm 2}\, \kappa_i\kappa_{i\pm 1}
\kappa_{i\pm 2} \ .\ee
Now, for $n>2$,
\ba\nonumber c^{(2n)}\sigma_1&=&\eta^{-1} c^{(2n-2)\uparrow 1}
\left(\kappa_{2n-1}\kappa_{2n-2}\dots\kappa_{n+1}\right)\left(\kappa_1\kappa_2\dots\kappa_n\right)\sigma_1
\\[1em]\nonumber &=&\eta^{-1} c^{(2n-2)\uparrow 1}\left(\kappa_{2n-1}\kappa_{2n-2}\dots\kappa_{n+1}\right)
\sigma_3\left(\kappa_1\kappa_2\dots\kappa_n\right)\\[1em]\nonumber
&=&\eta^{-1} c^{(2n-2)\uparrow 1}\sigma_3\left(\kappa_{2n-1}\kappa_{2n-2}\dots\kappa_{n+1}\right)
\left(\kappa_1\kappa_2\dots\kappa_n\right)\\[1em]\nonumber &=&\eta^{-1} c^{(2n-2)\uparrow 1}\sigma_{2n-3}
\left(\kappa_{2n-1}\kappa_{2n-2}\dots\kappa_{n+1}\right)\left(\kappa_1\kappa_2\dots\kappa_n\right)
\\[1em]\nonumber &=&\eta^{-1} c^{(2n-2)\uparrow 1}
\left(\kappa_{2n-1}\kappa_{2n-2}\dots\kappa_{n+1}\right)\sigma_{2n-1}
\left(\kappa_1\kappa_2\dots\kappa_n\right) =c^{(2n)}\sigma_{2n-1}\ .\ea

\medskip
The relation (\ref{idemp-c3}) follows from eq.(\ref{bmw2}) and eq.(\ref{kappa-i2}) with $i=n$. Then,
orthogonality of  the contractors $c^{(2n)}$ and the (anti)symmetrizers $a^{(m)}$, $s^{(m)}$,  $m>n$ is
a corollary of eqs.(\ref{idemp-1}) and (\ref{idemp-c3}).

\medskip
A statement of the primitivity of the idempotent $c^{(2n)}\in{\cal W}_i(q,\mu)$, $i=2n,2n+1$, goes beyond
the needs of the present paper and we mention it for a sake of completeness only.
To prove the primitivity, one has to check,
using the equalities (\ref{idemp-c1})--(\ref{idemp-c3}),
that a relation
$c^{(2n)}\alpha^{(2n+1)}c^{(2n)}\sim c^{(2n)}$ is fulfilled for any element
$\alpha^{(2n+1)}\in{\cal W}_{2n+1}(q,\mu)$.
\hfill$\blacksquare$

\section{R-matrices}\lb{sec3}

Let $V$ denote a finite dimensional ${\Bbb C}$-linear space, $\dim V = \mbox{\sc n}$. Fixing some
basis $\{v_i\}_{i=1}^{\mbox{\footnotesize \sc n}}$ in $V$ we identify elements $X\in
{\rm End}(V^{\otimes n})$ with  matrices $X_{i_1 i_2 \dots i_n}^{j_1 j_2 \dots j_n}$.

\medskip
In this section we investigate properties of certain elements in ${\rm Aut}(V^{\otimes 2})$ generating
representations of the braid groups ${\cal B}_n$ or, more specifically, of the Birman-Murakami-Wenzl
algebras ${\cal W}_{n}(q,\mu)$ on the spaces $V^{\otimes n}$. Traditionally such operators are
called "R-matrices".

\subsection{Definition and notation}\lb{subsec3.1}

Let $X\in {\rm End}(V^{\otimes 2})$. For any $n=2,3,\dots$ and $1\leq m \leq n-1$, denote by $X_m$ an
operator whose action on the space $V^{\otimes n}$ is given by the matrix
\be\lb{X-k}(X_m)_{i_1 \dots i_n}^{j_1 \dots j_n}\ :=\ I_{i_1\dots i_{m-1}}^{j_1\dots j_{m-1}}\
X_{i_m i_{m+1}}^{j_m j_{m+1}}\ I_{i_{m+2}\dots i_n}^{j_{m+2}\dots j_n}\ .\ee
Here $I$ denotes the identity operator. In some cases below (see, e.g., eq.(\ref{s-inv})$\,$) we will also use
a notation $X_{mr}\in {\rm End}(V^{\otimes n})$, $1\leq m<r\leq n-1$, referring to an operator given by
a matrix
\be\lb{X-kl}(X_{mr})_{i_1 \dots i_n}^{j_1 \dots j_n}\ :=\ X_{i_m i_r}^{j_m j_r}\
I_{i_1\dots i_{m-1}i_{m+1}\dots i_{r-1} i_{r+1}\dots i_n}^{j_1\dots j_{m-1}
j_{m+1}\dots j_{r-1} j_{r+1}\dots j_n}\ .\ee
Clearly, $X_m= X_{m\,  m+1}$.

\medskip
We reserve the symbol $P$ for the permutation operator: $P(u\otimes v)= v\otimes u \;\;\; \forall\; u, v\in V\,$. Below we repeatedly make use of relations
$$P^2 = I\, ;\quad P_{12} X_{12} = X_{21} P_{12}\, \ \;\; \forall\;\  X\in {\rm End}(V\otimes V)\, ;\quad
\tr_{(1)} P_{12} = \tr_{(3)} P_{23} = I_2\, ,$$
where the symbol ${\rm Tr}_{(i)}$ stands for the trace over an $i$-th component space in the tensor power
of the space $V$.

\medskip
An operator $X\in {\rm End}(V^{\otimes 2})$ is called {\it skew invertible} if there exists an operator
${\Psi_X}\in {\rm End}(V^{\otimes 2})$ such that
\be\tr_{(2)} X_{12} {\Psi_X}_{23} =\tr_{(2)} {\Psi_X}_{12} X_{23} = P_{13}\, .\lb{s-inv}\ee
Define two elements of $\mbox{End}(V)$
\be\lb{CandD}C_X:={\rm Tr}_{(1)}{\Psi_X}_{12}\, ,\qquad D_X:={\rm Tr}_{(2)}{\Psi_X}_{12}\, .\ee
By (\ref{s-inv}),
\be\lb{traceCD-X}\tr_{(1)} {C_X}_1 X_{12} = I_2\, ,\qquad\tr_{(2)} {D_X}_2 X_{12} = I_1\, .\ee
A skew invertible operator $X$ is called {\em strict skew invertible} if one of the matrices,
$C_X$ or $D_X$, is invertible (by the lemma \ref{lemma3.5}, if $C_X$ or $D_X$ is invertible then they are
both invertible).

\medskip
An element $R\in {\rm Aut}(V^{\otimes 2})$ that fulfills condition
\be\lb{YBE}R_{1}\, R_{2}\, R_{1}\, = \, R_{2}\, R_{1}\, R_{2}\ .\ee
is called an {\em R-matrix}. Eq.(\ref{YBE}) is  called the {\em Yang-Baxter equation}.
Clearly, $P$ is the R-matrix; $R^{-1}$ is the R-matrix iff $R$ is. Any R-matrix $R$ generates
representations $\rho_R$ of the series of braid groups ${\cal B}_n$, $n=2,3,\dots$
\be\lb{rhoR}\rho_R:\, {\cal B}_n\rightarrow {\rm Aut}(V^{\otimes n})\ ,\quad
\sigma_i \mapsto \rho_R(\sigma_i) = R_i, \quad 1\leq i\leq n-1 .\ee
If additionally the R-matrix $R$ satisfies a third order {\em minimal}
characteristic polynomial (c.f. with eq.(\ref{cubic})$\,$)
\be\lb{charR}(qI-R)(q^{-1}I+R)(\mu I-R)=0\ ,\ee
and an element
\be\lb{K}K := \mu^{-1} (q-q^{-1})^{-1}\, (qI-R)(q^{-1}I+R)\ee
fulfills  conditions
\be\lb{bmwR}K_2\, K_1 \, = \,  R_1^{\pm 1}\, R_2^{\pm 1}\, K_1\, , \qquad
K_1\, K_2\, K_1 \, = \, K_1\, ,\ee
then we call $R$ an R-matrix of a {\em BMW-type}~ (c.f. with eqs.(\ref{kappa})--(\ref{bmw5});
we make a different but equivalent choice of defining relations).

{}For an R-matrix of the BMW-type, the formulas (\ref{rhoR}) define representations of the algebras ${\cal W}_n(q,\mu)
\rightarrow {\rm End}(V^{\otimes n})$, $n=2,3,\dots $ . In particular, $\rho_R(\kappa_i)=K_i$.

\medskip
An ordered pair $\{ R, F\}$ of two operators $R$ and $F$ from ${\rm End}(V^{\otimes 2})$ is
called {\em a compatible pair} if conditions
\be R_1\, F_2\, F_1\, =\, F_2\, F_1\, R_2\, ,\qquad R_2\, F_1\, F_2\, =\, F_1\, F_2\, R_1\, ,\lb{sovm}\ee
are satisfied. If, in addition, $R$ and $F$ are R-matrices, the pair $\{ R, F\}$ is called
a compatible pair of R-matrices. The equalities (\ref{sovm}) are called {\em twist relations} (on the
notion of the twist see \cite{D2,Resh2,IOP2}). Clearly, $\{ R,P\}$ and $\{ R,R\}$ are compatible pairs
of R-matrices;
pairs $\{ R^{-1},F\}$ and $\{ R,F^{-1}\}$ are compatible iff the pair $\{R,F\}$ is.

\begin{defin}\lb{definition3.1}
Consider a  space of ~$\mbox{\sc n}\times \mbox{\sc n}$ matrices
${\rm Mat}_{\mbox{\footnotesize\sc n}}(W)$, whose entries belong to some $\Bbb C$-linear space $W$.
Let $R$ be a skew invertible R-matrix. A linear map
\be{\rm Tr\str{-1.3}}_R:\; {\rm Mat}_{\mbox{\footnotesize\sc n}}(W)\,\rightarrow \, W ,\qquad
{\rm Tr\str{-1.3}}_R(M) =\sum_{i,j=1}^{\mbox{\footnotesize\sc n}}{(D_R)}_i^jM_j^i\, , \qquad
M\in{\rm Mat}_{\mbox{\footnotesize\sc n}}(W)\,,\lb{r-sled}\ee
is called an  R-trace. \end{defin}

The relation (\ref{traceCD-X}) in this notation reads
\be\lb{traceR}\Tr{2} R_{12} = I_1\, ,\ee
which determines $D_R$ in (\ref{r-sled}) uniquely provided $R$ is skew invertible.

\subsection{R-technique}\lb{subsec3.2}

In this subsection we develop a technique for dealing with the R-matrices, their compatible pairs
and the R-trace. Most of results reported here, like the lemma \ref{lemma3.5} and, in a particular
case of a compatible pair $\{ R,R\}$ -- the lemmas \ref{lemma3.2} and  \ref{lemma3.3} and the
corollary \ref{corollary3.4} -- are rather well known (see, e.g., \cite{I,O}).
However, we often use them in a more general settings and so, when necessary, we present sketches
of proofs.

The lemma \ref{lemma3.6} contains some new results.
Here we derive a different, from the standard one, expression for the
twisted R-matrix, which helps to investigate its skew-invertibility.

\smallskip
A universal (i.e., quasi-triangular Hopf algebraic) content
of the matrix relations derived in this subsection is discussed in the appendix.

\begin{lem}\lb{lemma3.2}
Let $\{X,F\}$ be a compatible pair,
where $X$ is skew invertible. Let ${\rm Mat}_{\mbox{\footnotesize\sc n}}(W)$ be as in the definition
\ref{definition3.1}. For any $M\in {\rm Mat}_{\mbox{\footnotesize\sc n}}(W)$, one has
\ba\lb{inv-trC}\tr_{(1)} \left( {C_X}_1 F_{12}^{\varepsilon}\, M_2\, F_{12}^{-\varepsilon}\right)
&=& I_2\, \tr (C_X M)\, ,\\[1em]\lb{inv-trD}\tr_{(2)} \left({D_X}_2 F_{12}^{-\varepsilon}\,
M_1\, F_{12}^{\varepsilon}\right) &=& I_1\ \tr(D_X M)\,  \quad \mbox{for}\;\; \varepsilon =\pm 1\, .\ea\end{lem}

\nin {\bf Proof.~} We use the twist relations  (\ref{sovm}) in a form
$F_{23}^{\varepsilon}\, X_{34}\, F_{23}^{-\varepsilon}\, =\, F_{34}^{-\varepsilon}\,
X_{23}\, F_{34}^{\varepsilon}\, ,\quad \varepsilon=\pm 1\, .$
Multiplying it by $({\Psi_X}_{12}{\Psi_X}_{45})$ and taking the traces in the spaces 2 and 4, we get
\be\lb{trpsifpf}\tr_{(2)}({\Psi_X}_{12}\ F_{23}^{\varepsilon}\ P_{35}\ F_{23}^{-\varepsilon})\ =\
\tr_{(4)}({\Psi_X}_{45}\ F_{34}^{-\varepsilon}\ P_{13}\ F_{34}^{\varepsilon})\ .\ee
Here eq.(\ref{s-inv}) was applied to calculate the traces.
Now taking the trace in the space number 1 or number 5, we obtain (after relabeling)
\ba\lb{tr-0}\tr_{(1)}({C_X}_1\ F_{12}^{\varepsilon}\ P_{23}\ F_{12}^{-\varepsilon})&=&
{C_X}_3\, I_2\ ,\\[1em]\lb{tr-4}\tr_{(3)}({D_X}_3\ F_{23}^{-\varepsilon}\ P_{12}\ F_{23}^{\varepsilon})
&=&  {D_X}_1\, I_2\, .\ea
These two relations are equivalent forms of eqs.(\ref{inv-trC}) and (\ref{inv-trD}). E.g., the formula
(\ref{inv-trC}) is obtained by multiplying eq.(\ref{tr-0}) by $M_3$  and taking the trace in the space 3.
\hfill$\blacksquare$

\begin{lem}\lb{lemma3.3}
Let $\{X,F\}$ be a compatible pair of skew invertible operators $X$ and $F$. Then following relations
\ba\lb{psi-C} {C_X}_1\, {\Psi_F}_{12}\, =\, F_{21}^{-1}\, {C_X}_2\, , &&
{\Psi_F}_{12}\, {C_X}_1\, =\, {C_X}_2\, F_{21}^{-1}\, ,\\[1em]\lb{psi-D}
{\Psi_F}_{12}\, {D_X}_2\, = \, {D_X}_1\, F_{21}^{-1}\, , &&
{D_X}_2\, {\Psi_F}_{12}\, =\, F_{21}^{-1}\, {D_X}_1\, \ea
hold.\end{lem}

\nin {\bf Proof.~} For a skew invertible $F$, eqs.(\ref{psi-C}) and (\ref{psi-D}) are
equivalent to eqs.(\ref{tr-0}) and (\ref{tr-4}). Let us demonstrate how  the left one of eqs.(\ref{psi-C})
is derived from eq.(\ref{tr-0}) with $\varepsilon = 1$.

Multiply eq.(\ref{tr-0}) by $(P_{23}{\Psi_F}_{24})$ from the right, take the trace in the space 2 and simplify
the result using eq.(\ref{s-inv}) for $X=F$ and properties of the permutation
$$\tr_{(1)}({C_X}_1\, P_{14}\, F_{13}^{-1})\, =\,
{C_X}_3\, \tr_{(2)}(P_{23}\, {\Psi_F}_{24})\, =\, {C_X}_3\, {\Psi_F}_{34}\, .$$
Then simplify the left hand side of the equality using the cyclic property of the trace
$$\tr_{(1)}({C_X}_1\, P_{14}\, F_{13}^{-1})\, =\, \tr_{(1)}( P_{14}\, F_{13}^{-1}\, {C_X}_1)\,
=\, F_{43}^{-1}\, {C_X}_4\, \tr_{(1)}P_{14}\, =\, F_{43}^{-1}\, {C_X}_4\, .$$
This proves the left relation in (\ref{psi-C}).\hfill$\blacksquare$

\begin{cor}\lb{corollary3.4}
Let $\{X,F\}$ and $\{Y,F\}$ be compatible pairs of skew invertible
operators  $X$,  $Y$ and $F$. Then following relations
\ba\lb{FCC} &F_{12}\, {C_X}_1 {C_Y}_2  = {C_Y}_1 {C_X}_2 F_{12}\ ,
\qquad\quad  F_{12}\, {D_X}_1 {D_Y}_2 \, =\, {D_Y}_1 {D_X}_2 F_{12}\ ,&
\hspace{10mm}\\[1em]
\lb{FCD} &F_{12}\, (C_X D_Y)_2 = (C_X D_Y)_1\, F_{12}\ ,\qquad\qquad F_{12}\,
(D_Y C_X)_1 \, =\, (D_Y C_X)_2\, F_{12}\, ,&\ea
\ba\lb{CxDf}
&{\tr}_{(1)}({C_X}_1 F_{12}^{-1}) = (C_X D_F)_2 = (D_F C_X)_2\ ,&\\[1em]
\lb{DxCf}
&{\tr}_{(2)}({D_X}_2 F_{12}^{-1}) = (C_F D_X)_1 = (D_X C_F)_1\,\ &\ea
hold.\end{cor}

\nin {\bf Proof.~} A calculation
$( F_{12}^{-1} {C_Y}_1) {C_X}_2 ={C_Y}_2 ( {\Psi_F}_{21} {C_X}_2) =
{C_Y}_2 {C_X}_1 F_{12}^{-1}={C_X}_1 {C_Y}_2 F_{12}^{-1}$
proves the left one of eqs.(\ref{FCC}). Here the relations (\ref{psi-C}) were applied.

\smallskip
A calculation $(F_{12}^{-1} {C_X}_1) {D_Y}_1 = {C_X}_2 (\Psi^F_{21} {D_Y}_1) = {C_X}_2 {D_Y}_2 F_{12}^{-1}$
proves  the left one of eqs.(\ref{FCD}). Here one uses subsequently the left equations from (\ref{psi-C})
and (\ref{psi-D}).

\smallskip
The relations (\ref{CxDf}) follow by taking  $\tr_{(2)}$ of the equations (\ref{psi-C}).

\smallskip
The rest of the relations in (\ref{FCC})--(\ref{DxCf}) are derived in a similar way.\hfill$\blacksquare$

\begin{lem}\lb{lemma3.5}
Let $X$ be a skew invertible R-matrix. Then statements\par
a) the R-matrix $X^{-1}$ is skew invertible;\par
b) the R-matrix $X$ is strict skew invertible,\par\noindent
are equivalent.

Provided these statements are satisfied, both $C_X$ and $D_X$ are invertible and one has
\ba\lb{CDinv} &&C_{X^{-1}} =\,  D_X^{-1}\, , \qquad D_{X^{-1}} =\,  C_X^{-1}\, .\ea\end{lem}

\nin {\bf Proof.~} See \cite{O}, the section 4.1, statements after eq.(4.1.77), or \cite{I}, the
proposition 2 in the section 3.1.

Under an assumption of an existence, for an R-matrix $X$, of the operators $X^{-1}$, ${\Psi_{\! X\,}}$ and
${\Psi_{\! X^{-1}\,}}$, the relations (\ref{CDinv}) were proved in \cite{Resh}.  \hfill$\blacksquare$

\medskip
Since, for a compatible pair $\{ X,F\}$, the pair $\{ X,F^{-1}\}$ is also compatible,
eqs.(\ref{CDinv}) together with eqs.(\ref{CxDf}-\ref{DxCf}) imply that $C_XC_F=C_FC_X$ and $D_XD_F=D_FD_X$.

\vskip .3cm
Let $\{R,F\}$ be a compatible pair of R-matrices. Define a {\em twisted} operator
\be\lb{R_f} R_f := F^{-1} R F\, .\ee It is well known that
$R_f$ is an R-matrix and the pair $\{ R_f , F\}$ is compatible. Therefore, one can twist again;
in \cite{IOP1} it was shown that if $F$ is skew invertible then
\be\lb{XffD}{D_F}_1 {D_F}_2 ((R_f)_f)_{12} =R_{12} {D_F}_1 {D_F}_2\, ;\ee

\begin{lem}\lb{lemma3.6}
Let $\{R,F\}$ be a compatible pair of R-matrices.
Following statements
\begin{itemize}
\item[{\em a)}] if $F$ is strict skew invertible then $R_f$ (\ref{R_f}) can be expressed in a form
\be\lb{R_f-fin} {R_f}_{12}\, =\,
\tr_{(34)}\left( F_{32}^{-1} {C_{F^{-1}}}_3 R_{34} {D_F}_4 F_{14}\right) \, ;\ee
\item[{\em b)}]
if $R$ is skew invertible and $F$ is strict skew invertible then $R_f$ is skew invertible; its skew
inverse is
\be\lb{Psi_R_f}{\Psi_{\! R_f}}_{12}\, =\,
{C_{F^{-1}}}_2\, \tr_{(34)}\left( F_{23}^{-1} {\Psi_{\! R\,}}_{34} F_{41}\right) {D_F}_1\, ;\ee
moreover, $\Psi_{\! R_f}$ can be expressed in a form
\be\lb{Psi_R_f-another} {\Psi_{\! R_f}}_{12}\, =\, {C_{F^{-1}}}_2 F_{21}{D_{F^{-1}}}_2
{\Psi_{\! R\,}}_{12}{C_F}_1 F^{-1}_{21}{D_F}_1\ ;\ee
\item[{\em c)}] under the conditions in {\em b)},
\be\lb{CDtwist} C_{R_f} = C_{F^{-1}} D_R\, C_F\, , \quad D_{R_f} = D_{F^{-1}} C_R D_F\, \ee
(thus, if, in addition to the conditions in {\em b)}, $R$ is strict skew invertible then
$R_f$ is strict skew invertible as well)
\end{itemize}
hold.\end{lem}

\nin {\bf Proof.~} To verify the assertion a) we calculate
\be\!\!\begin{array}{ccl} {R_f}_{12}&=& (F^{-1} R F)_{12}
= F^{-1}_{12}\left(\tr_{(4)} F^{-1}_{41} {C_{F^{-1}}}_4\right)(R F)_{12}\\[1em]
&=&\tr_{(4)} \left( (R F)_{41} F^{-1}_{12} F^{-1}_{41} {C_{F^{-1}}}_4\,\right) =
\left(\tr_{(3)}P_{13}\right)\tr_{(4)} \left( (R F)_{41} F^{-1}_{12} {C_{F^{-1}}}_1
{\Psi_F}_{14}\right)\\[1em] &=&\tr_{(34)} \left( (R F)_{43} F^{-1}_{32} {C_{F^{-1}}}_3
P_{13} {\Psi_F}_{14}\right) =\tr_{(3)}\left( F^{-1}_{32} {C_{F^{-1}}}_3
\underline{P_{13} \tr_{(4)}{\Psi_F}_{14}(R F)_{43}}\right) \ ,\end{array}\lb{R_f-alt} \ee
where in the second equality we used the relation (\ref{traceCD-X}) for $X=F^{-1}$; in the third
equality we applied the twist relations; in the fourth equality we applied eq.(\ref{psi-C}) for
$X=F^{-1}$ and inserted the identity operator $\tr_{(3)}P_{13}$; in the fifth equality we permuted
$P_{13}$ rightwards and then, in the sixth equality, used the cyclic property of the trace to move
$(RF)_{43}$ to the right.

\smallskip
To complete the transformation, we derive an alternative form for the underlined
expression in the last line in (\ref{R_f-alt}). Multiplying the twist relation $R_2 F_3 F_2 = F_3 F_2 R_3$ by
$({\Psi_F}_{12}{D_F}_4)$ and taking the traces in the spaces 2 and 4, we obtain (using eqs.(\ref{s-inv}) and
(\ref{traceCD-X}) for $X=F$)
$$\tr_{(2)}\left( {\Psi_F}_{12}(R F)_{23}\right) \, =\,
\tr_{(4)} \left( {D_F}_4 F_{34} P_{13} R_{34}\right) \, ,$$
which is equivalent  (multiply by $P_{13}$ from the left and use the cyclic property of the trace) to
\be\lb{alternative} P_{13}\tr_{(2)}\left( {\Psi_F}_{12}(R F)_{23}\right) \, =\,
\tr_{(4)}\left( R_{34} {D_F}_4 F_{14}\right) \, .\ee
Now, substituting (\ref{alternative}) into (\ref{R_f-alt}), we finish the transformation and
obtain (\ref{R_f-fin}).

\medskip
Given the formula for $R_f$, the calculation of ${\Psi_{R_f}}$ becomes straightforward
and one finds the formula (\ref{Psi_R_f}).

\smallskip
Thus, the skew invertibility of $R_f$ is established.

\medskip
Now we derive the expression (\ref{Psi_R_f-another}) for $\Psi_{\! R_f}$. Multiplying (\ref{trpsifpf})
with $\varepsilon =1$ by $P_{35}{D_{F^{-1}}}_5$ from the right and taking the trace in the space 5, we obtain
$$\begin{array}{rcl}\tr_{(2)}({\Psi_R}_{12}F_{23})&=&\tr_{(45)}({\Psi_R}_{45}F^{-1}_{34}P_{13}F_{34}
P_{35}{D_{F^{-1}}}_5)\\[1em] &=&\tr_{(4)}(F^{-1}_{34}P_{13}F_{34}{D_{F^{-1}}}_3{\Psi_R}_{43})\ .
\end{array}$$
Substituting this into the expression (\ref{Psi_R_f}), we find
\be\lb{lexfpsif}\begin{array}{rcl} {\Psi_{\! R_f}}_{12}\, &=&\,
{C_{F^{-1}}}_2\, \tr_{(34)}\left( F_{23}^{-1}F_{14}^{-1}P_{13} F_{14}{D_{F^{-1}}}_1 {\Psi_{\! R\,}}_{41}
\right) {D_F}_1\\[1em] &=&{C_{F^{-1}}}_2\, \tr_{(4)}\left(F_{14}^{-1}\tr_{(3)}(F_{23}^{-1}P_{13})
{}F_{14}{D_{F^{-1}}}_1 {\Psi_{\! R\,}}_{41}\right) {D_F}_1\\[1em]
&=&{C_{F^{-1}}}_2\, \tr_{(4)}\left(F_{14}^{-1}F_{21}^{-1}
{}F_{14}{D_{F^{-1}}}_1 {\Psi_{\! R\,}}_{41}\right) {D_F}_1\\[1em]
&=&{C_{F^{-1}}}_2\, F_{21}\, \tr_{(4)}\left(F_{14}^{-1}F_{21}^{-1}
{D_{F^{-1}}}_1 {\Psi_{\! R\,}}_{41}\right) {D_F}_1\\[1em]
&=&{C_{F^{-1}}}_2\, F_{21}\, {D_{F^{-1}}}_2\tr_{(4)}\left(F_{14}^{-1}{\Psi_F}_{12}
{\Psi_{\! R\,}}_{41}\right) {D_F}_1\ .\end{array}\ee
We used the Yang-Baxter equation for the operator $F$ in the fourth equality and eqs.(\ref{psi-C}) in the
fifth equality.

\smallskip
Multiplying eq.(\ref{trpsifpf}) with $\varepsilon =-1$ by ${\Psi_F}_{01}P_{13}$ from the left and
by $P_{35}{\Psi_F}_{56}$ from the right and taking the traces in the spaces 1 and 5, we find
$$F_{14}^{-1}{\Psi_F}_{12}{\Psi_{\! R\,}}_{41}={\Psi_{\! R\,}}_{12}{\Psi_F}_{41} F_{21}^{-1}\ .$$
Substituting this into the last line of (\ref{lexfpsif}), we obtain the equality (\ref{Psi_R_f-another}).

\medskip
{}Finally, expressions (\ref{CDtwist}) for $C_{R_f}$ and $D_{R_f}$ are obtained by taking the trace in
the space 1 or the space 2 of eq.(\ref{Psi_R_f}) and the subsequent use of eq.(\ref{traceCD-X}) for
$X=F^{\pm 1}$ and eqs.(\ref{CandD}), (\ref{CxDf}) and (\ref{DxCf}) for $X=R$. \hfill$\blacksquare$

\subsection{BMW type R-matrices}\lb{subsec3.3}

In this subsection we discuss the R-matrices of the BMW-type in more detail.

\smallskip
In the lemma \ref{lemma3.7} we collect additional relations specific to the BMW-type
R-matrices. Based on these formulas, we introduce, in the lemmas \ref{lemma3.8} and
\ref{lemma3.9}, an invertible operator $G\in {\rm Aut}(V)$ and linear maps $\phi$ and
$\xi$, which will be used in the section \ref{sec4} for a definition of a product of
quantum matrices and for a quantum matrix inversion.

\begin{lem}\lb{lemma3.7}
Let $R$ be a skew invertible R-matrix of the BMW-type. Then
\begin{itemize}
\item[$\bullet$] the operator $R$ is strict skew invertible;
\item[$\bullet$] the rank of the operator $K$ equals one, ${\rm rk}\, K = 1$;
\item[$\bullet$] following relations
\be\lb{traceK}\tr_{(2)} K_{12} = \mu^{-1} {D_R}_1\ \  ,\ \ \tr_{(1)} K_{12} = \mu^{-1} {C_R}_2\ ,\ee
\be\lb{traceDK} \Tr{2} K_{12} = \mu\, I_1\, ,\ee
\be\lb{traceD}\tr_{\!\! R} I =
\mu\, \eta \,\equiv\,{\displaystyle {(q-\mu)(q^{-1}+\mu)\over (q-q^{-1})}}\, ,\ee
\be\lb{C*D}
C_R D_R = \mu^2 I\, ,\ee
\be\lb{KDD} K_{12} {D_R}_1 {D_R}_2 = {D_R}_1 {D_R}_2 K_{12}\
=\ \mu^2 K_{12}\, \ee
hold.
\end{itemize}\end{lem}

\nin {\bf Proof.~} The proof of all the statements in the lemma but the last one is given in \cite{IOP3}.
The last relation (\ref{KDD}) (which, in another form, figures in \cite{IOP3}, in the proposition 2) can
be established in a following way.

The first equality in (\ref{KDD}) is a consequence of a relation
\be\lb{RDD}R_{12}{D_R}_1 {D_R}_2\, =\, {D_R}_1 {D_R}_2 R_{12}\, ,\ee
which is just the equality (\ref{FCC}) written for the pair $\{R,R\}$. Then the conditions
$K^2 \sim K$ and ${\rm rk}\, K = 1$ together imply
$K_{12} {D_R}_1 {D_R}_2 \sim K_{12} {D_R}_1 {D_R}_2 K_{12} \sim K_{12}\,$.
A coefficient of proportionality in this relation is recovered by taking the trace of it in the space 2
and the subsequent use of eqs.(\ref{traceK}) and (\ref{traceDK}). \hfill$\blacksquare$

\medskip
In \cite{IOP3}, a pair of mutually inverse matrices
\be
{E}_2 :=\tr_{(1)} (K_{12}P_{12})
\ \ {\mathrm{and}}\ \ E^{-1}_1 :=\tr_{(2)} (K_{12}P_{12})
\lb{defcaxy}\ee
were introduced (see eqs.(32) and (33) and the proposition 2 in \cite{IOP3}). In the
following lemma, we define analogues of the matrices $E$ and $E^{-1}$ for a compatible pair
$\{R,F\}$ of R-matrices. When the operator $F$ is the permutation operator, $F=P$, the
matrix $G$ of the definition-lemma \ref{lemma3.8} coincides with the matrix $E$.

\begin{def-lem}\lb{lemma3.8}
Let $\{R,F\}$ be a compatible pair of R-matrices, where
$R$ is skew-invertible of the BMW-type and $F$ is strict skew-invertible.
Define an element $G\in {\rm End}(V)$ by
\be\lb{G}G_1 \, :=\, \tr_{(23)} K_2 F_1^{-1} F_{2}^{-1}\, .\ee
The operator $G$ is invertible, the inverse operator reads
\be\lb{G-inv}G_1^{-1}\, =\, \tr_{(23)} F_2 F_1 K_2\, .\ee

\noindent
{}Following relations
\ba\lb{RF-G} &&R_{12} G_1 G_2\, =\, G_1 G_2 R_{12}\, ,  \qquad
{}F_{12}^{\varepsilon} G_1 = G_2 F_{12}^{\varepsilon}\ \ \ \mbox{for}\;\;\varepsilon=\pm 1\, ,
\\[1em]\lb{comm-G} &&[D_R, G]\, =\, [ C_F, G]\, =\, [D_F, G]\, =\, 0\, \ea
are satisfied.
\end{def-lem}

\nin {\bf Proof.~} A check of the invertibility of $G$ is a direct calculation
\be\begin{array}{ccc} G_1\, G_1^{-1} &=& (\tr_{(23)}K_2 F_1^{-1} F_2^{-1})
(\tr_{(23)}F_2 F_1 K_2)\, =\,\tr_{(23)} K_2 F_1^{-1} F_2^{-1} K_2 F_2 F_1\\[1em] &=& \tr_{(23)}K_2
{}F_1^{-1}{K_f}_2 F_1\, =\,\tr_{(23)} K_2 F_2 {K_f}_1 F_2^{-1}\, =\, \tr_{(23)} {K_f}_2 {K_f}_1 \, =\, I\, .
\end{array}\lb{gginvpr}\ee
Here in the first line we used the formulas (\ref{G}) and (\ref{G-inv}) and the property ${\rm rk}\, K=1$:
if $\Pi =|\zeta\rangle\langle\psi |$ is a rank one projector then $\tr (\Pi A)=\langle\psi |A|\zeta\rangle$
for any operator $A$ and
$$\tr (\Pi A)\,\tr (\Pi B)= \langle\psi |A|\zeta\rangle\, \langle\psi |B|\zeta\rangle=
\langle\psi |A\Pi B|\zeta\rangle =\tr (\Pi A\Pi B)\ $$
for any $A$ and $B$; in the second line of (\ref{gginvpr}) we passed from $K$ to $K_f = F^{-1} K F$ and
applied the twist relations (for $K_f$ and $F$) and the cyclic property of the trace. In the last equality
of (\ref{gginvpr}) we evaluated the traces using the relations (\ref{traceK}) and then the relation
(\ref{traceDK}) for $K_f$ (we are allowed to use these relations because $R_f$ is skew-invertible
by the lemma \ref{lemma3.6}).

\smallskip
Notice that, in view of eq.(\ref{KDD}), we can rewrite the formula for $G$  using the R-traces instead of
the ordinary ones
\be\lb{RG}G_1 \, =\, \mu^{-2}\,\Tr{23} K_2 F_1^{-1} F_{2}^{-1}\, .\ee
Applying the formula (\ref{inv-trD}) (written for $F^{\varepsilon}=X=R$) twice to this equality, we begin
our next calculation
\be\begin{array}{ccl} G_1 I_2 &=& \mu^{-2}\, \Tr{34} (R_2 R_3) K_2 F_1^{-1} F_2^{-1}
(R_3^{-1} R_2^{-1})\, =\,  \mu^{-2}\, \Tr{34} K_3 K_2 F_1^{-1} F_2^{-1}
R_3^{-1} R_2^{-1}\\[1em] &=& \mu^{-2}\, \Tr{34}K_2 F_1^{-1} F_2^{-1} K_2 K_3
\, =\, \mu^{-1}\, \Tr{3}K_2 F_1^{-1} F_2^{-1} K_2\, .\end{array}\lb{GI}\ee
Here we used eq.(\ref{bmwR}) in the last equality of the first line. In the second line we again
applied eq.(\ref{bmwR}) after moving $K_3$ to the right (for that we need eq.(\ref{KDD}) and the cyclicity
of the trace) and then we evaluated one R-trace with the help of eq.(\ref{traceDK}).

\smallskip
Now we use the formula (\ref{GI}) for the product $G_1 G_2$ in a transformation
\ba\nonumber && \hspace{-5mm} G_1 G_2 R_1\, =\, \mu^{-2}\, \Tr{34}(K_3 F_2^{-1} F_3^{-1} K_3)
(K_2 F_1^{-1} F_2^{-1} K_2) R_1\\[1em]\nonumber && =\, \mu^{-2}\,
\Tr{34}F_2^{-1} F_3^{-1} K_2 K_3 K_2 F_1^{-1} F_2^{-1} K_2 R_1\, =\,
\mu^{-2}\, \Tr{34}F_2^{-1} F_3^{-1} F_1^{-1} F_2^{-1} K_1 K_2 R_1\\[1em]\nonumber
&& =\,  \mu^{-2}\,\Tr{34}F_2^{-1} F_3^{-1} F_1^{-1} F_2^{-1} K_1 R_2^{-1}\, =\,
\mu^{-2}\,\Tr{34} K_3 F_2^{-1} F_3^{-1} F_1^{-1} F_2^{-1}  R_2^{-1}\\[1em]\nonumber
&& =\, \mu^{-2}\,\Tr{34}F_2^{-1} F_1^{-1} F_3^{-1} F_2^{-1}  R_3 K_2 K_3\, =\, \mu^{-2}
R_1\, \Tr{34}F_2^{-1} F_1^{-1} F_3^{-1} F_2^{-1} K_2 K_3\\[1em]\nonumber
&& =\, \mu^{-2} R_1\, \Tr{34}K_3 F_2^{-1} F_3^{-1} F_1^{-1} F_2^{-1}  K_2\, =\, R_1 G_1 G_2\, ,\ea
which demonstrates the first of the relations in (\ref{RF-G}). While doing this calculation, we repeatedly
used the twist relations for the pairs $\{K,F^{-1}\}$ and $\{R,F^{-1}\}$, applied the formulas (\ref{bmwR})
and exploited the cyclic property of the trace to move the operator $K_3$ to the right/left
in the fourth/fifth line, respectively.

\smallskip
The second relation in (\ref{RF-G}) is proved as follows
\be\begin{array}{c} G_1 F_1^{\varepsilon}=\mu^{-2}\, \Tr{23}\left( K_2 F_1^{-1} F_2^{-1}\right)
{}F_1^{\varepsilon}\, =\,\mu^{-2}\,\Tr{34}\left( (F_2^{-\varepsilon}F_3^{-\varepsilon})
K_2 F_1^{-1}F_2^{-1}(F_3^{\varepsilon} F_2^{\varepsilon})\right) F_1^{\varepsilon}\hspace{3mm}
\\[1em] =\mu^{-2}\,\Tr{34} (F_2^{-\varepsilon}F_3^{-\varepsilon}) F_3^{\varepsilon}
{}F_2^{\varepsilon}F_1^{\varepsilon} K_3 F_2^{-1} F_3^{-1}\, =\, F_1^{\varepsilon}\mu^{-2}\,
\Tr{34}K_3F_2^{-1}F_3^{-1}\, =\, F_1^{\varepsilon} G_2\ .\end{array}\lb{FG}\ee
Here we subsequently used eqs.(\ref{RG}) and (\ref{inv-trD}) for $X=R$, the twist relations for the
pair $\{K,F^{\varepsilon}\}$, the Yang-Baxter equations for $F$ and again eq.(\ref{RG}).

\smallskip
Vanishing of the last two commutators in eq.(\ref{comm-G}) follow from the above proved equality. To
find these commutators, transform eq.(\ref{FG}) to
\be\lb{PsiFG} G_1 {\Psi_F}_{12} \, =\, {\Psi_F}_{12} G_2\, , \qquad
G_2 {\Psi_F}_{12} \, =\, {\Psi_F}_{12} G_1\, ,\ee
(multiply (\ref{FG}) by ${\Psi_F}_{41}{\Psi_F}_{23}$ and take $\tr_{(12)}$)
and then apply the trace in the space 1 or the space 2 to these relations and compare results.

\smallskip
The first of the relations (\ref{comm-G}) is approved by a calculation
\ba\nonumber G_1 {D_R}_1 &=& \mu^{-2}\, \Tr{23} K_2 F_1^{-1} F_2^{-1} {D_R}_1\, =\,
\mu^{-2}\, \tr_{(23)} K_2 F_1^{-1} F_2^{-1} {D_R}_1{D_R}_2{D_R}_3\\[1em]\nonumber
&=& \mu^{-2}\, \tr_{(23)} {D_R}_1{D_R}_2{D_R}_3 K_2 F_1^{-1} F_2^{-1}\, =\, {D_R}_1 G_1\, .\ea
Here eqs.(\ref{RG}) and (\ref{FCC}) for $X=Y=R$ and eq.(\ref{KDD}) were used. \hfill$\blacksquare$

\begin{rem}{\rm
One can rewrite further the expression (\ref{RG}) for $G$:
\ba\nonumber G_1&=&\mu^{-2}\,\Tr{23}F_1^{-1}F_{2}^{-1}K_1=\mu^{-2}\Tr{2}F_1^{-1}{C_F}_2{D_R}_2K_1
\\[1em]\lb{anothG} &=&\Tr{2}F_1^{-1}{C_F}_2{D_R}_1^{-1}K_1=\mu^{-2}\Tr{2}F_1^{-1}{C_F}_2{C_R}_1K_1
\\[1em]\nonumber &=&\mu^{-2}{C_F}_1\Tr{2}{C_R}_2F_1^{-1}K_1={C_F}_1\tr_{(2)}F_1^{-1}K_1\ .\ea
Here we used subsequently: the twist relation, the relations (\ref{DxCf}), (\ref{KDD}), (\ref{C*D}), (\ref{FCC}) and then again (\ref{C*D}).}\end{rem}

\begin{def-lem}\lb{lemma3.9}
Let $\{R,F\}$ be a compatible pair of skew invertible R-matrices, where the operator $R$ is of the
BMW-type and the operator $F$ is
strict skew invertible. Define two endomorphisms $\phi$ and $\xi$ of the space ${\rm Mat}_{\mbox{\footnotesize\sc n}}(W)$:
\ba\lb{phi} \phi(M)_1 &:=& \Tr{2} \left( F_{12}M_1 F^{-1}_{12} R_{12}\right),
\qquad M\in {\rm Mat}_{\mbox{\footnotesize\sc n}}(W),\ea
and
\ba\lb{xi}
\xi(M)_1 &:=& \Tr{2} \left( F_{12}
M_1 F^{-1}_{12} K_{12}\right), \qquad M\in {\rm Mat}_{\mbox{\footnotesize\sc n}}(W)\, .\ea
The mappings $\phi$ and $\xi$ are invertible; their inverse mappings read
\ba\lb{phi-inv} \phi^{-1}(M)_1 &=& \mu^{-2}\TR{2}{R_f} \left(  F_{12}^{-1}
M_1 R^{-1}_{12} F_{12}\right) \ea
and
\ba\lb{xi-inv}\xi^{-1}(M)_1 &=& \mu^{-2} \TR{2}{R_f}\left( F^{-1}_{12} M_1 K_{12} F_{12}\right)\, .\ea
{}Following relations for the R-traces
\be\lb{phi-xi-traces}
{\rm Tr}_{\! R_f^{\rule{0pt}{6pt}}}\phi(M)\, =\, {\rm Tr}_{\! R^{\rule{0pt}{6pt}}} M\, , \qquad
{\rm Tr}_{\! R_f^{\rule{0pt}{6pt}}}\xi (M)\, =\, \mu\, {\rm Tr}_{\! R^{\rule{0pt}{6pt}}} M\, .\ee
are satisfied.
\end{def-lem}

\nin {\bf Proof.~} The expressions in the right hand sides of eqs.(\ref{phi-inv}) and
(\ref{xi-inv}) are well defined, since, by the lemma \ref{lemma3.6} b), the R-matrix
$R_f$ is skew invertible.

\smallskip
Let us check the relation $\phi^{-1}(\phi(M))=M$ directly.

\smallskip
Using the formulas (\ref{phi}) and (\ref{phi-inv}) and applying the relation (\ref{inv-trD}) for the
pair $\{R,F\}$ we begin a calculation
\ba\nn {\phi^{-1}(\phi(M))}_1 &=& \mu^{-2}\, \TR{2}{R_f}\left( F_{12}^{-1}
(\Tr{2'}  F_{12'} M_1 F^{-1}_{12'} R_{12'}) R^{-1}_{12} F_{12}\right)\\[1em]\nn &=&
\mu^{-2}\, \TR{2}{R_f}\Tr{3}\left( F_1^{-1} F_2^{-1} \underline{F_1} M_1
{}F_1^{-1} R_1 F_2 R_1^{-1} F_1\right)\ .\ea
In the next step we move the element $F_1$ underlined in the expression above to the left and then
transport it to the right hand side using the cyclic property of the trace (when $F_2$ moves cyclically,
$\TR{2}{R_f}\Tr{3}$ becomes $\Tr{2}\TR{3}{R_f}$ due to eq.(\ref{FCC})$\,$).
Applying the Yang-Baxter equation for $F$ and the relation (\ref{FCC}) in the case $X=R$ and $Y=R_f$, we continue
\be\begin{array}{ccl} {\phi^{-1}(\phi(M))}_1&=& \mu^{-2}\, \Tr{2}\TR{3}{R_f}\!\left(  F^{-1}_1 F^{-1}_2 M_1
\underline{F_1^{-1} R_1} F_2 \underline{R_1^{-1} F_1}F_2\right)\\[1em] &=&\mu^{-2}\, \Tr{2}\TR{3}{R_f}\!
\left(  F^{-1}_1M_1F_2^{-1} {R_f}_1\, \underline{F_1^{-1} F_2 F_1} {R_f}^{-1}_1\, F_2\right)\\[1em] &=&
\mu^{-2}\, \Tr{2}\TR{3}{R_f}\!\left(  F^{-1}_1M_1\underline{F_2^{-1}{R_f}_1F_2}\, F_1\,\underline{F_2^{-1}
{R_f}^{-1}_1\, F_2}\right)\\[1em] &=&\mu^{-2}\, \Tr{2}\!\left(  F^{-1}_1 M_1 F_1
\underline{(\TR{3}{R_f} {R_f}_2 F_1 {R_f}_2^{-1}\, )} F_1^{-1}\right)\ .\end{array}\lb{step}\ee
Here we consequently transformed the underlined expressions using the definition of $R_f$, the
Yang-Baxter equation for $F$ and the twist relations for the pair $\{R_f,F\}$. To calculate the trace
underlined in the last line of (\ref{step}), we apply the relation (\ref{inv-trD}) for the
pair $\{R_f, R_f\}$ and then use the relation (\ref{DxCf}) written for the pair $\{R_f,F^{-1}\}$.
The result reads
\ba\nn {\phi^{-1}(\phi(M))}_1&=& \mu^{-2}\, \Tr{2}\!\left(  F^{-1}_1 M_1 F_1
\underline{(D_{R_f} C_{F^{-1}})_1 F^{-1}_1}\right)\ .\ea
Now, using eq.(\ref{FCD}), written for the  pairs $\{R_f, F\}$ and $\{F^{-1},F\}$,
the relations (\ref{CDtwist}) and (\ref{CDinv}) for $X=F$, eq.(\ref{C*D}) and the relation (\ref{traceCD-X})
for $X=F^{-1}$, we complete the calculation
\ba\nn &&\hspace{-22mm}{\phi^{-1}(\phi(M))}_1\, =\, \mu^{-2}\, \tr_{(2)}\!\!\left((D_{R_f}
C_{F^{-1}} D_R)_2 F^{-1}_1\right)\! M_1\\[1em]
&&=\, \mu^{-2}\, \tr_{(2)}\!\!\left((D_{F^{-1}}
C_R D_F C_{F^{-1}} D_R)_2 F^{-1}_1\right)\! M_1
\, =\, \tr_{(2)}({D_{F^{-1}}}_2 F^{-1}_{12})M_1 = M_1\, .\ea

\smallskip
A proof of the equality $\xi^{-1}(\xi(M))=M$ proceeds quite similarly until
the line (\ref{step}), where one has to use a relation
\be\lb{traceKMK} \Tr{2}(K_1 M_1 K_1) =  (\tr_{\!\! R} M) I_1\, \quad \forall\;
M\in {\rm Mat}_{\mbox{\footnotesize\sc n}}(W)\, \ee
instead of the relation (\ref{inv-trD}). This in turn follows from eqs.(\ref{traceK}) and
(\ref{traceDK}) and the property ${\rm rk}\, K=1$.

\smallskip
The relations (\ref{phi-xi-traces}) can be directly checked starting from the definitions (\ref{phi}) and
(\ref{xi}), applying the relation (\ref{FCC}) in the case $X=R$ and $Y=R_f$ and then using the formulas (\ref{traceR})
and (\ref{traceDK}). \hfill$\blacksquare$

\begin{rem}\lb{remark3.10}
{\rm For the mapping $\phi$, the statement of the lemma \ref{lemma3.9} remains valid if one weakens
the conditions, imposed on the R-matrix $R$, replacing the BMW-type condition by the strict skew
invertibility. In this case, one should substitute the term $\mu^{-2} D_{R_f}$ by $D_{R_f^{-1}}$ in the
expression (\ref{phi-inv}) for the inverse mapping $\phi^{-1}$. The proof repeats
the proof of the formula (\ref{phi-inv}).}\end{rem}

\subsection{Orthogonal and symplectic type R-matrices}\lb{subsec3.4}

Most of results of this paper will be derived for two particular families of the BMW-type
R-matrices. We are going to describe them now.

\medskip
Consider R-matrix realizations $\rho_R(a^{(i)})$ of the antisymmetrizers  (\ref{a^k}).
We impose  additional constraints on a skew invertible BMW-type R-matrix $R$ demanding that
\be\lb{spec4} {\rm rk}\rho_R(a^{(i)})\neq 0\, \quad \forall\; i=2,3,\dots ,k\quad\mbox{and} \quad
\rho_R\left( a^{(k)\uparrow 1}\sigma^-_1(q^{-2k})a^{(k)\uparrow 1} \right)\equiv 0\,\ee
for some $k\geq 2$. Here we assume that parameters $q$ and $\mu$ fulfill conditions (c.f. with
eq.(\ref{mu})$\,$)
\be\lb{mu1} i_q\neq 0\,  \;\;\forall\; i=2,3,\dots ,k;
\qquad \mu\neq -q^{-2i+1}\,  \;\;\forall\; i=1,2,\dots ,k\, .\ee
Notice that, in the last condition in eq.(\ref{spec4}), the unnormalized form of the $(k+1)$ order
antisymmetrizer $a^{(k+1)}$ is used and so there is no need to demand $(k+1)_q\neq 0$.

\smallskip
Let us derive some consequences of the relations (\ref{spec4}). Applying ${\Tr{i}}$ to $\rho_R(a^{(i)})$
and using the relations (\ref{a^k}), (\ref{traceR}), (\ref{traceK}) and (\ref{traceD}), we calculate
\be\lb{spec1}\Tr{i} \rho_R(a^{(i)})\, =\, \delta_i\, \rho_R(a^{(i-1)})\, ,\quad \mbox{where} \quad
\delta_i\equiv\delta_i(q,\mu) :=\, - {q^{i-1}(\mu + q^{1-2i})(\mu^2 - q^{4-2i})\over
(\mu + q^{3-2i})(q-q^{-1}) i_q }\, .\ee
In view of eqs.(\ref{spec1}), the conditions (\ref{spec4}) imply, in particular, $\delta_{k+1} = 0$,
wherefrom one specifies three admissible values of $\mu$: $\mu\in\{-q^{-1-2k},\pm q^{1-k}\}$.

\smallskip
Notice that the choice $\mu = -q^{1-k}$ contradicts the conditions (\ref{mu1}) in the case when the
number $k$ is even.

In the case when $k$ is odd, the choices $\mu = -q^{1-k}$ and $\mu = q^{1-k}$ are related by a
substitution $R \mapsto  - R$. On the algebra level, this corresponds to an algebra isomorphism
${\cal W}_n(q,\mu)\rightarrow {\cal W}_n(-q,-\mu):\; \sigma_i\mapsto -\sigma_i$,~ $i=1,\dots ,n-1$. The
(anti)symmetrizers are invariant under this map: $s^{(i)}\mapsto s^{(i)}$ and $a^{(i)}\mapsto a^{(i)}$.

\smallskip
Therefore we are left with only
two essentially different choices of the parameter $\mu$: either $\mu=-q^{-1-2k}$ or $\mu=q^{1-k}$.

\smallskip
Notice that, with these
choices, the consistency of the conditions on $\mu$ in eq.(\ref{mu1}) follows from the conditions
on $q$.\footnote{We stress again that the cases $q=\pm 1$ can be consistently treated as the limiting
cases $q\rightarrow \pm 1$.}

\begin{defin}\lb{definition3.11}
Let $R$ be a skew invertible BMW-type R-matrix defining representations $\rho_R$ of the algebras
${\cal W}_i(q,\mu)$. Assume that $R$ additionally satisfies the conditions (\ref{spec4}) for some $k\geq 2$.
This implies, in particular, restrictions on $q$: $i_q\ne 0$ for $i=2,\dots ,k$. Then

\medskip
\noindent
\hspace{10mm}a)~ $R$ is called a $Sp(2k)$-type R-matrix in the case when $\mu=-q^{-1-2k}$;

\medskip
\noindent
\hspace{10mm}b)~ $R$ is called an $O(k)$-type R-matrix in the case when $\mu=q^{1-k}$
and~ $\mbox{rk}\,\rho_R(a^{(k)})=1$.\end{defin}

{}For the standard R-matrices related to the quantum groups of the series $Sp_q(2k)$ and $SO_q(k)$ \cite{FRT}, the
conditions a) and b), respectively,  as well as the relations (\ref{spec4}) are fulfilled. This explains
our notation.

\begin{rem}\lb{remark3.12}
{\rm Functions
\be\lb{Delta-i} \Delta^{(i)}(q,\mu):=\Tr{1,2,\dots ,i}\rho_R(a^{(i)}) = \prod_{j=1}^i \delta_j(q,\mu)\ee
are (up to an overall factor) particular elements of a set of rational functions $Q_\lambda(\mu^{-1},q)$
labelled by partitions $\lambda\vdash i$. The functions (\ref{Delta-i}) were introduced in the theorem 5.5 in \cite{W}
($\Delta(q,\mu)=\mu^i Q_{[1^i]}(\mu^{-1},q)$). These functions describe the q-dimensions of the
highest weight modules $V_{\lambda}$ for the orthogonal and symplectic quantum groups (see \cite{W}, the section 5
and \cite{OW}, the lemma 3.1). E.g., in the orthogonal case, the q-dimensions of the highest weight modules
$V_{[1^i]}$ are (up to an overall factor) given by a formula
\be\lb{qdim-O(k)}\Delta^{(i)} (q,q^{1-k}) =
q^{i(1-k)}\, {q^i+q^{k-i}\over q^k+1}\, {k_q!\over i_q! (k-i)_q!} \, \ ,\ee
which we shall several times use below.}\end{rem}

\begin{rem}\lb{remark3.12.1}
{\rm We note, without going into details, particular values of $k$ for the R-matrices
of the orthogonal and symplectic types.

\smallskip
For the symplectic type, the case $k=1$ is particular: the antisymmetrizer
$\rho_R(a^{(2)})$ vanishes and the spectral decomposition of $R$ contains only two
projectors:
$$R=q\,\rho_R(s^{(2)})-q^{-3}\,\rho_R(c^{(2)})\ .$$
This R-matrix is of the Hecke type $GL(2)$ (to normalize the eigenvalues, multiply by
$q$ and replace $q^2$ by $q$), which is a manifestation
of the accidental isomorphism $SL(2)\sim Sp(2)$.

Accidental isomorphisms for quantum
groups, corresponding to the standard deformation, are discussed in \cite{JO}.

\smallskip
For the orthogonal type, the case $k=4$ is particular. We just mention
that in this case the antisymmetrizer $\rho_R(a^{(2)})$ splits into two projectors
(the self-dual and anti-self-dual parts), $\rho_R(a^{(2)})=\rho_R(a^{(2)})^{(+)}
+\rho_R(a^{(2)})^{(-)}$.
They satisfy
$$\rho_R(a^{(2)})^{(\pm)}_{12}\, R_{23}\, R_{12}=
R_{23}\, R_{12}\,\rho_R(a^{(2)})^{(\pm)}_{23}$$
and
$$\rho_R(a^{(2)})^{(\pm)}_{23}\, R_{12}\, R_{23}=
R_{12}\, R_{23}\,\rho_R(a^{(2)})^{(\pm)}_{12}\ .$$
This shows that, in this case, the structure of the centralizer algebra is not of the BMW
type.

The geometry of the four-dimensional euclidean and Minkowski spaces is considered,
for the standard deformation, in \cite{OSWZ}.}\end{rem}

\subsection{Rank-one projectors}\lb{subsec3.5}

In the definition of type $O(k)$ R-matrices we demanded the operator $\rho_R(a^{(k)})
\in {\rm End}(V^{\otimes k})$ to be of rank one. In this subsection we investigate specific
consequences the rank-one condition results in.

\smallskip
We first observe that any skew invertible BMW-type R-matrix  gives rise to a series of
rank-one projectors.

\begin{lem}\lb{lemma3.13}
Assume that $R$ is a skew invertible R-matrix of the BMW-type. Then  $\rho_R(c^{(2i)})\in {\rm End}(V^{\otimes 2i})$
(for the definition of the idempotents $c^{(2i)}$, see eqs.(\ref{kappa-i})$\,$) are rank-one operators fulfilling
relations
\ba\lb{trace-c2i}\Tr{2i}\rho_R(c^{(2i)}) &=& \eta^{-1}\mu\, I_1\, \rho_R(c^{(2i-2)\uparrow 1})
\quad \forall\; i\geq 1\, ,\\[1em]\lb{traces-c2i}\Tr{i+1,i+2,\dots ,2i}\rho_R(c^{(2i)}) &=&
(\eta^{-1}\mu)^i\, I_{1,2,\dots ,i}\, .\ea\end{lem}

\nin{\bf Proof.~} We will prove eq.(\ref{trace-c2i}) performing induction on $i$. For $i=1$, the relation
(\ref{trace-c2i}) is just the equality (\ref{traceDK}). Assuming that the equation (\ref{trace-c2i}) is valid for all
$i=1,2,\dots j-1$, we check the case $i=j$
\ba\begin{array}{ccl} \Tr{2j} \rho_R(c^{(2j)}) &=& \eta^{-1}\,\rho_R(c^{(2j-2)\uparrow 1})
\left(\Tr{2j} K_{2j-1}\right) \rho_R\bigl( (\kappa_{2j-2}\kappa_{2j-3} \dots\kappa_{j+1})
(\kappa_1\kappa_2\dots \kappa_j)\bigr)\, \\[1em] &=&\eta^{-1}\mu\, \rho_R\bigl( c^{(2j-2)\uparrow 1}
\kappa_{2j-2}\kappa_{2j-3}\dots\kappa_{j+1}(\kappa_1\kappa_2\dots \kappa_j)\bigr)\, \\[1em] &=&
\eta^{-1}\mu\, \rho_R\bigl( c^{(2j-2)\uparrow 1}\kappa_{j-1}\dots\kappa_3\kappa_{2}(\kappa_1\kappa_2\dots
\kappa_j)\bigr)\, \\[1em] &=&\eta^{-1}\mu\, I_1\, \rho_R\bigl( c^{(2j-2)\uparrow 1}\kappa_{j-1}\kappa_j
\bigr)\, =\, \eta^{-1}\mu\, I_1\, \rho_R\bigl(c^{(2j-2)\uparrow 1}\bigr)\, .\end{array}\nonumber\ea
Here in the first line we substituted the expression (\ref{kappa-i2}) for $2i$-th order contractors; in the
second line we used eq.(\ref{traceDK}); in the third line we applied repeatedly eq.(\ref{idemp-c2}) for
$\kappa$'s (they are quadratic polynomials in $\sigma$'s); in the last line we used
eq.(\ref{bmw5}) several times. Then, noticing that $\kappa_i$ divides (from both left and right sides)
$c^{(2j)}$, we complete the transformation.

\smallskip
The relation (\ref{traces-c2i}) follows by a repeated application of eq.(\ref{trace-c2i}).

\smallskip
To calculate the rank of the operators $\rho_R(c^{(2i)})$, we notice that by eqs.(\ref{kappa-i}) and
(\ref{idemp-c1}) one can rewrite the formula (\ref{kappa-i}) for the contractors in a form
\be\lb{kappa-i3} c^{(2i)}\, =\, c^{(2i-2)\uparrow 1} {\displaystyle \left(\prod_{j=1}^{i} \kappa_{2j-1}\right)}
c^{(2i-2)\uparrow 1}\, .\ee
By the rank-one property of operators $K_i\in{\rm End}(V^{\otimes 2})$, a composite operator
$\rho_R(\prod_{j=1}^{i} \kappa_{2j-1}) = \prod_{j=1}^{i} K_{2j-1}$ also has rank one. Hence, the
operator $\rho_R(c^{(2i)})$ is either of rank-one or identically vanishes.
The latter possibility contradicts eq.(\ref{traces-c2i}). \hfill$\blacksquare$

\medskip
We next describe specific properties of the rank-one projectors.

\begin{lem}\lb{lemma3.14}
Let an R-matrix $R$ be of either skew invertible BMW-type or $O(k)$-type. Then, respectively,
\ba\lb{spec-c1}
&& {\displaystyle \left(\prod_{j=1}^{2i}{D_R}_j\right)}
\rho_R(c^{(2i)})\, =\, \mu^{2i}\, \rho_R(c^{(2i)})\, \ea
in the BMW case and
\ba\lb{spec-a1}
&&
{\displaystyle \left(\prod_{j=1}^{k}{D_R}_j\right)} \rho_R(a^{(k)})\, =\, q^{k(1-k)}\, \rho_R(a^{(k)})\, \ea
in the $O(k)$ case.

Let additionally $F$ be a strict skew invertible R-matrix and assume that $\{R,F\}$
is a compatible pair. Then
one has, respectively,
\ba\lb{spec-c2} && {\displaystyle \left(\prod_{j=1}^{2i}{D_{R_f}}_j\right)}
\rho_R(c^{(2i)})\, =\, \mu^{2i}\, \rho_R(c^{(2i)})\, ;\ea
in the BMW case and
\ba\lb{spec-a2} &&
{\displaystyle\left(\prod_{j=1}^{k}{D_{R_f}}_j\right)}\rho_R(a^{(k)})\, =\, q^{k(1-k)}\, \rho_R(a^{(k)})\, \ea
in the $O(k)$ case.
\end{lem}

\nin {\bf Proof.~} To prove the equality (\ref{spec-c1}), take the expression (\ref{kappa-i3}) for
the $2i$-th order contractor, change the order of commuting (by eq.(\ref{RDD})$\,$) operators
$\rho_R(c^{(2i-2)\uparrow 1})$ and $\prod_{j=1}^{2i}{D_R}_j$, and evaluate the factor
$\prod_{j=1}^{2i}{D_R}_j$ on $\rho_R(\prod_{j=1}^{i} \kappa_{2j-1}) = \prod_{j=1}^{i} K_{2j-1}$
using the relation (\ref{KDD}).

\medskip
The equality (\ref{spec-a1}) is confirmed by a calculation
$$\begin{array}{ccl} {\displaystyle \left(\prod_{j=1}^{k}{D_R}_j\right)}\rho_R(a^{(k)})=\rho_R(a^{(k)})
{\displaystyle \left(\prod_{j=1}^{k}{D_R}_j\right)} \rho_R(a^{(k)})\,\,{ =
\left(\Tr{1,2,\dots k}\rho_R(a^{(k)})\right)\rho_R(a^{(k)}) }\\[1em] =
\Delta^{(k)}(q,q^{1-k})\, \rho_R(a^{(k)}) = q^{k(1-k)}\, \rho_R(a^{(k)})\ .\end{array}\nonumber $$
Here the idempotency of $\rho_R(a^{(k)})$, its commutativity with $\prod_{j=1}^k {D_R}_j$ and its
rank-one property were taken into account in the first line; eqs.(\ref{spec1})--(\ref{qdim-O(k)})
were used for the evaluation of the R-traces in the second line.

\smallskip
The conditions of the second part of the proposition guarantee that the twisted R-matrix $R_f=F^{-1}RF$ is of
the same type as $R$ and is skew invertible (see the lemma \ref{lemma3.6}b).
Hence, the relations (\ref{spec-c1}) and (\ref{spec-a1}) are satisfied for $R_f$ as well.
Now, to prove the equalities (\ref{spec-c2}) and (\ref{spec-a2}), it is enough to demonstrate that the
two operators $\rho_{R_f}(c^{(2i)})$ and $\rho_R(c^{(2i)})$ (respectively, $\rho_{R_f}(a^{(k)})$
and $\rho_R(a^{(k)})$) are related by a similarity transformation, which commutes with
$\prod_{j=1}^{2i} {D_{R_f}}_j$ (respectively, $\prod_{j=1}^k {D_{R_f}}_j$).

\smallskip
Define a set of operators $Z^{(i)}\in {\rm Aut}(V^{\otimes i})$, $i=1,2,\dots ,$ by an iteration
\be\lb{Z-i} Z^{(1)}:=I, \qquad Z^{(2)}:=F_1,\qquad Z^{(i+1)}:=\left(F_1 F_2 \dots F_i\right) Z^{(i)} =
Z^{(i)} \left(F_i \dots F_2 F_1\right) \ee
(two expressions for the element $Z^{(i+1)}$ in eq.(\ref{Z-i}) are two different decompositions of the lift, to the
braid group, of the longest element of the symmetric group).

\smallskip
By induction, one can readily prove relations
\be\lb{Z-R} R_i\, Z^{(k)}\, =\, Z^{(k)}\, {R_f}_{(k-i)}\,\qquad \forall\; i=1,2,\dots ,k-1\, \ee
(compare with the formulas in eq.(\ref{innalis});
they correspond to a particular case, when $F=R$, of the relations (\ref{Z-R})$\,$).

E.g., for $i=2,\dots ,k-1$, one has, by induction on $k$,
\ba\nonumber R_i Z^{(k)} &=& R_i F_1 F_2 \dots F_{k-1} Z^{(k-1)} = F_1 F_2 \dots F_{k-1} R_{i-1}Z^{(k-1)}
\\[1em]\nonumber &=&F_1 F_2 \dots F_{k-1} Z^{(k-1)}{R_f}_{(k-i)} = Z^{(k)}{R_f}_{(k-i)}\ .\ea
Here we used the twist relation for the pair $\{ R,F\}$ and the first iterative definition
of the element $Z{(k)}$ in eq.(\ref{Z-i}).

{}For $i=1$, one has, again by induction on $k$,
\ba\nonumber R_1 Z^{(k)} &=& R_1 Z^{(k-1)} F_{k-1}\dots F_2 F_1 = Z^{(k-1)}{R_f}_{(k-2)}F_{k-1}\dots F_2 F_1
\\[1em]\nonumber &=&Z^{(k-1)}F_{k-1}\dots F_2 F_1{R_f}_{(k-1)}=Z^{(k)}{R_f}_{(k-1)}\ .\ea
Here we used the twist relation for the pair $\{ R_f,F\}$ and the second iterative definition
of the element $Z{(k)}$ in eq.(\ref{Z-i}).

\smallskip
The relations eq.(\ref{Z-R}) imply following equalities
\be\lb{Z-ac} \rho_R(a^{(k)})\, Z^{(k)} =  Z^{(k)} \rho_{R_f}(a^{(k)})\, ; \qquad
\rho_R(c^{(2i)})\,  Z^{(2i)} =  Z^{(2i)} \rho_{R_f}(c^{(2i)})\,\qquad\forall\; i=1,2,\dots\ee
To derive these equalities, substitute two different expressions from eq.(\ref{a^k})
for the antisymmetrizers in  terms $\rho_R(a^{(k)})$ and $\rho_{R_f}(a^{(k)})$; for the contractors,
the same expression (\ref{kappa-i}) is used in both sides of eq.(\ref{Z-ac}) since it is invariant with
respect to the reflection $\sigma_j\leftrightarrow\sigma_{i-j}$.

\smallskip
By eq.(\ref{FCC}) (for $X=Y=R_f$), the operators $Z^{(i)}$ and $\prod_{j=1}^i {D_{R_f}}_j$
commute and, hence, $Z^{(2i)}$ (respectively, $Z^{(k)}$) realizes the
similarity transformation we are looking for.\hfill$\blacksquare$

\medskip
Due to the rank one property of the operator $K$, there is an expression for the idempotents
$\rho_R(c^{(2j)})$, which is shorter than their original expression (\ref{kappa-i}) and which
fits well with the diagrammatic representation (see, e.g., \cite{BW})
of the elements of the Birman-Murakami-Wenzl algebras.

\begin{lem} Let $E$ be the matrix defined by eq.(\ref{defcaxy}). We have
\be\rho_R(c^{(2j)})=\prod_{s=1}^{j}
\Bigl( \eta^{-1}\; (E_s)^{s-j}\ K_{s,2j+1-s}\ (E_s)^{j-s}\Bigr)\ \lb{efoi}\ee
(the factors, corresponding to different values of $s$ in the right hand side of eq.(\ref{efoi}),
live in different spaces, they thus commute, so we don't need to specify the order of the product).
\lb{lemma3.14.1}\end{lem}

\nin {\bf Proof.~}
An identity
\be K_{12}K_{23}=E_3\; K_{12}P_{12}P_{13}\lb{k12k23}\ee
follows from the rank one property of the operator $K$ (written explicitly, with indices,
the relation (\ref{k12k23}) becomes evident; we have, similarly,
$K_{12}K_{31}=E^{-1}_3\; K_{12}P_{23}P_{13}$, $K_{13}K_{23}=\mu^{-1}{D_R}_2\; K_{13}P_{12}$
and $K_{12}K_{13}=\mu^{-1}{C_R}_3\; K_{12}P_{23}$).

\smallskip
An identity
\be K_{23}K_{14}P_{12}P_{34}=K_{23}K_{14}\ .\lb{kkpp}\ee
is again a consequence of the rank one property of the operator $K$.

\smallskip
Next, we have
\be K_{12}E_1 E_2= E_1 E_2K_{12}=K_{12}\ .\lb{kcax}\ee
To verify, for instance, that $K_{12}E_1^{-1} E_2^{-1}=K_{12}$, use the definition (\ref{defcaxy}) of the matrix
$E_2^{-1}$,  $E_2^{-1}=\tr_{(3)} (K_{23}P_{23})$, and then the relation (\ref{k12k23}) to remove the trace.

\medskip
The proof of the lemma is based on a following calculation
\be\begin{array}{l}
K_{23}\, E_2^m \, K_{12}\,  K_{34}\, E_2^{-m}\,  K_{23}=
K_{23}\, K_{34}\, E_2^m\,  K_{12}\, E_1^m\,  K_{23}
\\[1em]\ \ \
=E_4\,  K_{23}\, P_{23}\, P_{24}\, E_2^m \, E_3\,  K_{12}\, P_{12}\,
P_{13}\, E_1^m
= E_4^{m+1}\, K_{23}\,  K_{14}\, P_{23}\, P_{24}\, P_{12}\, P_{13}\, E_1^{m+1}
\\[1em]\ \ \
=E_4^{m+1}\,  K_{23}\, K_{14}\, P_{12}\, P_{34}\, E_1^{m+1}
=E_4^{m+1}\, K_{23}\,  K_{14}\, E_1^{m+1}
=K_{23}\, \Bigl( E_1^{-(m+1)}\, K_{14}\, E_1^{m+1}\Bigr)\ .
\end{array}\lb{kvtor}\ee
Here in the first equality we moved $K_{34}$ leftwards to $K_{23}$ and replaced $K_{12}\, E_2^{-m}$ by
$K_{12}\, E_1^m$ by eq.(\ref{kcax}); in the second equality we moved $E_1^m$ to the right
and used eq.(\ref{k12k23}) to rewrite products $K_{23}K_{34}$ and $K_{12}K_{23}$; in the third equality
we collected together all the permutation matrices as well as the powers of $E_4$ and
$E_1$; in the fourth equality we replaced $P_{23}P_{24}P_{12}P_{13}$ by $P_{12}P_{34}$;
in the fifth equality we used eq.(\ref{kkpp}); in the sixth equality we moved $E_{4}^{m+1}$
to $K_{14}$ and replaced $E_{4}^{m+1}K_{14}$ by $E_{1}^{-(m+1)}K_{14}$.

\medskip
Now we prove eq.(\ref{efoi}) by induction on $j$ (the case when $j=1$ is just the definition
of the element $c^{(2)}$).

Denote by $U[j,s]$ the $s$-th factor in the product in the left hand side of eq.(\ref{efoi}),
$U[j,s]=\eta^{-1}\; E_s^{s-j}\ K_{s,2j+1-s}\ E_s^{j-s}$.
Assume that $\rho_R(c^{(2j)})=\prod_{s=1}^{j}U[j,s]$ holds. By the recursive definition
(\ref{kappa-i}),
\be\rho_R(c^{(2j+2)})=\Bigl( \prod_{s=1}^{j}(U[j,s])^{\uparrow 1}\Bigr)\;
K_{12}K_{2j+1,2j+2}\;\Bigl( \prod_{s=1}^{j}(U[j,s])^{\uparrow 1}\Bigr)\ .\lb{reckne}\ee
Each factor $(U[j,s])^{\uparrow 1}$ is a projector, the factors $(U[j,s])^{\uparrow 1}$ pairwise commute
and there is only one factor, namely $(U[j,1])^{\uparrow 1}$, in each product in the right hand side of
eq.(\ref{reckne}), which does not commute with the operator $K_{12}K_{2j+1,2j+2}$. We have
$U[j+1,s]=(U[j,s])^{\uparrow 1}$ thus the expression (\ref{reckne}) reduces to
\be\rho_R(c^{(2j+2)})=\Bigl( \prod_{s=3}^{j+1}U[j+1,s]\Bigr)\;
\underline{U[j+1,2]\; K_{12}K_{2j+1,2j+2}\; U[j+1,2]}\ .\lb{reckne2}\ee
Now use the calculation (\ref{kvtor}) (with the needed relabeling) to rewrite the underlined expression
$X[j]=U[j+1,2]\; K_{12}K_{2j+1,2j+2}\; U[j+1,2]$ in the right hand side of eq.(\ref{reckne2}):
\be\begin{array}{l} X[j]=
\eta^{-2} E_{2}^{1-j}\Bigl( K_{2,2j+1}\; E_{2}^{j-1}\;  K_{12}\; K_{2j+1,2j+2}\;
E_{2}^{1-j}\; K_{2,2j+1}\Bigr)\; E_{2}^{j-1}\\[1em]\ \ \ \ \ \ \ \ \ \
=\eta^{-2} E_{2}^{1-j}\; K_{2,2j+1}\Bigl( E_{1}^{-j}\; K_{1,2j+2}\; E_{1}^{j}\Bigr)\;
E_{2}^{j-1}=U[j+1,1]\; U[j+1,2]\end{array}\ee
and to finish the proof.
\hfill$\blacksquare$

\begin{rem}{\rm The operator $K$ is skew invertible; we leave, as an exercise for the reader, to verify
that ${\Psi_K}_{12}=E_1 K_{12} E_1^{-1}=\mu^{-2}{C_R}_2 K_{21}{C_R}_2$.}\end{rem}

\section{Quantum matrix algebra}\lb{sec4}

In this section we deal with the main objects of our study, the quantum matrix algebras,
and construct the $\star$-product for them. We mainly discuss the quantum matrix algebras of
the type BMW (and of the subtypes $O$ and $Sp$).

\subsection{Definition}\lb{subsec4.1}

Consider a linear space ${\rm Mat}_{\mbox{\footnotesize\sc n}}(W)$, introduced in the definition
\ref{definition3.1}. For a fixed element $F\in {\rm Aut}(V\otimes V)$, we consider series of "copies"
$M_{\overline{i}}$, $i=1,2,\dots ,n,$ of a matrix $M\in {\rm Mat}_{\mbox{\footnotesize\sc n}}(W)$.
They are defined recursively by
\be M_{\overline 1}:=M_1, \quad M_{\overline{i}}:= F^{\phantom{-1}}_{i-1}M_{\overline{i-1}}F_{i-1}^{-1}\ .
\lb{kopii}\ee
{}For $F=P$, these are usual copies, $M_{\overline{i}}=M_i$, but, in general,
$M_{\overline{i}}$ can be nontrivial in all the spaces $1,\dots ,i$.

\smallskip
We shall, slightly abusing notation, denote by the same symbol $M_{\overline{i}}$
an element in ${\rm Mat}_{\mbox{\footnotesize\sc n}}(W)^{\otimes k}$ for any $k\geq i$, which is
defined by an inclusion of the spaces
\be\lb{kopii2} {\rm Mat}_{\mbox{\footnotesize\sc n}}(W)^{\otimes j}\hookrightarrow
{\rm Mat}_{\mbox{\footnotesize\sc n}}(W)^{\otimes (j+1)}:\quad
 {\rm Mat}_{\mbox{\footnotesize\sc n}}(W)^{\otimes j}\ni X\mapsto
X\otimes I\in{\rm Mat}_{\mbox{\footnotesize\sc n}}(W)^{\otimes (j+1)}\, .\ee

{}From now on we specify $W = {\Bbb C}\langle 1,M_a^b\rangle$, $1\le a,b\le \mbox{\sc n}$, ---
an associative $\Bbb C$-algebra freely generated by the unity and by $\mbox{\sc n}^2$ elements $M_a^b$.

\begin{defin}\lb{definition4.1}
Let $\{R,F\}$ be a compatible pair of strict skew invertible R-matrices (see the section \ref{subsec3.1}).
A quantum matrix algebra ${\cal M}(R,F)$ is a quotient algebra of the algebra ${\Bbb C}\langle M_i^j\rangle$ by a
two-sided ideal generated by entries of the matrix relation
\be R_1M_{\overline 1}M_{\overline 2} = M_{\overline 1}M_{\overline 2}R_1\, ,\label{qma}\ee
where $M = \|M_a^b\|_{a,b=1}^{\mbox{\footnotesize\sc n}}$ is a matrix of the generators
of ${\cal M}(R,F)$ and the matrix copies $M_{\overline i}$ are constructed with the help
of the R-matrix $F$ as in eq.(\ref{kopii}).

\medskip
If $R$ is a BMW-type R-matrix (see (\ref{charR})--(\ref{bmwR})$\,$) then ${\cal M}(R,F)$ is
called a BMW-type quantum matrix algebra.

\smallskip
If $R$ is a $Sp(2k)$- , respectively, an $O(k)$-type R-matrix (see the definition \ref{definition3.11}) then
${\cal M}(R,F)$ is called a $Sp(2k)$- , respectively, an $O(k)$-type quantum matrix algebra.\end{defin}

\begin{rem}\lb{remark4.2}
{\rm The quantum matrix algebras were introduced in Ref. \cite{Hl} under the name "quantized
braided groups". In the context of the present paper they have been first investigated
in \cite{IOP1}. The matrix $M'$ of the generators of the algebra ${\cal M}(R,F)$ used in \cite{IOP1}
is different from the matrix $M$ that we use here. A relation between these two matrices is explained
in the section 3 of \cite{IOP2}: $M' =D_R M (D_F)^{-1}$.}\end{rem}

\begin{lem}{\rm \bf \cite{IOP1}} \lb{lemma4.3}
The matrix copies of the matrix $M= \|M_a^b\|_{a,b=1}^{\mbox{\footnotesize\sc n}}$ of the generators of
the algebra ${\cal M}(R,F)$ satisfy relations
\ba F_i\, M_{\overline j} &=&\,M_{\overline j}\, F_i\,\;\;\quad\qquad \mbox{for}\;\; j\neq i,i+1,
\lb{fm-k}\\[1em] R_i\, M_{\overline j} &=&\,M_{\overline j}\, R_i\,\;\;\quad\qquad \mbox{for}\;\;
j\neq i,i+1,\lb{rm-k}\\[1em] R_j\,M_{\overline j}\,M_{\overline {j+1}} &=& M_{\overline j}\,
M_{\overline {j+1}}\, R_j\,\quad \mbox{for}\;\; j=1,2,\dots \ .\lb{rmm-k}\ea\end{lem}

\subsection{Characteristic subalgebra}\lb{subsec4.2}

{}From now on we assume that $M$ is the matrix of generators of the quantum matrix algebra
${\cal M}(R,F)$ and its copies $M_{\overline n}$ are calculated by the rule (\ref{kopii}).

\medskip
Denote by ${\cal C}(R,F)$ a vector subspace of the quantum matrix algebra ${\cal M}(R,F)$ linearly
spanned by the unity and elements
\be\lb{char} ch(\alpha^{(n)}) := \Tr{1,\dots ,n}(M_{\overline 1}\dots M_{\overline n}\,
\rho_R(\alpha^{(n)}))\ ,\quad n =1,2,\dots\ ,\ee
where $\alpha^{(n)}$ is an arbitrary element of the group algebra ${\Bbb C}{\cal B}_{n}$ of the braid
group ${\cal B}_n$ (we put, by definition, ${\cal B}_1 :=\{1\}$).

\smallskip
Notice that elements of the space ${\cal C}(R,F)$ satisfy a {\em cyclic property}
\be\lb{cyclic} ch(\alpha^{(n)}\beta^{(n)}) = ch(\beta^{(n)}\alpha^{(n)})\,  \quad
\forall\; \alpha^{(n)}, \beta^{(n)}\in {\Bbb C}{\cal B}_{n}\, ,\quad n=1,2,\dots\ ,\ee
which is a direct consequence of the relations (\ref{rm-k}), (\ref{rmm-k}) and (\ref{RDD}) and the cyclic property of
the trace.

\begin{def-prop} {\rm\bf \cite{IOP1}}\lb{proposition4.4}
The space ${\cal C}(R,F)$ is a commutative subalgebra of the quantum matrix algebra ${\cal M}(R,F)$,
the multiplication rule is
\be\lb{multip-rule} ch(\alpha^{(n)})\, ch(\beta^{(i)}) =ch(\alpha^{(n)}\, \beta^{(i)\uparrow n}) =
ch(\alpha^{(n)\uparrow i}\, \beta^{(i)})\, .\ee
Recall that $\alpha^{(n)\uparrow i}$ denotes the image of an element $\alpha^{(n)}$ under the embedding
${\cal B}_n \hookrightarrow {\cal B}_{n+i}$ defined in (\ref{h-emb2}).
We shall call ${\cal C}(R,F)$ the characteristic subalgebra of ${\cal M}(R,F)$.\end{def-prop}

A proof of the proposition given in \cite{IOP1} is based in particular on a following lemma:

\begin{lem}{\rm\bf \cite{IOP1}}\lb{lemma4.5}
Consider an arbitrary element $\alpha^{(n)}$ of the braid group ${\cal B}_n$. Let $\{R,F\}$
be a compatible pair of R-matrices, where $R$ is skew invertible. Then relations
\be\lb{char1}\Tr{i+1,\dots ,i+n}(M_{\overline{i+1}}\dots
M_{\overline{i+n}}\ \rho_R(\alpha^{(n)\uparrow i}))\, =\, I_{1,2,\dots ,i}\  ch(\alpha^{(n)})\,\ee
hold for any matrix $M\in {\rm Mat}_{\mbox{\footnotesize\sc n}}(W)$\footnote{
Here there is no need to specify $M$ to be the matrix of the generators of the algebra ${\cal M}(R,F)$.}.\end{lem}

We will make use of the lemma \ref{lemma4.5} several times below.
\medskip

Let us introduce a shorthand notation for certain elements of ${\cal C}(R,F)$
\ba\lb{P-01} p_0& :=& \tr_{\!\! R}\, I\; (=\mu \eta \mbox{~in the BMW case})\,  ,\qquad
p_1 := \tr_{\!\! R}\, M\, ,\\[1em]\lb{P_k} p_i& := & ch(\sigma_{i-1}\dots\sigma_2\sigma_1) =
ch(\sigma_1\sigma_2\dots\sigma_{i-1})\, ,  \quad i=2,3,\dots .\ea
(the last equality in eq.(\ref{P_k}) is due to the inner automorphism (\ref{innalis}) and the
cyclic property (\ref{cyclic})$\,$).

The elements $p_i$ are called {\em traces of powers of $M$} or, shortly, {\em power sums}.

\medskip
{}From now on in this subsection we assume the R-matrix $R$ and, hence, the algebra  ${\cal M}(R,F)$
to be of the BMW-type. Denote
\be\lb{tau}\textstyle g\, :=\, ch(c^{(2)}) \,\equiv\, \eta^{-1} ch(\kappa_1)\, \equiv\,  \eta^{-1}\,
\Tr{1,2} \left( M_{\overline{1}}M_{\overline{2}}\, K_1\right)\, .\ee
The notation used here was introduced in eqs.(\ref{kappa}), (\ref{bmw5}), (\ref{kappa-i}) and
(\ref{K}). We call $g$  a {\em contraction of two matrices} $M$ or, simply, a {\em 2-contraction}.

\begin{lem}\lb{lemma4.6}
Let $M$ be the matrix of generators of the BMW-type quantum matrix algebra ${\cal M}(R,F)$. Then
its copies, defined in eq.(\ref{kopii}), fulfill relations
\be\lb{tau2}K_{n}\, M_{\overline{n}}M_{\overline{n+1}}\, =\,
M_{\overline{n}}M_{\overline{n+1}}\, K_{n}\, =\,\mu^{-2} K_{n}\, g\, \quad
\forall\; n\geq 1\, .\ee\end{lem}

\nin {\bf Proof.~} We employ induction on $n$. Due to the property ${\rm rk}\, K=1$, one has
\be\lb{tau3} K_1\, M_{\overline{1}}M_{\overline{2}}\, =\, M_{\overline{1}}M_{\overline{2}}\, K_1\, =\,
K_1\, t\, ,\ee
where $t\in {\cal M}(R,F)$ is a scalar. Evaluating the R-trace of eq.(\ref{tau3}) in the spaces 1 and 2
and using eqs.(\ref{traceDK}) and (\ref{traceD}), one finds  $t=\mu^{-2} g$, which
proves the relation (\ref{tau2}) in the case $i=1$. It remains to check an induction step
$n\rightarrow (n+1)$:
$$\begin{array}{ccl} K_{n+1}M_{\overline{n+1}}M_{\overline{n+2}}&=& K_{n+1}(F_n M_{\overline{n}}
{}F^{-1}_{n})M_{\overline{n+2}}\, =\, K_{n+1} F_n M_{\overline{n}} (F_{n+1}M_{\overline{n+1}}F^{-1}_{n+1})
{}F^{-1}_{n} \\[1em] &=&(K_{n+1}F_n F_{n+1}) M_{\overline{n}} M_{\overline{n+1}}F^{-1}_{n+1}F^{-1}_{n}\, =\,
{}F_n F_{n+1} (K_{n}M_{\overline{n}} M_{\overline{n+1}})F^{-1}_{n+1}F^{-1}_{n}\\[1em] &=& \mu^{-2} F_n
{}F_{n+1} K_n F^{-1}_{n+1}F^{-1}_{n} g\, =\, \mu^{-2} K_{n+1}\, g\, .\end{array}$$
Here eqs.(\ref{kopii}) and (\ref{fm-k}), the twist relation (\ref{sovm}) for the pair $\{K,F\}$
and the induction assumption were used for the transformation. \hfill$\blacksquare$

\begin{prop}\lb{proposition4.7}
Let ${\cal M}(R,F)$ be the quantum matrix algebra of the BMW-type. Its characteristic subalgebra
${\cal C}(R,F)$ is generated by the set $\{g,p_i\}_{i\geq 0}$.\end{prop}

\nin {\bf Proof.~} Consider the chain of the BMW algebras monomorphisms (\ref{h-emb})--(\ref{h-emb2}).
We adapt, for $n\geq 3$, a following presentation for  an element $\alpha^{(n)}\in {\cal W}_n$
\be\lb{char8}\alpha^{(n)} = \beta\sigma_1 \beta' + \gamma \kappa_1 \gamma' + \delta\, ,\ee
where $\beta,\beta',\gamma,\gamma',\delta\in {\rm Im}({\cal W}_{n-1})\subset {\cal W}_n$. For $n=3$,
the formula (\ref{char8}) follows from the relations (\ref{braid})--(\ref{bmw7}). For  $n>3$, it can be proved
by induction on $n$.

\smallskip
Using the expression (\ref{char8}) for $\alpha^{(n)}$ and the cyclic property (\ref{cyclic}),
we conclude that, in the BMW case, any element (\ref{char}) of the characteristic subalgebra
can be expressed as a linear combination of terms
\be\lb{char8a} ch(\alpha_1\alpha_2\dots\alpha_{n-1})\, , \quad
\mbox{where}\quad \alpha_i\in\{1,\sigma_i,\kappa_i\}\, .\ee

Let us analyze the expressions (\ref{char8a}) for different choices of $\alpha_i$.

\medskip\noindent
i)~ If $\alpha_i=1$ for some value of $i$,  then, applying the relation (\ref{char1}), we get
\be\lb{char8b} ch(\alpha_1\dots\alpha_{i-1}\alpha_{i+1}\dots\alpha_{n-1}) =
ch(\alpha_1\dots\alpha_{i-1})\, ch((\alpha_{i+1}\dots\alpha_{n-1})^{\downarrow i})\, ,\ee
where $(\alpha_{i+1}\dots\alpha_{n-1})^{\downarrow i}\in {\cal W}_{n-i}$ is the preimage
of $(\alpha_{i+1}\dots\alpha_{n-1})\in {\cal W}_n$.

\medskip\noindent
ii)~ In the case when $\alpha_{n-1}=\kappa_{n-1}$, we apply eq.(\ref{tau2}) and then
eqs.(\ref{traceR}), (\ref{traceDK}) or (\ref{traceD}) to reduce the expression (\ref{char8a}) to
\be\lb{char8c} ch(\alpha_1\dots\alpha_{n-2}\kappa_{n-1}) =f(\alpha_{n-2})\,
ch(\alpha_1\dots\alpha_{n-3})\, g\, ,\ee
where $f(\sigma_{n-2})=\mu^{-1}$, $f(\kappa_{n-2})=1$ and $f(1)=\eta$.

\medskip\noindent
iii)~ In the case when $\alpha_i=\kappa_i$ for some $i$, and $\alpha_j=\sigma_j$ for all $j=i+1,\dots ,n-1$,
we perform following transformations
\be\begin{array}{ccl} ch(\alpha_1\dots\alpha_{i-1}\underline{\kappa_i\sigma_{i+1}}\sigma_{i+2}\dots \sigma_{n-1}) =
ch(\alpha_1\dots\alpha_{i-2}\,\sigma_i^{-1}\alpha_{i-1}\kappa_i\underline{\kappa_{i+1}\sigma_{i+2}}\dots
\sigma_{n-1}) &&\\[1em] =\dots =ch(\alpha_1\dots\alpha_{i-2} (\sigma_{n-2}^{-1}\dots \sigma_i^{-1})
\alpha_{i-1}\kappa_i\kappa_{i+1}\dots\kappa_{n-1}) .\hspace{35mm}&&\end{array}\lb{char10}\ee
Here the relations (\ref{bmw3}) and the cyclic property (\ref{cyclic}) are repeatedly used; expressions
suffering a transformation are underlined.

\smallskip
Now, depending on a value of $\alpha_{i-1}$, we proceed in different ways.
If $\alpha_{i-1}=\kappa_{i-1}$ then by  eqs.(\ref{bmw3}) and (\ref{char8c}) we have
\be\begin{array}{ccl} (\ref{char10}) &=& ch(\alpha_1\dots\alpha_{i-2}\,
\sigma_{i-1} \sigma_i\dots \sigma_{n-3} \kappa_{n-2}\kappa_{n-1})\\[1em] &=&
ch(\alpha_1\dots\alpha_{i-2}\,\sigma_{i-1} \sigma_i\dots \sigma_{n-3})\, g\, .\end{array}\lb{char11}\ee
If $\alpha_{i-1}=\sigma_{i-1} = \sigma_{i-1}^{-1}+(q-q^{-1})(1-\kappa_{i-1})$ then, using
the relations $\sigma_i^{-1}\sigma_{i-1}^{-1}\kappa_i=\kappa_{i-1}\kappa_i$ and applying
the previous results (\ref{char8c}) and (\ref{char8b}), we obtain
\be\begin{array}{ccl} (\ref{char10}) &=&ch(\alpha_1\dots\alpha_{i-2}\, \kappa_{i-1}\sigma_i\dots \sigma_{n-3})\, g
\\[1em]\nonumber &&+\, (q-q^{-1})\, \mu^{-1}\, ch(\alpha_1\dots\alpha_{i-2})\, p_{n-i-1}\, g\\[1em]
&&-\, (q-q^{-1})\, ch(\alpha_1\dots\alpha_{i-2}\,
\sigma_{i-1} \sigma_i\dots \sigma_{n-3})\, g\, .\end{array}\lb{char12}\ee

Iterating transformations i)---iii)  finitely many times, we eventually prove the assertion
of the proposition. \hfill$\blacksquare$

\medskip
We keep considering the BMW-type quantum matrix algebra ${\cal M}(R,F)$ with the R-matrix $R$ generating
representations of the algebras ${\cal W}_n(q,\mu)$, $n=1,2,\dots$. Assume that the antisymmetrizers $a^{(i)}$
and symmetrizers $s^{(i)}$ in these latter algebras are consistently defined (see eqs.(\ref{a^k}),
(\ref{s^k}) and (\ref{mu})$\,$). In this case, we can introduce two following sets of elements in the characteristic
subalgebra ${\cal C}(R,F)$
\ba\lb{SA_0}a_0& :=&1\ \ \ {\mathrm{and}}\ \ \ s_0\, :=\, 1\, ;\\[1em]\lb{SA_k} a_i &:=& ch(a^{(i)})\ \ \
{\mathrm{and}}\ \ \ s_i\, :=\,
ch(s^{(i)})\, ,\quad i=1,2,\dots  .\ea

\begin{prop}\lb{proposition4.8}
Let ${\cal M}(R,F)$ be the quantum matrix algebra of the BMW-type. Assume that $j_q\neq 0 ,\;\; \mu\neq -q^{-2j+3}$ (respectively,
$j_q\neq 0,$ $\mu\neq q^{2j-3}$) for all $j=2,3,\dots\; .$ Then the characteristic subalgebra
${\cal C}(R,F)$ is generated by the set $\{g,a_i\}_{i\geq 0}$ (respectively, $\{g,s_i\}_{i\geq 0}$).

\smallskip
Let the quantum matrix algebra ${\cal M}(R,F)$ be of the types either  $Sp(2k)$ or $O(k)$ (this implies restrictions on $q$:
$j_q\neq 0$ for all $j=2,3,\dots ,k$). Then the characteristic subalgebra ${\cal C}(R,F)$ is
generated by the set $\{g,a_i\}_{i=0}^k$.\end{prop}

\nin {\bf Proof.~} These statements are byproducts of the previous proposition and the Newton relations,
which are proved in the section 6, the theorem \ref{theorem6.1}
(the proof of
the theorem relies on the basic identities from the lemma \ref{lemma5.1}). \hfill$\blacksquare$

\subsection{Reciprocal relations in ${\cal C}(R,F)$ in orthogonal case}\lb{subsec4.3}

{}For the orthogonal $O(k)$-type quantum matrix algebras, the rank-one condition on the operator
$\rho_R(a^{(k)})\in {\rm End}(V^{\otimes k})$ gives rise to algebraic dependencies among
$a_i$ and the 2-contraction $g$. In this subsection we construct explicitly
a family of equalities --- we call them {\em reciprocal relations} --- for the set of generators
of the characteristic subalgebra.

\begin{theor}\lb{theorem4.9}
Let the quantum matrix algebra ${\cal M}(R,F)$ be of the orthogonal $O(k)$-type. Then
following reciprocal relations for the generators of its
characteristic subalgebra $\{g,a_i\}_{i=0}^k$ (see the proposition \ref {proposition4.8})
\be\lb{reciprocal} g^{k-i}\, a_i\, =\, a_k\, a_{k-i}\,\qquad \forall\; i=0,1,\dots ,k\, \ee
are satisfied.\end{theor}

\nin {\bf Proof.~} For the proof, we use two auxiliary statements.

\begin{lem}\lb{lemma4.10}
Let the quantum matrix algebra ${\cal M}(R,F)$ be of the BMW-type, respectively, of the $O(k)$-type.
Then we have, respectively,
\ba\lb{spec-cM} && M_{\overline{1}}\, M_{\overline{2}}\dots M_{\overline{2i}}\,\,
\rho_R(c^{(2i)})\, =\, (\mu^{-2} g)^i\, \rho_R(c^{(2i)})\, \ea
in the BMW case and
\ba\lb{spec-aM} && M_{\overline{1}}\,M_{\overline{2}}\dots M_{\overline{k}}\,\,
\rho_R(a^{(k)})\, =\, q^{k(k-1)}\, a_k\, \rho_R(a^{(k)})\, \ea
in the $O(k)$ case.\end{lem}

\nin {\bf Proof.~} The relations (\ref{spec-cM}) and (\ref{spec-aM}) are proved in exactly the same way as
eqs.(\ref{spec-c1}) and (\ref{spec-a1}) of the lemma \ref{lemma3.14}.
One has to use additionally a relation (\ref{tau2}) for the proof of eq.(\ref{spec-cM}).
{}For eq.(\ref{spec-aM}), an equality
$$\tr_{(1,2,\dots k)} (M_{\overline{1}}\,M_{\overline{2}}\dots M_{\overline{k}}
\,\rho_R(a^{(k)}))\, =\, q^{k(k-1)}\,\Tr{1,2,\dots k}
(M_{\overline{1}}\,M_{\overline{2}}\dots M_{\overline{k}} \,\rho_R(a^{(k)}))\, =\, q^{k(k-1)}\, a_k\ ,$$
following from eq.(\ref{spec-a1}), also should be taken into account. \hfill$\blacksquare$

\begin{lem}\lb{lemma4.11}
Let the operators $Z^{(i)}\in {\rm Aut}(V^{\otimes i})$ be defined as in the formulas (\ref{Z-i}). For a given matrix
$M\in {\rm Mat}_{\mbox{\footnotesize\sc n}}(W)$, consider one more set of its "copies" defined as (c.f.
with eq.(\ref{kopii})$\,$)
\be\lb{kopii-2} M_{\underline{1}}\, :=\, M_1 ,\qquad M_{\underline{i+1}}\, :=\,
{}F_i^{-1}M_{\underline{i}}\, F_i , \qquad i=2,3,\dots .\ee
Then relations
\be \lb{Z-M}M_{\overline{j}}\, Z^{(i)}\, =\, Z^{(i)}\, M_{\underline{i-j+1}}\ee
hold for any $i=1,2,\dots$ and for any $j=1,2,\dots ,i$.\end{lem}

\nin {\bf Proof.~} The equalities (\ref{Z-M}) can be proved by induction. First, one considers the case
$j=i$ and carries out induction on $i$. Then, performing induction on $(i-j)$, one proves
eqs.(\ref{Z-M}) in a full generality. \phantom{a}\hfill$\blacksquare$

\medskip
Now, we will prove the relations (\ref{reciprocal}) by expressing an element $ch(c^{(2k-2i)}a^{(k)\uparrow k-i})
\in {\cal C}(R,F)$ in terms of elementary generators in two different ways.

\smallskip
The first way is
$$\begin{array}{l} ch(\lefteqn{ c^{(2k-2i)} a^{(k)\uparrow k-i}) = (\mu^{-2} g)^{k-i}\, \Tr{1,2,\dots ,2k-i}
\left( M_{\overline{2k-2i+1}}\dots M_{\overline{2k-i-1}} M_{\overline{2k-i}}\,\rho_R
(c^{(2k-2i)}a^{(k)\uparrow k-i})\right) }\\[1em] =(\mu^{-2} g)^{k-i}\,\Tr{1,2,\dots ,2k-i}
\left(Z^{(2k-i)}M_{\underline{i}}\dots M_{\underline{2}} M_{\underline{1}}\,\,\rho_{R_f}(c^{(2k-2i)
\uparrow i} a^{(k)}) (Z^{(2k-i)})^{-1}\right) \hspace{20mm}\\[1em] =(\mu^{-2} g)^{k-i}\,
\Tr{1,2,\dots ,i} M_{\underline{i}}\dots M_{\underline{2}} M_{\underline{1}}\,\Bigl( \TR{i+1,\dots ,2k-i}
{R_f}\, \rho_{R_f}(c^{(2k-2i)\uparrow i} a^{(k)})\Bigr) \\[1em] =(\eta^{-1}\mu^{-1} g)^{k-i}\,
\Tr{1,2,\dots ,i} M_{\underline{i}}\dots M_{\underline{2}} M_{\underline{1}}\, \Bigl( \TR{i+1,\dots ,k}
{R_f}\,\rho_{R_f}(a^{(k)})\Bigr) \\[1em] ={\displaystyle \bigl(\prod_{j=1}^{k-i}\,\eta^{-1}\mu^{-1}
\delta_{k-j+1}\, g\bigr)}\,\Tr{1,2,\dots ,i}\left( (Z^{(i)})^{-1} M_{\overline{1}}M_{\overline{2}}\dots
M_{\overline{i}}\, \rho_R(a^{(i)}) Z^{(i)}\right) \\[1.5em]={\displaystyle\bigl(\prod_{j=1}^{k-i}\,
\eta^{-1}\mu^{-1}\delta_{k-j+1}\bigr)}\, g^{k-i}\, a_i\, .\end{array}$$
Here eq.(\ref{spec-cM}) was applied  in the first line; eqs.(\ref{Z-M}) and
direct generalizations
\ba
\nonumber
\rho_R(c^{(2i)}) Z^{(j)}& =& Z^{(j)} \rho_{R_f}(c^{(2i)\uparrow j-2i})\qquad \forall\; j\geq 2i\geq 0\, ,
\\
\nonumber
\rho_R(a^{(i)\uparrow j-i}) Z^{(j)}& =&  Z^{(j)}\rho_{R_f}(a^{(i)})\qquad\qquad\; \forall\; j\geq i\geq 0
\ea
of eqs.(\ref{Z-ac}) were used in the second line; the commutativity of $\prod_{j=1}^{2k-i}{D_R}_j$ with $Z^{(2k-i)}$
allows to cancel the latter when passing to third line; eqs.(\ref{spec-c1}) and (\ref{spec-c2}) were used
to substitute the R-trace $\tr_{\!\!\raisebox{-3pt}{$\scriptstyle R$}}$ by the R-trace
$\tr_{\!\!\raisebox{-1.5pt}{$\scriptstyle R_f$}}$ in the spaces with labels ~$\scriptstyle i+1,i+2,\dots ,
2k-i$~ in the third line; eqs.(\ref{traces-c2i}) and (\ref{spec1}) were applied for the evaluation of the R-traces
in the fourth and fifth lines; a similarity transformation by $Z^{(i)}$ was performed in the argument of the
R-trace and, then, the relations (\ref{Z-M}) and (\ref{Z-ac}) were used to complete the transformation.

\smallskip
The second way is
$$\begin{array}{l} ch(\lefteqn{ c^{(2k-2i)} a^{(k)\uparrow k-i}) =
q^{k(k-1)}\, \Tr{1,2,\dots ,2k-i}\Bigl(M_{\overline{1}} M_{\overline{2}}\dots M_{\overline{k-i}}\,
\rho_R(a^{(k)\uparrow k-i} c^{(2k-2i)})\Bigr)\, a_k } \\[1em] =q^{k(k-1)}\,
\bigl({\displaystyle \prod_{j=k-i+1}^{k} \delta_j}\bigr)\, \Tr{1,2,\dots ,2k-2i}\Bigl( M_{\overline{1}}
M_{\overline{2}}\dots M_{\overline{k-i}}\,\rho_R(a^{(k-i)} c^{(2k-2i)})\Bigr)\, a_k \hspace{10mm}\\[1.5em]
=q^{k(k-1)}\, \bigl({\displaystyle \prod_{j=k-i+1}^{k} \delta_j}\bigr)\,
\bigl( \eta^{-1}\mu\bigr)^{k-i}\ \Tr{1,2,\dots ,k-i}\Bigl( M_{\overline{1}}
M_{\overline{2}}\dots M_{\overline{k-i}}\,\rho_R(a^{(k-i)})\Bigr) \, a_k \\[1em] =
\bigl({\displaystyle \prod_{j=1}^{k-i} \eta^{-1}\mu\, \delta_j^{-1}}\bigr)\,a_{k-i}\, a_k\, .\end{array}$$
Here (a version of) eq.(\ref{spec-aM}) was used in the first line; eq.(\ref{spec1}) was applied for the
evaluation of the R-traces in the second line; there also, using eqs.(\ref{idemp-c2}) and (\ref{a^k}), we
substituted $a^{(k-i)\uparrow k-i}\, c^{(2k-2i)}$ by $a^{(k-i)} c^{(2k-2i)}$ in the argument of $\rho_R$;
in the third line, some R-traces were evaluated with the help of eq.(\ref{traces-c2i}); the relations (\ref{Delta-i})
and (\ref{qdim-O(k)}) were used to transform the expression to its final form in the fourth line.

\smallskip
Comparing the results of these two calculations and substituting the values, which $\mu$ and $\delta_j$ take in
the $O(k)$ case (see the definition \ref{definition3.11} and eq.(\ref{spec1})$\,$), we verify eq.(\ref{reciprocal}).
\hfill$\blacksquare$

\subsection{Matrix $\star\,$-product}\label{subsec4.4}

Consider the quantum matrix algebra ${\cal M}(R,F)$ of the generic type (no additional conditions on an R-matrix $R$).

\smallskip
Denote by ${\cal P}(R,F)$ a linear subspace of ${\rm Mat}_{\mbox{\footnotesize\sc n}}({\cal M}(R,F))$
spanned by ${\cal C}(R,F)$-multiples of the identity matrix, $I\, ch$ $\forall\, ch\in {\cal C}(R,F)$,
and by elements
\be\lb{pow} M^1 := M , \qquad (M^{\alpha^{(n)}})_{1} := \Tr{2,\dots ,n}(M_{\overline 1}
\dots M_{\overline n}\,\rho_R(\alpha^{(n)}))\ ,\quad n =2,3,\dots\ ,\ee
where $\alpha^{(n)}$ belongs to the group algebra of the braid group ${\cal B}_n$. The space ${\cal P}(R,F)$
inherits a structure of a right ${\cal C}(R,F)$--module
\be\lb{r-module} M^{\alpha^{(n)}} ch(\beta^{(i)}) =M^{(\alpha^{(n)}\beta^{(i)\uparrow n})}\, \quad
\forall\, \alpha^{(n)}\in {\cal B}_n, \; \beta^{(i)}\in {\cal B}_i\, ,\quad n,i=1,2,\dots\, ,\ee
which is just a component-wise multiplication of the matrix $M^{\alpha^{(n)}}$ by $ch(\beta^{(i)})$
(use eq.(\ref{char1}) to check this). The ${\cal C}(R,F)$--module structure agrees with an R-trace map
$\tr_{\!\! R}$ (which means that $\tr_{\!\! R} (Xa)=\tr_{\!\! R} (X)a\ \ \forall\ X\in{\cal P}(R,F)$
and $\forall\ a\in {\cal C}(R,F)$)
\be\lb{Rtrace-map}{\cal P}(R,F)\stackrel{\tr_{\!\! R}}{\longrightarrow}
{\cal C}(R,F)\, : \quad\left\{ \begin{array}{rcl} M^{\alpha^{(n)}}&\mapsto& ch(\alpha^{(n)})\, ,
\\[1em] I\, ch(\alpha^{(n)})&\mapsto& (\tr_{\!\! R} I)\, ch(\alpha^{(n)})\, , \end{array}\right.\ee
where $\alpha^{(n)}\in {\cal B}_n\, , \;\; n=1,2,\dots$

\smallskip
Besides, elements of the space ${\cal P}(R,F)$ satisfy a {\em reduced cyclic property}
\be\lb{red-cycl} M^{(\alpha^{(n)}\beta^{(n-1)\uparrow 1})} =
M^{(\beta^{(n-1)\uparrow 1}\alpha^{(n)})}\, \quad\forall\, \alpha^{(n)}\in {\cal B}_n , \;
\beta^{(n-1)}\in {\cal B}_{n-1} ,\quad n=2,3,\dots \, .\ee

\begin{def-prop}\lb{proposition4.12}
{}Formulas
\be\lb{MaMb}M^{\alpha^{(n)}}\! \star  M^{\beta^{(i)}} \,
:=\,  M^{(\alpha^{(n)}\star \beta^{(i)})}\, ,\ee
where
\ba
\lb{a*b}
\qquad\alpha^{(n)}\star \beta^{(i)} &:=&
\alpha^{(n)}\beta^{(i)\uparrow n} (\sigma_n\dots \sigma_2 \sigma_1\sigma_2^{-1}\dots \sigma_n^{-1})\, ,
\\[1em]\lb{MaI} (I\, ch(\beta^{(i)})) \star  M^{\alpha^{(n)}} &:=&
M^{\alpha^{(n)}}\! \star  (I\, ch(\beta^{(i)}))\, :=\, M^{\alpha^{(n)}}  ch(\beta^{(i)})\, ,
\\[1em]
\lb{II}
(I\, ch(\alpha^{(i)}))\star (I\, ch(\beta^{(n)}))& :=& I\,( ch(\alpha^{(i)})\, ch(\beta^{(n)}))\ ,
\ea
define an associative multiplication on the space ${\cal P}(R,F)$, which agrees with the ${\cal C}(R,F)$--module
structure (\ref{r-module}).\footnote{In other words,
a map
~$ch(\alpha^{(n)})\mapsto I\, ch(\alpha^{(n)})$~ is an algebra monomorphism
~${\cal C}(R,F)\hookrightarrow {\cal P}(R,F)$.}\end{def-prop}

\nin {\bf Proof.~} To prove the associativity of the multiplication (\ref{MaMb}), it is enough to check
$$(\alpha^{(n)} \star  \beta^{(i)}) \star  \gamma^{(m)}\ =\alpha^{(n)} \star  (\beta^{(i)} \star
\gamma^{(m)})\, ,$$
which is a staightforward exercise in an application of eqs.(\ref{braid}) and (\ref{braid2}).
It is less trivial to prove a compatibility condition for eqs.(\ref{MaMb}) and (\ref{MaI})
$$\left\{M^{\alpha^{(n)}}\! \star  (I\, ch(\beta^{(i)}))\right\} \star  M^{\gamma^{(m)}}\, =\, M^{\alpha^{(n)}} \star
\left\{(I\, ch(\beta^{(i)}))\star  M^{\gamma^{(m)}}\right\}\, ,$$
which, in terms of the matrix "exponents", amounts to
\be\begin{array}{l} \alpha^{(n)}\beta^{(i)\uparrow n} \gamma^{(m)\uparrow (i+n)}
(\sigma_{i+n}\dots \sigma_2\sigma_1\sigma_2^{-1}\dots \sigma_{i+n}^{-1})\hspace{50mm}\\[1em]
\ \ \ \ \ \ \ \ \ \ \ \ \ \ \ \stackrel{\rm mod\, (\ref{red-cycl})}{=}\
\alpha^{(n)}\gamma^{(m)\uparrow n} \beta^{(i)\uparrow (m+n)}
(\sigma_n\dots \sigma_2\sigma_1\sigma_2^{-1}\dots \sigma_n^{-1})\, .\end{array}\lb{lesstriv}\ee
Here the symbol $\stackrel{\rm mod\, (\ref{red-cycl})}{=}$ means the equality modulo the reduced cyclic property
(\ref{red-cycl}).

\smallskip
To check eq.(\ref{lesstriv}), we apply a technique, which was used in \cite{IOP1} to prove
the commutativity of the characteristic subalgebra. Consider an element
\be\lb{u}
\begin{array}{rcl}
u_{i,m}^{(i+m)} &:=& (\sigma_i\dots \sigma_2\sigma_1)
(\sigma_{i+1}\dots \sigma_3\sigma_2) \dots (\sigma_{i+m-1}\dots \sigma_{m+1}\sigma_m)\\[1em]
&=&(\sigma_i\sigma_{i+1}\dots\sigma_{i+m-1})
(\sigma_{i-1}\sigma_i\dots \sigma_{i+m-2}) \dots (\sigma_{1}\sigma_2\dots\sigma_m)\, ,
\end{array}
\ee
which intertwines certain elements of the braid group ${\cal B}_{(i+m)}$:
\be\lb{u1} \beta^{(i)}\, u_{i,m}^{(i+m)}\, =\, u_{i,m}^{(i+m)}\, \beta^{(i)\uparrow m}\, ,\qquad
\gamma^{(m)\uparrow i}\, u_{i,m}^{(i+m)}\, =\, u_{i,m}^{(i+m)}\, \gamma^{(m)}\, .\ee
Substitute an expression $(u_{i,m}^{(i+m)\uparrow n}\gamma^{(m)\uparrow n}\beta^{(i)\uparrow (n+m)}
(u_{i,m}^{(i+m)\uparrow n})^{-1})$ for the factor $(\beta^{(i)\uparrow n} \gamma^{(m)\uparrow (i+n)})$
in the left hand side of the equation (\ref{lesstriv}), move the element $u_{i,m}^{(i+m)\uparrow n}$
cyclically to the right and then use an equality
\be\lb{u2} (\sigma^{-1}_1\sigma^{-1}_2\dots \sigma^{-1}_i) u_{i,m}^{(i+m)}
= u_{i,m-1}^{(i+m-1)\uparrow 1}\ee
to cancel it on the right hand side. Such transformation results in the right hand side of
the equation (\ref{lesstriv}).

\smallskip
Consistency of the multiplication and the ${\cal C}(R,F)$--module structures on ${\cal P}(R,F)$
follows  obviously from the last equality in
(\ref{MaI}). \hfill$\blacksquare$

\medskip
To illustrate the relation between the $\star \, $-product and the usual matrix multiplication,
we present the formulas (\ref{MaMb}) and (\ref{a*b}) in the case $n=1$ ($\alpha^{(1)}\equiv 1$) in a form
\be\lb{M*} M \star  N = M\cdot \phi(N)  \quad \forall N\in {\cal P}(R,F)\, ,\ee
where "$\cdot$" denotes the usual matrix multiplication and the map $\phi$ is defined by eq.(\ref{phi}).

\smallskip
The noncommutative analogue of the matrix power is given by a repeated $\star \, \, $-multi\-pli\-ca\-tion
by the matrix $M$
\be\lb{M^k} M^{\overline{0}} := I\, , \qquad M^{\overline{n}}\, :=\,
\underbrace{M\star  M\star \dots \star  M}_{\mbox{\small $n$ times}}\, =\, M^{(\sigma_1\sigma_2\dots \sigma_{n-1})}\, =\,
M^{(\sigma_{n-1}\dots \sigma_2\sigma_1)}\, .\ee
Here we introduce symbol $M^{\overline{n}}$ for the {\em $n$-th power of the matrix $M$}. The standard
matrix powers multiplication formula follows immediately from the definition
\be\lb{M^k* M^p}M^{\overline{n}} \star  M^{\overline{i}}\, =\, M^{\overline{n+i}}\, .\ee

\begin{prop}\lb{proposition4.13}
A ${\cal C}(R,F)$--module, generated by the matrix powers $M^{\overline{n}}$, $n=0,1,\dots$,
belongs to the center of the algebra ${\cal P}(R,F)$.\end{prop}

\nin {\bf Proof.~} It is sufficient to check a relation $M\star  M^{\alpha^{(i)}} = M^{\alpha^{(i)}}
\star  M$, which, in turn, follows from a calculation\\[3pt]
\phantom{a}\hfill
$ \hspace{15mm} \alpha^{(i)}\sigma_i\dots \sigma_2\sigma_1\sigma_2^{-1}\dots
\sigma_i^{-1}\, =\, \sigma_i\dots \sigma_2\sigma_1\alpha^{(i)\uparrow 1}\sigma_2^{-1}\dots
\sigma_i^{-1}\,\stackrel{\rm mod\, (\ref{red-cycl})}{=}\, \alpha^{(i)\uparrow 1} \sigma_1\, .$
\hfill$\blacksquare$

\bigskip
It is natural to expect that the algebra ${\cal P}(R,F)$ is commutative as all of its elements are
generated by the matrix $M$ alone. We can prove the commutativity in the BMW case. Notice that
(unlike in the Hecke case),
in the BMW case, the algebra ${\cal P}(R,F)$ can not be generated by the $\star\,$-powers of $M$ only.

\smallskip
By an analogy with the formula (\ref{M*}), we define a ${\cal C}(R,F)$--module map~
$\Mt : {\cal P}(R,F)$ $\rightarrow$ $ {\cal P}(R,F)$
\be\lb{Mt} \Mt (N) := M\cdot \xi(N), \qquad  N\in {\cal P}(R,F)\, ,\ee
where the endomorphism $\xi$ is defined by eq.(\ref{xi}). Equivalently, we can write
\be\lb{Mt-2} \Mt (M^{\alpha^{(n)}}) = M^{(\alpha^{(n)\uparrow 1}\kappa_1 )}
\qquad \forall\; \alpha^{(n)}\in {\cal W}_n ,\quad n=1,2,\dots \, .\ee

\begin{prop}\lb{proposition4.14}
Let the quantum matrix algebra ${\cal M}(R,F)$ be of the BMW-type. Then the algebra ${\cal P}(R,F)$
is commutative. As a ${\cal C}(R,F)$--module, it is spanned by matrices
\be\lb{P-gen} M^{\overline{n}}\,  \quad \mbox{and}\quad
\Mt(M^{\overline{n+2}})\, , \quad n=0,1,\dots\, .\ee\end{prop}

\nin{\bf Proof.~} A proof of the last statement of the proposition goes essentially along
the same lines as the
proof of the proposition \ref{proposition4.7} and we will not repeat it. The only modification is a reduction
of the cyclic property (c.f., eqs.(\ref{cyclic}) and (\ref{red-cycl})$\,$), which finally leads to
an appearance of the additional elements $\{\Mt(M^{\overline{n}})\}_{n\geq 2}$ in the generating set.

\smallskip
To prove the commutativity of ${\cal P}(R,F)$, we derive an alternative expression for the exponent in
the matrix product formula (\ref{MaMb})
\be\lb{a*b-2}\alpha^{(n)}\star \beta^{(i)}\, =\, (\sigma_i^{-1}\dots \sigma_2^{-1} \sigma_1
\sigma_2\dots \sigma_i)\alpha^{(n)\uparrow i}\beta^{(i)}\, .\ee
The calculation proceeds as follows
\ba\nonumber &&\hspace{-6mm}\alpha^{(n)}\star \beta^{(i)} =\alpha^{(n)}\beta^{(i)\uparrow n}
(\sigma_n\dots\sigma_1\sigma_2^{-1}\dots\sigma_n^{-1}) =
u_{n,i}^{(n+i)}\alpha^{(n)\uparrow i}\beta^{(i)}(u_{n,i}^{(n+i)})^{-1}
(\sigma_n\dots\sigma_1\sigma_2^{-1}\dots\sigma_n^{-1})\\[1em]\nonumber &&\hspace{-5mm}
\stackrel{\rm mod (\ref{red-cycl})}{=}\, (u_{n,i-1}^{(n+i-1)\uparrow 1})^{-1}(\sigma_2^{-1}\dots\sigma_n^{-1})
u_{n,i}^{(n+i)}\alpha^{(n)\uparrow i}\beta^{(i)} = (u_{n,i-1}^{(n+i-1)\uparrow 1})^{-1}\sigma_1
u_{n,i-1}^{(n+i-1)\uparrow 1}\alpha^{(n)\uparrow i}\beta^{(i)}\\[1em]\nonumber &&
= (\sigma_i^{-1}\dots\sigma_2^{-1}) (u_{n-1,i-1}^{(n+i-2)\uparrow 2})^{-1}\sigma_1
u_{n-1,i-1}^{(n+i-2)\uparrow 2}(\sigma_2\dots\sigma_i)\alpha^{(n)\uparrow i}\beta^{(i)} =
\mbox{right hand side of eq.(\ref{a*b-2})}.\ea
Here  we applied again the intertwining operators (\ref{u}) and used their properties (\ref{u1}) and
(\ref{u2}) and the reduced cyclicity. One more property
\be\lb{u3} u_{n,i}^{(n+i)} = u_{n-1,i}^{(n+i-1)\uparrow 1} (\sigma_1\sigma_2\dots \sigma_i)\ee
is used in the last line of the calculation.

\medskip
Due to the proposition \ref{proposition4.13}, to prove the commutativity of the algebra ${\cal P}(R,F)$, it
remains to check the commutativity of the set $\{\Mt(M^{\overline{n}})\}_{n\geq 2}$.

\smallskip
Notice that the
factors of the exponents of the matrices $\Mt(M^{\overline{n}})$ can be taken in an opposite order,
$\Mt(M^{\overline{n}}) = M^{(\kappa_1 \sigma_2 \sigma_3\dots\sigma_n)}
=M^{(\sigma_n\dots \sigma_3 \sigma_2\kappa_1)}\, .$
This observation, together with the formula (\ref{a*b-2}), allow us to choose the exponents of two matrices
$\Mt(M^{\overline{n}})\star \Mt(M^{\overline{i}})$ and $\Mt(M^{\overline{i}})\star \Mt(M^{\overline{n}})$
to be mirror (left-right) images of each other. Finally, the exponentiation of the matrix $M$ with mirror symmetric
arguments gives the identical results since they can be expanded into linear combinations of the generators
(\ref{P-gen}), which are invariant with respect to the mirror reflection of their exponents, and since
the expansion rules (that is, the defining relations for the BMW algebras) are mirror symmetric as well.\hfill$\blacksquare$

\begin{lem}\lb{lemma4.15}
{}For the BMW-type quantum matrix algebra ${\cal M}(R,F)$, one has
\ba\lb{Mt-IM}\Mt (I)\, =\, \mu M\, , \qquad \Mt (M)\, =\, \mu^{-1} I\, g\, ,\\[1em]\lb{Mt2-N}
\Mt (\Mt (N))\, =\, N\, g\,  \quad\forall\; N\in {\cal P}(R,F)\, .\ea\end{lem}

\nin {\bf Proof.~} The relations (\ref{Mt-IM}) follow immediately from eqs.(\ref{tau2}) and (\ref{traceDK}) and
the definitions (\ref{Mt}) and (\ref{xi}).

\smallskip
As for the equality (\ref{Mt2-N}), it is enough to check it in the case when the matrix $N$ is a power
of the matrix $M$.

\smallskip
To evaluate the expression
$\Mt(\Mt(M^{\overline n}))=M^{(\kappa_1\kappa_2\sigma_3\dots\sigma_{n+1})}$, we transform its exponent,
using the relations (\ref{bmw3}) in the BMW algebra and the reduced cyclic property, to
\be\begin{array}{ccl}\kappa_1\kappa_2\sigma_3\dots\sigma_{n+1}&=&\kappa_1(\kappa_2\kappa_3\sigma_2^{-1})
\sigma_4\dots\sigma_{n+1}\,\stackrel{\rm mod\, (\ref{red-cycl})}{=}\, (\sigma_2^{-1}\kappa_1\kappa_2)
\kappa_3\sigma_4\dots\sigma_{n+1}\\[1em] &=&\sigma_1\kappa_2\kappa_3\sigma_4\dots\sigma_{n+1}
\, =\; \dots\;\stackrel{\rm mod\, (\ref{red-cycl})}{=}\,\sigma_1\sigma_2\dots\sigma_{n-1}\kappa_n
\kappa_{n+1}\, .\end{array}\lb{trans-exp}\ee
{}For the exponent (\ref{trans-exp}), the matrix power is easily calculated, again with the help of
eqs.(\ref{tau2}) and (\ref{traceDK}), and gives the expression $M^{\overline{n}}g$. \hfill$\blacksquare$

\medskip
The last relation in (\ref{Mt-IM}) shows that to introduce the inverse matrix to the matrix $M$ it is
sufficient to add the inverse $g^{-1}$ of the 2-contraction $g$ to the algebra ${\cal M}(R,F)$.
This is realized in the next subsection.

\subsection{Matrix inversion}\label{subsec4.5}

In this subsection we define an {\em extended} quantum matrix algebra, to which the
inverse of the quantum matrix belongs. Then we solve the reciprocal relations from
the subsection \ref{subsec4.3}.

\begin{lem}\lb{lemma4.16}
Let ${\cal M}(R,F)$ be the BMW-type quantum matrix algebra. Its 2-contraction $g$ fulfills a relation
\be\lb{Mj} M\, g\, =\, g\, ( G^{-1} M G)\, ,\ee
where $G$ is defined by eq.(\ref{G}).\end{lem}

\nin {\bf Proof.~} The proof consists of a calculation
\be\begin{array}{ccl}\phantom{a}\hspace{-10mm} M_1\left(g K_2\right) &=&\mu^2 M_{\overline{1}}
M_{\overline{2}}M_{\overline{3}} K_2\, =\,\mu^2 M_{\overline{1}}M_{\overline{2}}M_{\overline{3}}K_2K_1K_2\,
=\,\mu^2 K_2 \left( M_{\overline{1}}M_{\overline{2}} K_1\right) M_{\overline{3}} K_2\\[1em] &=& gK_2K_1
M_{\overline{3}} K_2 \, =\, \left(g K_2\right) \tr_{(2,3)}\left(K_2 K_1 M_{\overline{3}}\right)
\\[1em] &=& \left(g K_2\right)\tr_{(2,3)}\left(K_2 F_2 F_1 K_2 M_1 F_1^{-1} F_2^{-1}\right)
\\[1em] &=& \left(g K_2\right)\tr_{(2,3)}(F_2 F_1 K_2) M_1
\tr_{(2,3)}(K_2 F_1^{-1} F_2^{-1}) =\left(g K_2\right) (G^{-1} M G)_1 .\lb{MjK}\end{array}\ee
Here the relations (\ref{tau2}) and (\ref{bmw5}) were used in the first two lines; the property
${\rm rk}\, K=1$ was
used in the last/first equality of the second/fourth line; the definition of $M_{\overline{3}}$ was
substituted and the twist relation for the pair $\{K,F\}$ was used in the third line; the formulas (\ref{G})
and (\ref{G-inv}) for $G$ and $G^{-1}$ were substituted in the last equality. \hfill$\blacksquare$

\begin{def-prop}\lb{proposition4.17}
Let ${\cal M}(R,F)$ be the BMW-type quantum matrix algebra. Consider an extension of ${\cal M}(R,F)$ by a
generator $g^{-1}$ subject to relations
\be\lb{j-inv} g^{-1}\, g\, =\, g\, g^{-1}\, =\, 1\, , \qquad g^{-1}\, M\, =\, (G^{-1}MG)\, g^{-1}\, .\ee
The extended algebra, which we shall further denote by ${\cal M^{^\bullet\!}}(R,F)$, contains
an inverse matrix to the matrix $M$
\be\lb{M-inv} M^{-1}\, :=\, \mu\, \xi(M)\, g^{-1}\, :\qquad
M\cdot M^{-1}\, =\, M^{-1}\cdot M\, =\, I\, .\ee\end{def-prop}

\nin {\bf Proof.~} The lemma \ref{lemma4.16} ensures the consistency of the relations (\ref{j-inv}). The equality
$M\cdot M^{-1}=I$ for the inverse matrix (\ref{M-inv}) follows immediately from eqs.(\ref{Mt-IM})
and (\ref{Mt}).
To prove the equality $M^{-1}\cdot M = I$, consider a mirror partner of the map $\xi$:
\be\lb{Pi}\theta(M)\, :=\, \mu^{-2}\, \Tr{2} K_1 M_{\overline{2}}\, .\ee
By the (left-right) symmetry arguments in the assumptions of the lemma \ref{lemma3.9},
the map $\theta$ is invertible and the inverse map reads
\be\lb{Pi-inv}\theta^{-1}(M)\, =\,\TR{2}{R_f}\left( F^{-1}_1 K_1 M_1  F_1\right)\, .\ee
Applying in a standard way the transformation formula (\ref{inv-trD}), we calculate a composition of
the maps $\xi$ and $\theta$,
\ba\lb{composition}\xi(\theta(M))_1 = \theta(\xi(M))_1 =
\mu^{-2}\, \Tr{2,3} K_2 K_1 M_{\overline{3}} =\tr_{(2,3)} K_2 K_1 M_{\overline{3}} =(G^{-1} M G)_1\, .\ea
Here eq.(\ref{KDD}) was used to substitute the R-traces by the usual traces; the last equality follows
from a comparison of the second and the last lines in the calculation (\ref{MjK}).

\smallskip
Now we observe that, in view of eqs.(\ref{tau2}) and (\ref{traceDK}), a matrix
$(^{-1}M)\, :=\, \mu\, g^{-1}\, \theta^{-1}(M)$ fulfills the relation $(^{-1}M)\cdot M=I$. The identity
$(^{-1} M) = M^{-1}$ follows then from the relations (\ref{composition}) and (\ref{Mj}). \hfill$\blacksquare$

\begin{rem}\lb{remark4.18}
{\rm One can generalize the definitions of the characteristic subalgebra and of the matrix powers to the case of
the extended quantum matrix algebra  ${\cal M^{^\bullet\!}}(R,F)$. Not going into details, we just mention
that the {\em extended} characteristic subalgebra ${\cal C^{^\bullet\!}}(R,F)$ is generated by the set
$\{g,g^{-1},p_i\}_{i\geq 0}$ and the {\em extended} algebra ${\cal P^{^\bullet\!}}(R,F)$, as a
${\cal C^{^\bullet\!}}(R,F)$--module, is spanned by matrices
$$M^{\overline{n}}\, \quad \mbox{and}\quad\Mt(M^{\overline{n}})\, \quad \forall\; n\in{\Bbb Z}\, .$$
Here inverse powers of $M$ are defined through the repeated $\star\, $-multiplication by
$M^{\overline{-1}}$, which is given by
\be\lb{Minv*N} M^{\overline{-1}}\star  N \, :=\, N \star  M^{\overline{-1}}\, :=\,\phi^{-1}(M^{-1}\cdot N)
\,  \quad
\forall\; N\in{\cal P^{^\bullet\!}}(R,F)\, .\ee
Explicitly, one has
$$M^{\overline{-n}}\, :=\,\underbrace{M^{\overline{-1}}\star \dots \star  M^{\overline{-1}}\star }_{n\ times}I\, =\,
{\tr_{\!\!R_f^{-1}}}_{(2,\dots,n+1)}\!\left( M^{-1}_{\underline{2}} M^{-1}_{\underline{3}}\dots
M^{-1}_{\underline{n+1}}\, \rho_{R_f^{-1}}(\sigma_n\dots\sigma_2\sigma_1)\right) ;$$
the copies $M^{-1}_{\underline{i}}$ of the matrix $M^{-1}$ were introduced earlier in the lemma
\ref{lemma4.11} (see eqs.(\ref{kopii-2})$\,$). Notice that in general $M^{\overline{-1}} =
\phi^{-1}(M^{-1})\neq M^{-1}$. Here are some particular examples of the  multiplication
by $M^{\overline{-1}}$
$$M^{\overline{-n}}\star M^{\overline{i}}\, =\, M^{\overline{i-n}}\, ,\qquad
M^{\overline{-1}}\star M^{\alpha^{(n)}\uparrow 1}\, =\, ch(\alpha^{(n)})\, I\, .$$ }\end{rem}

\medskip
Let us employ the invertibility of the 2-contraction  to  solve the reciprocal relations derived in the
subsection \ref{subsec4.3}.

\begin{def-prop}\lb{proposition4.19}
For the $O(2\ell)$-type extended quantum matrix algebra ${\cal M^{^\bullet\!}}(R,F)$,
its quotient algebra by a relation
\be\lb{quotients-O}
a_{2\ell}\, =\, + g^{\ell}
\qquad {or,\ respectively,}\qquad a_{2\ell}\, =\, - g^{\ell}
\ee
we call a component, which we shall denote by $O^{+}(2\ell)$ or $O^{-}(2\ell)$, respectively.

\smallskip
In a similar manner, the corresponding quotients of the extended algebra ${\cal P^{^\bullet\!}}(R,F)$ we
shall denote by ${\cal P}^{\pm}(R,F)$ and the corresponding quotients of the extended characteristic
subalgebra ${\cal C^{^\bullet\!}}(R,F)$ --- by ${\cal C}^{\pm}(R,F)$.

\medskip
The reciprocal relations (\ref{reciprocal}) for the components
$O^{\pm}(2\ell)$ reduce, respectively, to
\be\lb{reciprocal2} a_{\ell +i} = \pm g^i\, a_{\ell -i}\,\qquad \forall\; i=0,1,\dots ,\ell\, ,\ee
which, in the $O^-(2\ell)$ case, implies $a_{\ell}=0$.

\medskip
{}For the $O(2\ell -1)$-type extended quantum matrix algebra ${\cal M^{^\bullet\!}}(R,F)$, which we further
refer to as the $O^{^\bullet\!}{(2\ell -1)}$-type quantum matrix algebra, the reciprocal relations
(\ref{reciprocal}) are resolved in a following way
\be\lb{reciprocal3} a_{\ell+i}\, =\, g^{i+1/2}\ a_{\ell-1-i}\,\qquad
\forall\; i=0,1,\dots \ell-1\, ,\ee
where a notation
\be\lb{g-root} g^{1/2}\, :=\,  g^{1-\ell}\, a_{2\ell -1}\,\qquad\mbox{
$\Bigl($ so that $(g^{1/2})^2=g$ and $g^{\ell-1/2}=a_{2\ell -1}\Bigr)$}\ee
is used.

\end{def-prop}

\noindent {\bf Proof.}~ The relations (\ref{quotients-O}) are the two factors of the reciprocal relation
(\ref{reciprocal}) in the case $k=2\ell$ and $i=\ell$:
$0$ = $g^{2\ell}-a_{2\ell}^2$ = $(g^{\ell}+a_{2\ell})$ $(g^{\ell}-a_{2\ell})$.

\smallskip
A check of the proposition is straightforward. \hfill$\blacksquare$

\section{Cayley-Hamilton theorem}\lb{sec5}

In this section we derive, using auxiliary results from the subsection \ref{subsec5.1}, the characteristic
identities (the theorems \ref{theorem5.4}, \ref{theorem5.6} and \ref{theorem5.8}) for the quantum matrix
algebras ${\cal M}(R,F)$ of the types $Sp(2k)$ and $O(k)$.

\medskip
We begin by assuming that $\{R,F\}$ is a compatible pair of R-matrices, in which the operator $F$ is
strict skew invertible and the operator $R$ is skew invertible of the BMW-type (and, hence, strict skew
invertible). Specific $Sp(2k)$ and $O(k)$-type conditions are imposed in the subsections \ref{subsec5.2}
and \ref{subsec5.3}, \ref{subsec5.4}, respectively.

\medskip
Next, in the corollaries \ref{corollary5.5}, \ref{corollary5.7} and \ref{corollary5.9}, we parameterize
the coefficients of the characteristic polynomials by means of a $\Bbb C$-algebra
${\cal E}_{2i}$  of polynomials in $2i+1$ pairwise commuting variables
$\nu_j$, $~j=0,1,\dots ,2i$, satisfying conditions
\be\lb{specSp}\nu_{i+j}\, \nu_{i+1-j}\, =\, \nu_0^2\,\quad\forall\; j=1,2,\dots ,i\, .\ee
We call $\nu_i$, $i=0,1,\dots $, {\em orthosymplectic spectral variables} (below we shall be writing
"spectral variables" for short, always assuming that the relations (\ref{specSp}) are satisfied).

\medskip
For the classical groups $Sp(2k)$, $O^+(2k)$, $O^-(2k+2)$ and $O(2k+1)$,
the functions $a_i$, $i=1,\dots ,k$ on the group manifolds are
functionally independent. This justifies, at least perturbatively, the corresponding assumptions about the independence of the elements $a_i$ in the corollaries \ref{corollary5.5},
\ref{corollary5.7} and \ref{corollary5.9}.

\subsection{Basic identities}\lb{subsec5.1}

Let the eigenvalues $q$ and $\mu$ of the R-matrix $R$ (i.e., the parameters of the BMW algebras, whose
representations are generated by the matrix $R$) satisfy conditions
$i_q\neq 0,\;\mu\neq -q^{3-2i}\;\;\forall
\; i=2,3,\dots ,n$. Then the matrices $M^{a^{(i)}}$ are well defined for $0\leq i\leq n$.

\medskip
In this subsection we establish relations between "descendants" of the matrices $M^{a^{(i)}}$ in the
algebra ${\cal P}(R,F)$. These relations are used later in a derivation of the Cayley-Hamilton theorem
and the Newton relations.

\medskip
{}For $1\leq i\leq n$ ~and~ $m\geq 0$, denote
\be\lb{La} A^{(m,i)}\ :=\ i_q\, M^{\overline{m}}  \star  M^{a^{(i)}}\,\ , \quad B^{(m+1,i)}\ :=\
i_q\, M^{\overline{m}} \star  \Mt (M^{a^{(i)}})\,  .\ee

It is  suitable to define $A^{m,i}$ and $B^{m,i}$ for boundary values of their indices
\ba\lb{AB-boundary} & A^{(-1,i)}\, :=\, i_q\, \phi^{-1}\left(
\Tr{2,3,\dots i} M_{\overline{2}}M_{\overline{3}}\dots M_{\overline{i}}\, \rho_R(a^{(i)})
\right)\ , \quad B^{(0,i)}\ :=\ i_q\, \phi^{-1}\bigl(\xi\bigl(M^{a^{(i)}}\bigr)\bigr)\, &\ea
and
\ba\lb{LT0}
&A^{(m,0)}\ :=0\ \ \ {\mathrm{and}}\ \ \ \ B^{(m,0)}\ :=\ 0\,\qquad \forall\; m\geq 0\, .&\ea
Notice that although the elements $A^{(-1,i)}$ and $B^{(0,i)}$ do not, in general, belong to the algebra ${\cal P}(R,F)$,
their descendants  $A^{(-1,i)} g$ and $B^{(0,i)} g$ do (see eqs.(\ref{rek1}) and (\ref{rek2})  in the case $m=0$).

\smallskip
In the case when the contraction $g$ (and, hence, the matrix $M$) is invertible, the formulas (\ref{La}), with $m$ now an arbitrary integer,
can be used  to define descendants of $M^{a^{(i)}}$ in the extended algebra  ${\cal P^{^\bullet\!}}(R,F)$
(see the remark \ref{remark4.18}). In this case, the matrices $A^{(-1,i)}$ and $B^{(0,i)}$ are expressed uniformly:
$A^{(-1,i)} = i_q M^{\overline{-1}} \star  M^{a^{(i)}}$,~~ $B^{(0,i)} = i_q M^{\overline{-1}} \star  \Mt (M^{a^{(i)}})$.

\begin{lem}\lb{lemma5.1}
{}For ~$0\leq i\leq n-1$~ and ~$m\geq 0$, the matrices $A^{(m-1,i+1)}$ and $B^{(m+1,i+1)}$
satisfy recurrent relations
\ba\lb{rek1} A^{(m-1,i+1)} &=& q^i M^{\overline{m}}\, a_i\, -\,
A^{(m,i)}\, -\, { \mu q^{2i-1}(q-q^{-1})\over 1+\mu q^{2i-1}}\ B^{(m,i)}\, ,\\[1em]\lb{rek2}
B^{(m+1,i+1)} &=&\Bigl( \mu^{-1}q^{-i} M^{\overline{m}}\, a_i\, +\,
{q-q^{-1}\over 1+\mu q^{2i-1}}\ A^{(m,i)}\, -\,B^{(m,i)}\Bigr) g\,  .\ea\end{lem}

\noindent {\bf Proof.}~ We prove eqs.(\ref{rek1}) and (\ref{rek2}) by induction on $i$.
Using the equalities (\ref{M^k* M^p}) and (\ref{Mt-IM}), we find
\be\lb{rek-i=0} A^{(m-1,1)} = M^{\overline{m}}\, , \quad B^{(m+1,1)} = \mu^{-1} M^{\overline{m}} g\, .\ee
These formulas, in the setting (\ref{LT0}), coincide with the formulas (\ref{rek1}) and (\ref{rek2}) for $i=0$.

\smallskip
Let us check eqs.(\ref{rek1}) and (\ref{rek2}) for some fixed $i>0$ assuming that they are
satisfied for all smaller values of the index $i$.

\smallskip
{}For $A^{(m,i+1)}$ and $m \geq 0$, we calculate
\ba\nonumber A^{(m,i+1)} = (i+1)_q M^{(a^{(i+1)\uparrow m}\, \sigma_m\dots \sigma_2\sigma_1)}
\, =\, q^i M^{(a^{(i)\uparrow(m+1)}\,\sigma^{-}_{m+1}(q^{-2i})\,\sigma_m\dots \sigma_2\sigma_1)}
\\[1em]\nonumber =\, q^i M^{\overline{m+1}}\, a_i\, -\, A^{(m+1,i)}\, -\,
{\mu q^{2i-1}(q-q^{-1})\over 1+\mu q^{2i-1}}\, B^{(m+1,i)}\, .\hspace{20mm}\ea
Here in the first line we used the second formula from (\ref{a^k}) for $a^{(i+1)\uparrow m}$ and applied
the reduced cyclic property (\ref{red-cycl}) and eq.(\ref{idemp-1}) to cancel one of two terms
$a^{(i)\uparrow m}$. In the second line we substituted the formula (\ref{ansatz}) for the baxterized
elements $\sigma_{m+1}^-(q^{-2k})$ and applied eq.(\ref{char1}) to simplify the first term in the sum.

\smallskip
{}For $A^{(-1,i+1)}$, the relations (\ref{rek1}) are verified similarly
$$\begin{array}{ccl} A^{(-1,i+1)} &=& q^{i}\, \phi^{-1}\left( \Tr{2,3,\dots i+1} M_{\overline{2}}
M_{\overline{3}}\dots M_{\overline{i+1}}\rho_R(a^{(i)\uparrow 1}\sigma^{-}_1(q^{-2i}))\right)\, \\[1em]
&=&q^i\, \phi^{-1}(I)\,  a_i\, -\, i_q\, \phi^{-1}(\phi(M^{a{(i)}}))\, -\, {\mu q^{2i-1}(q-q^{-1})\over
1+\mu q^{2i-1}}\, \phi^{-1}(\xi(M^{a^{(i)}}))\,\\[1em] &=&q^i\, I\, a_i \, -\, A^{(0,i)}\, -\,
{\mu q^{2i-1}(q-q^{-1})\over 1+\mu q^{2i-1}}\, B^{(0,i)}\, .\end{array}$$
Here the definitions (\ref{phi}) and (\ref{xi}) of the endomorphisms $\phi$ and $\xi$ were additionally taken into account.

\smallskip
{}For the matrices $B^{(m+1,i+1)}$, the calculation starts similarly to the case of the matrices $A^{(m,i+1)}$, $m\geq 0$,
\be\lb{kusok}\!\!\!\!\!\!\!\begin{array}{l} B^{(m+1,i+1)} = (i+1)_q
M^{(a^{(i+1)\uparrow m+1}\, \kappa_{m+1}\sigma_{m}\dots \sigma_2\sigma_1)}\, =\, q^i
M^{(a^{(i)\uparrow m+2}\, \sigma^{-}_{m+2}(q^{-2i})\,\kappa_{m+1}\sigma_{m}\dots \sigma_2\sigma_1)}
\hspace{13mm}\\[1em] = q^{-i} M^{\overline{m}}\!\star \! \Mt(M) a_i-
i_q M^{(a^{(i)\uparrow m+2 }\,\sigma^{-1}_{m+2}\kappa_{m+1}\sigma_{m}\dots\sigma_1)} +
\displaystyle {q^i-q^{-i}\over 1+\mu q^{2i-1}}\, M^{\overline{m}}\! \star \!\Mt(\Mt(M^{a^{(i)}})).\hspace{3mm}\end{array}\ee
Here in the second line we used another expression for the baxterized generators
$$\sigma_{i}^{\varepsilon}(x)\, =\, x 1\, +\, {x-1\over q-q^{-1}}\, \sigma^{-1}_{i}\, -\,
{\alpha_{\varepsilon}x(x-1)\over \alpha_{\varepsilon}x+1}\,\kappa_{i}\, ,$$
which follows by a substitution $\sigma_i = \sigma_i^{-1} + (q-q^{-1})(1-\kappa_i)$ into the original
expression (\ref{ansatz}).

\smallskip
Now, notice that
\be\lb{kakaka}\sigma^{-1}_3\kappa_2\sigma_1\, =\,\sigma^{-1}_3\kappa_2\kappa_1\sigma_2^{-1}\,
\stackrel{\rm mod\, (\ref{red-cycl})}{=}\,\sigma_2^{-1}\sigma^{-1}_3\kappa_2\kappa_1\, =\,
\kappa_3\kappa_2\kappa_1\, ,\ee
and, hence, in the case $m\geq 1$, the second term in (\ref{kusok}) can be expressed as
\be\lb{kusok2} -\, i_q M^{(a^{(i)\uparrow m+2}\,\sigma^{-1}_{m+2}\kappa_{m+1}\sigma_{m}\dots\sigma_1)}\,
=\, - i_q M^{\overline{m-1}}\star \Mt(\Mt(\Mt(M^{a^{(i)}})))\, .\ee
Applying then the formulas (\ref{Mt-IM}) and (\ref{Mt2-N}) to the expressions (\ref{kusok}) and
(\ref{kusok2}), we obtain the recurrent relation (\ref{rek2}).

\smallskip
In the case $m=0$, the transformation of the second term in (\ref{kusok}) suffers a slight modification.
Notice that by eq.(\ref{kakaka}),
$$\phi(M^{a^{(i)\uparrow 2}\sigma_2^{-1}\kappa_1})\, =\, \xi(\Mt(\Mt(M^{a^{(i)}})))\, .$$
Inverting the endomorphism $\phi$ in this formula
and using the relation (\ref{Mt2-N}) and the definition (\ref{AB-boundary}) of
the matrices $B^{(0,i)}$, we can complete
the transformation of the second term in (\ref{kusok}) and, again, get the equality (\ref{rek2}).
\hfill$\blacksquare$

\medskip
In two corollaries below we use the recurrent relations  (\ref{rek1}) and (\ref{rek2}) to derive expansions,
in terms of the matrix powers $M^{\overline{j}}$, of the matrices $A^{(i,m)}$ and $B^{(i,m)}$ and of
their certain linear combinations, for some particular values of the indices $i$ and $m$. Notice that, by
the proposition \ref{proposition4.14}, one could expect an appearance of the terms $\Mt(M^{\overline{j}})$
in the expansions. Our aim is to prepare such combinations of the matrices $A^{(i,m)}$ and $B^{(i,m)}$,
whose expansions do not contain these unwanted terms.

\begin{cor}\lb{corollary5.2}
{}For $1\leq i\leq n$ and $m\geq i-2$, one has\footnote{When ~$m\geq i-2$~ (respectively, ~$m\geq i$),
all the $\star\, $-powers of $M$ in the right hand side of eq.(\ref{cor1a}) (respectively,
eq.(\ref{cor1b})$\,$) are non-negative. This is why we specify these restrictions on $m$. For
an invertible matrix $M$, the restrictions on $m$ can be removed.}
\be\lb{cor1a} A^{(m,i)}\;\, =\;\, (-1)^{i-1}\sum_{j=0}^{i-1} (-q)^j\Bigl\{
M^{\overline{m+i-j}} + {1-q^{-2}\over 1+\mu q^{2i-3}}\sum_{r=1}^{i-j-1}M^{\overline{m+i-j-2r}} (q^2 g)^r
\Bigr\} a_j\, .\ee

\hspace{21mm}For $1\leq i\leq n$ and $m\geq i$, one has
\ba\nonumber  B^{(m,i)}&=& (-1)^{i-1}\sum_{j=0}^{i-1} (-q)^j\Bigl\{
\mu^{-1} q^{-2j} M^{\overline{m-i+j}} g^{i-j}\hspace{56mm}\\[0em] &&\hspace{45mm}
-{q^{-1}(1-q^{-2})\over 1+\mu q^{2i-3}}\sum_{r=1}^{i-j-1}M^{\overline{m+i-j-2r}} (q^2 g)^r
\Bigr\} a_j\, .\lb{cor1b}\ea\end{cor}

\noindent {\bf Proof.}~ We employ induction on $i$.

\smallskip
In the case $i=1$, the relations (\ref{cor1a}) and (\ref{cor1b}) reproduce the formulas (\ref{rek-i=0}).

\smallskip
It is then straightforward to verify the induction step $i\rightarrow i+1$ with the help of the relations
(\ref{rek1}) and (\ref{rek2}). \hfill$\blacksquare$

\begin{rem}
{\rm In the Hecke case, the relation (\ref{cor1a}) simplifies drastically, the terms with $g$ disappear
and the condition $m\geq i-2$ weakens to $m\geq -1$.
{}For $m=0$, eq.(\ref{cor1a}) reproduces the Cayley--Hamilton--Newton identities found in \cite{IOP,IOP1}.
The R-trace maps (\ref{Rtrace-map}) of these identities are the Newton relations.
In the $GL(k)$-case (that is, if $\rho_R(a^{(k+1)})=0$), for $m=-1$ and $i=k+1$,
the left hand side of eq.(\ref{cor1a}) vanishes and (\ref{cor1a}) becomes
the Cayley--Hamilton identity.}\end{rem}

\medskip
{}For $m<i$, only linear combinations of the matrices $A^{(m-2,i)}$ and $B^{(m,i)}$ can be expanded in the matrix powers
$M^{\overline{j}}$. To prepare such combinations we introduce a family of  ${\cal P}(R,F)$-valued
$2\times 2$ matrices
\be\lb{Q(i)} Q^{(i)}\ :=\  {1\over 1+\mu q^{2i-5}} \left( \begin{array}{rr}
-(1+\mu q^{2i-3})\, g I, & \mu q^{2i-3} (q-q^{-1}) I\\[1em] - (q-q^{-1})\, g M^{\overline{2}}, &
-(1+\mu q^{2i-3}) M^{\overline{2}}\end{array}\right) , \quad i=1,2,\dots ,n\, ,\ee
and define $A^{(m,i|s)}$ and $B^{(m,i|s)}$ --- linear combinations of
$A^{(m+2j,i)}$ and $B^{(m+2j,i)}$, $j=0,1,\dots ,s$ --- recursively
\be\lb{AB(s)def}\begin{array}{c} A^{(m,i|0)} := A^{(m,i)} , \;\; B^{(m,i|0)} := B^{(m,i)} , \\[1em]
\left( \begin{array}{c} \!\!A^{(m-1,i|s+1)}\!\!\!\\[1em] \!\!B^{(m+1,i|s+1)}\!\!\!\end{array}
\right) := Q^{(i-s)} \star  \left( \begin{array}{c} \!A^{(m-1,i|s)}\!\!\\[1em] \!B^{(m+1,i|s)}\!\!
\end{array}\right) ,\;\; s=0,1,\dots ,i-1 .\end{array}\ee

\begin{cor}\lb{corollary5.3}
{}For $1\leq i\leq n$, ~$m\geq 0$ and for $s=1,2,\dots ,i$, one has
\ba\nonumber && \hspace{-3mm}\left( \begin{array}{c} \!A^{(m-1,i|s)}\!\!\\[1em] \!B^{(m+1,i|s)}\!\!
\end{array} \right) =(-1)^{s+1} {1+\mu q^{2i-1}\over 1+\mu q^{2i-2s-3}}\sum_{j=1}^s (-q)^{j}
\left( \begin{array}{r} \! q^{i-2j-2} M^{\overline{m+j-1}}\, a_{i-j}\, g^s \!\!\\[1em]
\! \mu^{-1} q^{-i} M^{\overline{m+2s-j+1}}\, a_{i-j} g^j\!\!\end{array}\right) \\[1em]\nonumber
&&+(-1)^s{ q^{i-2s-2}(1-q^{-2})(1+\mu q^{2i-1})\over (1+\mu q^{2i-2s-1}) (1+\mu q^{2i-2s-3})}
\sum_{j=1}^{s-1}\sum_{r=1}^{s-j} (-q)^{j+2r} M^{\overline{m+2s-j-2r+1}} a_{i-j} g^{j+r-1}\! \star \!
\left( \begin{array}{r}\!\!\! I~~ \!\!\\[1em] \!\!\! - q g I \!\!\end{array} \right) \\[1em]
\lb{AB(s)expr} && +{1+\mu q^{2i-1}\over 1+\mu q^{2i-2s-1}}\left( \begin{array}{c} \! A^{(m+s-1,i-s)}\!\!
\\[1em] \! B^{(m+s+1,i-s)}\!\! \end{array} \right)g^s\, .\ea\end{cor}

\noindent {\bf Proof.}~ For $s=0$, the relation (\ref{AB(s)expr}) is tautological. For $s=1$, it reads
\be\lb{AB(s=1)} Q^{(i)} \star \left( \begin{array}{c} \! A^{(m-1,i)}\!\! \\[1em] \! B^{(m+1,i)}\!\!\end{array}
\right) \, =\, -\, {1+\mu q^{2i-1}\over 1+\mu q^{2i-5}}\left( \begin{array}{c}
\! q^{i-3}M^{\overline{m}}\!\! \\[1em] \! \mu^{-1} q^{1-i}M^{\overline{m+2}}\!\!\end{array}\right)
a_{i-1} g \, +\, {1+\mu q^{2i-1}\over 1+\mu q^{2i-3}}\left( \begin{array}{c} \! A^{(m,i-1)}\!\! \\[1em]
\! B^{(m+2,i-1)}\!\! \end{array}\right) g\, .\ee

\noindent
This relation can be directly verified with the use of the definition (\ref{Q(i)})
and of the recurrent formulas (\ref{rek1}) and (\ref{rek2}) for $A^{(m-1,i)}$ and $B^{(m+1,i)}$.

\smallskip
To prove eq.(\ref{AB(s)expr}) in general, one carries out an induction on $s$,
using the relation (\ref{AB(s=1)}) at each step of the induction. Notice that
$$Q^{(i)} \star  \left( \begin{array}{r}\!\!\! I~~ \!\! \\[1em] \!\!\! - q g I \!\!\end{array}\right) \, =\,
-\, {1+\mu q^{2i-1}\over 1+\mu q^{2i-5}}\left( \begin{array}{c} \! -I \!\! \\[1em]
\! q^{-1} M^{\overline{2}}\!\!\end{array}\right) g\ ,$$
which makes calculations more comfortable.\hfill$\blacksquare$

\subsection{Cayley-Hamilton theorem: type $Sp(2k)$}\lb{subsec5.2}

Specifying to the case of the $Sp(2k)$-type quantum matrix algebra, we observe that the matrices $A^{(m-1,i)}$ and
$B^{(m+1,i)}$ are well defined for the values $i=1,2,\dots ,k$ and $m\geq 0$ of their indices. The condition (\ref{spec4}) implies
\be\lb{ahah} A^{(m-1,k+1)}=0\, , \quad B^{(m+1,k+1)}=0\, \quad \forall\; m\geq 0\, .\ee
Howether, these latter  equalities are not independent. Indeed, taking $i=k+1$ and
$\mu =-q^{-1-2k}$ in eq.(\ref{AB(s=1)}) and using the expression (\ref{Q(i)}) for $Q^{(k+1)}$, we find
$$(B^{(m+1,k+1)} + q A^{(m-1,k+1)} g )|_{\mu=-q^{-1-2k}} = 0\,\ \quad \forall\; m\geq 0\, .$$

Consider now one particular condition from the set (\ref{ahah}):
\be\lb{ahah-1} A^{(k-1,k+1)}\, =\, 0\, .\ee
By the corollary \ref{corollary5.2}, its left hand side can be expanded in the powers of the matrix $M$ only (i.e.,
the expansion does not contain terms of the form $\Mt(M^{\overline{j}})$). As we argued above, all the
other relations in (\ref{ahah}), which can be expressed as $\star\, $-polynomials in the matrix $M$, are
consequences of the equality (\ref{ahah-1}).

\smallskip
Substituting $\mu=-q^{-1-2k}$ into the condition (\ref{ahah-1}), using the relation (\ref{cor1a})
and reassembling it in a suitable way, we obtain:

\begin{theor}\lb{theorem5.4}
Let ${\cal M}(R,F)$ be the $Sp(2k)$-type quantum matrix algebra (see the definition \ref{definition4.1}).
Then the matrix $M$  of the algebra generators satisfies the Cayley-Hamilton identity
\be\lb{CHSp-1}\sum_{i=0}^{2k} (-q)^i M^{\overline{2k-i}} \epsilon_i\, =\,  0\, ,\ee

\noindent
where
\be\lb{CHSp-2}\epsilon_i\, :=\, \sum_{j=0}^{[i/2]} a_{i-2j}\, g^j\, , \qquad
\epsilon_{k+i}\, :=\, \epsilon_{k-i}\, g^{i}\,\qquad \forall\; i=1,2,\dots ,k\, .\ee\end{theor}

\medskip
Let us now use the Cayley-Hamilton identity to introduce "eigenvalues" of the matrix $M$. In other words, we are
going to factorize the polynomial in the left hand side of eq.(\ref{CHSp-1}). To this end, we realize elements of
the characteristic subalgebra ${\cal C}(R,F)$ as polynomials in the spectral variables
(see eq.(\ref{specSp})$\,$) and construct a
corresponding extention of the algebra ${\cal P}(R,F)$.

\begin{cor}\lb{corollary5.5}
In the setting of the theorem \ref{theorem5.4}, assume that the  elements $a_i$,
$i=1,2,\dots ,k$, are algebraically independent.

\medskip
Consider an algebra homomorphism of the characteristic subalgebra ${\cal C}(R,F)$ to the algebra of
the spectral variables ${\cal E}_{2k}$,
\be\lb{rep-charSp}
\pi_{Sp(2k)}: {\cal C}(R,F)\rightarrow {\cal E}_{2k}\ ,\ee
defined on the generators by
\be
g\mapsto \nu_0^2\, , \quad a_i\mapsto e_i(\nu_0,-\nu_0,\nu_1, \nu_2,\dots ,\nu_{2k})\,
\quad \forall\; i=1,2,\dots k\, ,\ee
where $e_i(\dots)$ are the elementary symmetric functions of their arguments (for the symmetric functions
we adopt a notation of \cite{Mac}).The map $\pi_{Sp(2k)}$ defines naturally a, say, left
${\cal C}(R,F)$--module structure on the algebra ${\cal E}_{2k}$. Consider a corresponding completion
of the algebra ${\cal P}(R,F)$,
$$ {\cal P}_{Sp(2k)}(R,F)\, :=\, {\cal P}(R,F)\raisebox{-4pt}{$\bigotimes\atop
{\cal C}(R,F)$} {\cal E}_{2k}\, ,$$
where the $\star \, $-product on the completed space is given by a formula
\be\lb{P-Sp2} (N\raisebox{-4pt}{$\bigotimes\atop {\cal C}(R,F)$}\nu)\star
(N'\raisebox{-4pt}{$\bigotimes\atop {\cal C}(R,F)$}\nu') :=
(N\star N')\raisebox{-4pt}{$\bigotimes\atop {\cal C}(R,F)$}(\nu\nu')\ \ \ \forall\;
N,N'\in {\cal P}(R,F)\ \; {\mathrm{and}}\ \; \forall\;\nu, \nu' \in {\cal E}_{2k}\, .\ee
Then, in the completed algebra  ${\cal P}_{Sp(2k)}(R,F)$,
the Cayley-Hamilton identity (\ref{CHSp-1}) acquires a factorized form
\be\lb{CHSp-factor} {\prod_{i=1}^{2k}}\hspace{-12.pt} \star \left( M - q\nu_i I\right)\, =\, 0\, ,\ee
where the  symbol $\displaystyle\prod\hspace{-11.pt} \star $ denotes the product with respect
to the $\star \, $-multiplication (\ref{P-Sp2}).\end{cor}

\noindent {\bf Proof.}~ Using the equalities
\be\nonumber e_i(\nu_0,-\nu_0,\nu_1, \nu_2,\dots ,\nu_{2k})\, =\, e_i(\nu_1, \nu_2,\dots ,\nu_{2k}) -
\nu_0^2\, e_{i-2}(\nu_1, \nu_2,\dots ,\nu_{2k})\ \ \ \forall\ \ i\geq 0
\ee
and
\be\lb{recip-nu} e_{k+i}(\nu_1,\nu_2,\dots ,\nu_{2k})\, =\,\nu_0^{2i}\, e_{k-i}(\nu_1,\nu_2,\dots ,
\nu_{2k})\ \ \ \forall\ \  i=1,2,\dots ,k\ \;\;
\ee
(if $\{\nu_j\}_{j=0}^{2k}$ verifies eq.(\ref{specSp})$\,$),
it is straightforward to check  that the map $\pi_{Sp(2k)}$ sends the coefficients (\ref{CHSp-2})
of the Cayley-Hamilton identity to the elementary symmetric functions in the spectral
variables: $\ \epsilon_i\ $ $\mapsto\ $
$e_i(\nu_1, \nu_2,\dots ,\nu_{2k})\ $ $\quad\forall\ i=1,2,\dots ,2k\, .$ \hfill$\blacksquare$

\subsection{Cayley-Hamilton theorem: type $O(k)$, $k$ even }\lb{subsec5.3}

In this subsection we set $k=2\ell$, $\ell=1,2,\dots\; .$

\medskip
Let us gradually impose conditions making the R-matrix $R$ orthogonal. In this way we observe effects,
which different orthogonality restrictions have for the Cayley-Hamilton identities.

\smallskip
{}First, assume that the conditions (\ref{spec4}) and (\ref{mu1}) are satisfied and
the R-matrix $R$ is not of the symplectic type,
$\mu\neq -q^{-1-2k}$. Hence, $\mu=q^{1-k}$~ and ~$q^{k+2}\neq -1$.

\smallskip
We have again the relations (\ref{ahah}), but this time they imply nontrivial (contrary to the symplectic case) consequences
\be\lb{ahah-2} A^{(m-1,k+1|s)}\, =\, B^{(m+1,k+1|s)}\, =\, 0
\ \ \ \forall\ m\geq 0\ \ {\mathrm{and}}\ \ \forall\ \  1\leq s\leq k-1 .\ee

Let us take $s=\ell$, $m=0$, $\mu=q^{1-k}$ and substitute the expression (\ref{cor1a})
for $A^{(\ell-1,\ell+1)}$ into the formula (\ref{AB(s)expr}) for $A^{(-1,k+1|\ell)}$.
After some recollections of summands, we get
\ba\lb{ortho-1} && A^{(-1,k+1|\ell)}\, =\, (-1)^{\ell}\, {1+q^{k+2}\over 1+q^2}\Bigl(
\sum_{j=0}^k (-q)^j M^{\overline{k-j}}\, a_j g^{\ell}\, +\, {1\over 2}\, U_k\Bigr) ,
\ea
where
\ba\lb{ortho-2} &&
U_k\, :=\, (1-q^{-2})\,\sum_{j=0}^{\ell-1}\sum_{r=1}^{\ell -j}
(-q)^{j+2r} M^{\overline{k-j-2r}} (a_j g^{\ell} - a_{k-j} g^j) g^r\, .\ea
Opening the bracket $(1-q^{-2})$, we represent the matrix $U_k$ as a difference of two double sums.
Changing the summation index $r\rightarrow r-1$ in the second double sum (the one with the common
factor $-q^{-2}$), one can rewrite the expression for $U_k$ in a form
\be\nonumber\begin{array}{l}
U_k\, =\, - {\displaystyle \sum_{j=0}^k} (-q)^j M^{\overline{k-j}}
(a_j g^{\ell}- a_{k-j}g^j)\,\\[1em]
\ \ \ \ \ \ \ \ \ \ \ \ \ \ \ \ \ \ \ \
+\, (M^{\overline{2}}- g I)\star
{\displaystyle \sum_{j=0}^{\ell -1} \sum_{r=0}^{\ell-j-1} } (-q)^{j+2r}
M^{\overline{k-j-2r-2}} (a_j g^{\ell}- a_{k-j}g^j) g^r .\end{array}\ee
Now, substituting this formula back into eq.(\ref{ortho-1}), we obtain
\ba\nonumber A^{(-1,k+1|\ell)} & =& (-1)^{\ell}\, {1+q^{k+2}\over 2(1+q^2)}\,\Bigl( \,
\sum_{j=0}^k (-q)^j M^{\overline{k-j}}( a_j g^{\ell} + a_{k-j} g^j)\, \,\\[1em]\lb{ortho-3}
&& + (M^{\overline{2}}-gI)\star  \sum_{j=0}^{k-2} (-q)^j\, M^{\overline{k-j-2}}
\hspace{-7mm}\sum_{r={\rm max}\{0,j+1-\ell\}}^{[j/2]}\hspace{-8mm}
(a_{j-2r}g^{\ell}-a_{k-j+2r}g^{j-2r}) g^r\Bigr) .\ea
Here, in the last term, a change of the summation index $j\rightarrow j+2r$ was made.

\smallskip
In a similar way, for $B^{(1,k+1|\ell)}$, we find
$$B^{(1,k+1|\ell)}\, =\, (-1)^{\ell}\, {q^{-1}(1+q^{k+2})\over 1+q^2}\,\Bigl( \,
\sum_{j=0}^k (-q)^j M^{\overline{k-j}}\, a_{k-j} g^{j+1}\, -\, {1\over 2}\,
M^{\overline{2}}\!\star U_{k}\Bigr) ,$$
and then, using another expression for $U_{k}$
\ba\nonumber M^{\overline{2}} \star  U_{k} &=& -\sum_{j=0}^k (-q)^j M^{\overline{k-j}}
(a_j g^{\ell}- a_{k-j}g^j) g \\[1em] \nonumber && + q^2 (M^{\overline{2}}- g I)\star
\sum_{j=0}^{\ell -1}\sum_{r=0}^{\ell-j-1} (-q)^{j+2r}
M^{\overline{k-j-2r-2}} (a_j g^{\ell}- a_{k-j}g^j) g^{r+1}\ ,\ea
we obtain
\ba\nonumber B^{(1,k+1|\ell)} & =& (-1)^{\ell}\, {q^{-1}(1+q^{k+2})\over 2(1+q^2)}\,\Bigl(\,
\sum_{j=0}^k (-q)^j M^{\overline{k-j}} ( a_j g^{\ell} + a_{k-j} g^j)\, \,\\[1em]\lb{ortho-4}
&& - q^2 (M^{\overline{2}}-gI)\star  \sum_{j=0}^{k-2} (-q)^j\, M^{\overline{k-j-2}}
\hspace{-7mm}\sum_{r={\rm max}\{0,j+1-\ell\}}^{[j/2]} \hspace{-8mm}
(a_{j-2r}g^{\ell}-a_{k-j+2r}g^{j-2r}) g^r\Bigr)g\, .\ea
Comparing the expressions (\ref{ortho-2}) and (\ref{ortho-4}), we conclude that the relations
(\ref{ahah-2}) imply following identities
\ba\lb{ortho-5} &&\sum_{j=0}^k (-q)^j M^{\overline{k-j}} (a_j g^{\ell} + a_{k-j} g^j) g  =  0\ea
and
\ba\lb{ortho-6} && (M^{\overline{2}}-gI)\star  \sum_{j=0}^{k-2} (-q)^j\, M^{\overline{k-j-2}}
\hspace{-7mm}\sum_{r={\rm max}\{0,j+1-\ell\}}^{[j/2]}\hspace{-8mm}
(a_{j-2r}g^{\ell}-a_{k-j+2r}g^{j-2r}) g^{r+1} =  0 .\ea
Up to now, these were effects of the conditions (\ref{spec4}). Next, imposing additionally the
rank-one condition ${\rm rk}\rho_R(a^{(k)})=1$ and using the consequent reciprocal relations
(\ref{reciprocal}), we can convert eqs.(\ref{ortho-5}) and (\ref{ortho-6}) into
\ba\lb{ortho-7} &&\Bigl\{\sum_{j=0}^k (-q)^j M^{\overline{k-j}} a_j\Bigr\} (g^{\ell} + a_k)g =0\ea
and
\ba \lb{ortho-8}&&\Bigl\{(M^{\overline{2}}-gI)\star \sum_{j=0}^{k-2} (-q)^j\, M^{\overline{k-j-2}}
\hspace{-7mm}\sum_{r={\rm max}\{0,j+1-\ell\}}^{[j/2]} \hspace{-8mm}
a_{j-2r}g^r\Bigr\} (g^{\ell}-a_k)g =  0 .\ea
The equations (\ref{ortho-8}) are the Cayley-Hamilton identities for the quantum matrix algebra of the $O(2\ell)$-type.

\smallskip
Finally, assuming the invertibility of the 2-contraction $g$ and specifying to the  components
$O^{\pm}(2\ell)$ (see the definition-proposition \ref{proposition4.19}), we obtain the Cayley-Hamilton
identities for the components in a more familiar form:

\begin{theor}\lb{theorem5.6}
Consider the components $O^{\pm }(2\ell)$ of the extended orthogonal type quantum matrix algebra
${\cal M^{^\bullet\!}}(R,F)$. Assume, in addition  to the standard restrictions on $q$ (see the
definition \ref{definition3.11}), that $q^{2\ell +2}\neq -1$. Then the matrix $M$ of the algebra
generators satisfies the Cayley-Hamilton  identity, which reads
\ba\lb{CH-O+} &&\sum_{i=0}^{2\ell}(-q)^i\, M^{\overline{2\ell -i}}\, a_i\, =\, 0\, \ea
for the $O^+(2\ell)$ component and
\ba\lb{CH-O-1}&&\left(M^{\overline{2}} - g I\right)\star
\sum_{i=0}^{2\ell-2}(-q)^i\, M^{\overline{2\ell -i-2}}\,\epsilon_i\, =\, 0\, \ea
for the $O^-(2\ell)$ component.

\smallskip
Here
\be\lb{CH-O-2}\epsilon_i\, :=\, \sum_{j=0}^{[i/2]} a_{i-2j}\, g^j\, , \qquad
\epsilon_{\ell -1+i}\, :=\, \epsilon_{\ell-1-i}\, g^{i}\,\qquad \forall\; i=0,1,\dots ,\ell -1\, .\ee
\end{theor}

\begin{cor}\lb{corollary5.7}
In the setting of the theorem \ref{theorem5.6}, assume that the elements $a_i$, $i=1,2,\dots ,2\ell$,
do not satisfy any other algebraic relations except for the relations
(\ref{reciprocal2}).

\medskip
Consider following algebra homomorphisms
of the characteristic subalgebras ${\cal C}^{\pm}(R,F)$ to the algebras of the spectral variables:
\be\lb{rep-charO-}\begin{array}{c}
\pi_{O^-(2\ell)}: {\cal C^-}(R,F)\rightarrow {\cal E}^{^\bullet}_{2\ell -2}\, ,\\[1em]
\ g\mapsto \nu_0^2 , \;\; a_i\mapsto e_i(\nu_0,-\nu_0,\nu_1, \nu_2,\dots ,\nu_{2\ell-2})\,
\ \forall\; i=1,2,\dots 2\ell\,\end{array}\ee
and
\be\lb{rep-charO+}\begin{array}{c}
\pi_{O^+(2\ell)}: {\cal C^+}(R,F)\rightarrow {\cal E}^{^\bullet}_{ 2\ell}\, ,\\[1em]
\ g\mapsto \nu_0^2 , \;\; a_i\mapsto e_i(\nu_1, \nu_2,\dots ,\nu_{2\ell})\,
\ \forall\; i=1,2,\dots 2\ell\, ,\end{array}\ee
where by ${\cal E}^{^\bullet}_i$ we understand an extention of the algebra ${\cal E}_i$
by an element $\nu_0^{-1}$, inverse to the element $\nu_0$: $\nu_0\, \nu_0^{-1}=1$.

\smallskip
The maps $\pi_{O^-(2\ell)}$ and $\pi_{O^+(2\ell)}$ define a left module structure
on the algebras ${\cal E}^{^\bullet}_{2\ell -2 }$ and ${\cal E}^{^\bullet}_{2\ell}$, respectively.

\smallskip
Then, in the corresponding completions of the algebras ${\cal P}^{\pm}(R,F)$,
$$ {\cal P}_{O^+(2\ell)}(R,F)\, :=\, {\cal P}^+(R,F)\raisebox{-4pt}{$\bigotimes\atop
{\cal C}^+(R,F)$} {\cal E}^{^\bullet}_{2\ell}\,$$
and
$$\qquad {\cal P}_{O^-(2\ell)}(R,F)\, :=\,
{\cal P}^-(R,F)\raisebox{-4pt}{$\bigotimes\atop {\cal C}^-(R,F)$} {\cal E}^{^\bullet}_{2\ell-2}\, ,$$
where the $\star \, $-product is given by a formula repeating eq.(\ref{P-Sp2}), the Cayley-Hamilton
identities (\ref{CH-O+}) and (\ref{CH-O-1}) acquire factorized forms
\ba\lb{CHO+factor} && {\prod_{i=1}^{2\ell}}\hspace{-12.pt} \star \left(
M - q\nu_i I \right)\, =\, 0\, \ea
for the $O^{+}(2\ell)$ component and
\ba\lb{CHO-factor} &&
(M-\nu_0 I)\star (M+\nu_0 I)\star \!\!{\prod_{i=1}^{2\ell-2}}\hspace{-14.5pt} \star \left( M - q\nu_i I
\right)\, =\, 0\, .\ea
for the $O^{-}(2\ell)$ component.
\end{cor}

\noindent {\bf Proof.}~ Similar to the proof of the corollary \ref{corollary5.5}. \hfill$\blacksquare$

\subsection{Cayley-Hamilton theorem: type $O(k)$, $k$ odd}\lb{subsec5.4}

In this subsection we set $k=2\ell-1$, $\ell=2,\dots\; .$

\medskip
Here again we analyze the expressions for $A^{(-1,k+1|\ell)}$ and
$B^{(1,k+1|\ell)}$. Setting $\mu =q^{1-k}$ and using the reciprocal relations
(\ref{reciprocal}) (thus, assuming that the rank-one condition is satisfied), we obtain eventually
$$A^{(-1,k+1|\ell)}\, =\, (-1)^{\ell}\, {1+q^{k+2}\over 1+q}\,
\sum_{j=0}^{\ell -1}\sum_{r=0}^{\ell-j-1} (-q)^{j+2r} M^{\overline{k-j-2r-1}}g^r\Bigl(
a_j (a_k g I - g^{\ell} M) + g a_{j-1} (a_k M - g^{\ell} I)\Bigr)\ .$$
The matrix $B^{(1,k+1|\ell)}$ is algebraically dependent:~
$B^{(1,k+1|\ell)} g^{\ell-1}\, =\, A^{(0,k+1|\ell)} a_k .$
Thus, the Cayley-Hamilton identities in the case $O(2\ell-1)$ and $q^{k+2}\neq -1$ read
\be\lb{ortho-9}\sum_{j=0}^{\ell -1}\sum_{r=0}^{\ell-j-1} (-q)^{j+2r}
M^{\overline{k-j-2r-1}}g^{r+1}\Bigl( a_j (a_k I - g^{\ell-1} M) +  a_{j-1} (a_k M - g^{\ell} I)
\Bigr)\, =\, 0\, .\ee
Here we have set $a_{-1}=0$. Assuming then the invertibility of the 2-contraction
$g$ and using the  notation (\ref{g-root}), we obtain the Cayley-Hamilton identities in a familiar form:

\begin{theor}\lb{theorem5.8}
Consider the extended $O^{^\bullet\!}(2\ell-1)$-type quantum matrix algebra ${\cal M^{^\bullet\!}}(R,F)$.
Assume, in addition to standard restrictions on $q$ (see the definition \ref{definition3.11}), that
$q^{2\ell +1}\neq -1$. Then the matrix $M$ of the algebra generators satisfies the Cayley-Hamilton identity
\be\lb{CH-O-odd}\left(M - g^{1/2} I\right)\star \sum_{i=0}^{2\ell-2}(-q)^i\, M^{\overline{2\ell -i-2}}\,
\epsilon_i\, =\, 0\, ,\ee
where
\be\lb{CH-O-odd2}\epsilon_i\, :=\, \sum_{j=0}^{i} a_{i-j}\, (-g^{1/2})^j\, , \qquad
\epsilon_{\ell -1 +i}\, :=\, \epsilon_{\ell-1-i}\, g^i\,\qquad\forall\; i=0,1,\dots ,\ell -1\, .\ee
\end{theor}

\begin{cor}\lb{corollary5.9}
In the setting of the theorem \ref{theorem5.8}, assume that the elements $a_i$, $i=1,2,\dots ,2\ell -1$,
do not satisfy any other algebraic relations except for the reciprocal relations (\ref{reciprocal}).

\medskip
Consider a following algebra homomorphism
of the characteristic subalgebra ${\cal C^{^\bullet\!}}(R,F)$ to the algebra of the spectral variables:
\be\lb{rep-charOodd}\begin{array}{c}
\pi_{O^{^\bullet\!}(2\ell -1)}: {\cal C^{^\bullet\!}}(R,F)\rightarrow
{\cal E}^{^\bullet}_{2\ell -2},\\[1em]
g^{1/2}\mapsto \nu_0 ,\;\;\; a_i\mapsto
e_i(\nu_0,\nu_1, \nu_2,\dots ,\nu_{2\ell-2})\,\;\;\; \forall\; i=1,2,\dots 2\ell-1\, .\end{array}\ee
The map $\pi_{O^{^\bullet\!}(2\ell -1)}$ defines a left module structure on the algebra
${\cal E}^{^\bullet}_{2\ell -2 }$.

\smallskip
Then, in the corresponding completion of the algebra ${\cal P^{^\bullet\!}}(R,F)$,
$$ {\cal P}_{O^{^\bullet\!}(2\ell -1)}(R,F)\, :=\,
{\cal P^{^\bullet\!}}(R,F)\raisebox{-4pt}{$\bigotimes\atop
{\cal C^{^\bullet\!}}(R,F)$} {\cal E}^{^\bullet}_{2\ell -2}\, ,$$
where the $\star \, $-product is given by a formula repeating eq.(\ref{P-Sp2})$\,$,
the Cayley-Hamilton identity (\ref{CH-O-odd}), acquires a factorized form
\ba\lb{CHOodd-factor} (M-\nu_0 I)\star \!\!{\prod_{i=1}^{2\ell-2}}
\hspace{-14.5pt} \star \left( M - q\nu_i I\right)\, =\, 0\, ,\ea\end{cor}

\noindent {\bf Proof.}~ Use relations
$$e_i(\nu_0,\nu_1, \nu_2,\dots ,\nu_{2\ell -1})\, =\,
e_i(\nu_1, \nu_2,\dots ,\nu_{2\ell -1}) +\nu_0\, e_{i-1}(\nu_1, \nu_2,\dots ,\nu_{2\ell -1})
\;\; \forall\; i\geq 0 ,$$
and the equations (\ref{recip-nu}) to check that the map $\pi_{O^+(2\ell -1)}\ $ sends the coefficients
(\ref{CH-O-odd2}) of the Cayley-Hamilton identity to the elementary symmetric functions in the
spectral variables:
$\ \epsilon_i\ $ $\mapsto\ $ $e_i(\nu_1, \nu_2,\dots ,\nu_{2\ell -2})$ $\forall$ $i=0,1,\dots ,2\ell -2\, .$
\hfill$\blacksquare$

\def\theequation{\thesection.\arabic{equation}}
\makeatletter
\@addtoreset{equation}{section}
\makeatother

\section{Newton relations}\lb{sec6}

In this last section we  use the basic identities from the subsection \ref{subsec5.1} to establish relations
between the three sets of elements in  the characteristic subalgebra ---
$\{a_i\}_{i\geq 0}$, $\{s_i\}_{i\geq 0}$ and the power sums $\{p_i\}_{i\geq 0}$. In particular,
we approve the statement of the proposition \ref{proposition4.8}.

\medskip
{}For the theorem \ref{theorem6.1}, we do not need to specify the subtype (symplectic or orthogonal)
of the quantum matrix algebra and the value of the parameter $\mu$.
The results of this theorem are valid for any BMW-type quantum matrix algebra.
Next, in the corollary \ref{corollary6.4}, we parameterize the power sums $p_n$ and
the elements $s_n$ in terms of the spectral variables.

\smallskip
Below, $[x]$ denotes the integer part of the number $x$.

\begin{theor}\lb{theorem6.1}
Let ${\cal M}(R,F)$ be a BMW-type quantum matrix algebra. Assume that its two parameters $q$
and $\mu$ satisfy the conditions (\ref{mu}), which allow to introduce either the set
$\{a_i\}_{i=0}^n$ or, respectively, the set $\{s_i\}_{i=0}^n$ in the characteristic subalgebra
${\cal C}(R,F)$ (see eqs.(\ref{SA_0}) and (\ref{SA_k})$\,$). Then, respectively, the following
Newton recurrent formulas
relating these sets to the set of the  power sums (see eq.(\ref{P-01}) and (\ref{P_k})$\,$)
\ba\lb{Newton-a}\sum_{i=0}^{n-1} (-q)^i a_i\, p_{n-i} &=& (-1)^{n-1} n_q\, a_n \, +\, (-1)^n
\sum_{i=1}^{[n/2]}\Bigl( \mu q^{n-2i} -q^{1-n+2i}\Bigr)\, a_{n-2i}\, g^i\, \ea
and
\ba\lb{Newton-s}
\sum_{i=0}^{n-1} q^{-i} s_i\, p_{n-i} &=&  n_q\, s_n \, +\,\sum_{i=1}^{[n/2]}\Bigl(
\mu q^{2i-n} + q^{n-2i-1}\Bigr)\, s_{n-2i}\, g^i\, \ea
are fulfilled.

\smallskip
In the case, when both sets $\{a_i\}_{i=0}^n$ and $\{s_i\}_{i=0}^n$
are consistently defined, they satisfy the Wronski relations
\be\lb{Wronski}\sum_{i=0}^n (-1)^i a_i\, s_{n-i}\, =\, \delta_{n,0} -\delta_{n,2}\, g\, ,\ee
where $\delta_{i,j}$ is a Kronecker symbol.\end{theor}

\begin{rem}\lb{remark6.2}
{\rm One can use the formulas (\ref{Newton-a}) and (\ref{Newton-s}) for an iterative definition of the
elements $a_i$ and  $s_i$ for $i\geq 1$, with initial conditions $a_0=s_0=1$. In this case, the elements
$a_n$ and $s_n$ are well defined, assuming that $i_q\neq 0 \;\; \forall\; i=2,3,\dots n$. The additional
restrictions on the parameter $\mu$, which appeared in their
initial definition (\ref{SA_k}), are artifacts of
the use of the antisymmetrizers and symmetrizers $a^{(n)}, s^{(n)}\in {\cal W}_n(q)$.}\end{rem}

\noindent {\bf Proof.}~ We prove the  relation (\ref{Newton-a}). Denote
\be\lb{J_i} J^{(0)} := 0, \quad J^{(i)} := {\displaystyle \sum_{j=0}^{i-1}} (-q)^j M^{\overline{i-j}} a_j,
\quad i=1,2,\dots ,n\ .\ee
We are going to find an expression for the matrix $J^{(n)}$ in terms of
the matrices $A^{(0,i)}$ and $B^{(0,i)}$, $1\leq i\leq n$.

\medskip
As we shall see, there exist matrices $H^{(i)}$, which fulfill  equations
\be\lb{H-1} (1-q^2) H^{(i)}\, g\, =\,\Bigl(J^{(i)} + (-1)^i A^{(0,i)} \Bigr)\,  , \quad
i=0,1,\dots ,n.\ee
To calculate $H^{(i)}$,
we substitute repeatedly   eqs.(\ref{rek1}) for
$A^{(0,i)}, A^{(1,i-1)}, \dots , A^{(i-1,1)}$ in the right hand side of
eq.(\ref{H-1}). It then transforms to
\be\lb{H-2} H^{(i)}\, g\, =\, -\mu q^{-1}\sum_{j=1}^{i-1} (-1)^j
{q^{2j-1}\over 1+\mu q^{2j-1}} B^{(i-j,j)}\, .\ee
Now, using the expressions (\ref{rek2}) for $B^{(i-j,j)}$, we see that matrices
\be\lb{H} H^{(i)}\, :=\,\sum_{j=0}^{i-2} {(-q)^{j}\over 1+\mu q^{2j+1}}\Bigl(
M^{\overline{i-j-2}} a_j\, +\, {\mu q^{j}(q-q^{-1})\over 1+\mu q^{2j-1}} A^{(i-j-2,i)}\, -\,
\mu q^j B^{(i-j-2,j)}\Bigr)\ee
satisfy eq.(\ref{H-1}) for $i\geq 2$. In the cases $i=0,1,$ solutions are trivial: $H^{(0)}=H^{(1)}=0$.

\smallskip
Next, consider a combination $(H^{(i+2)}- H^{(i)} g)$. Using eq.(\ref{H})
for the first term and eq.(\ref{H-2}) for the second term, we calculate
\ba\nonumber H^{(i+2)} - H^{(i)} g &=&\sum_{j=0}^{i-1} {(-q)^j\over 1 +\mu q^{2j+1}}\Bigl(
M^{\overline{i-j}} a_j\, +\, {\mu q^j (q-q^{-1})\over 1+ \mu q^{2j-1}} (A^{(i-j,j)}-q^{-1}B^{(i-j,j)})
\Bigr)\, \\[1em] \nonumber &&\hspace{24.4mm} + {(-q)^i\over 1+ \mu q^{2i+1}}\Bigl(
I a_i + {\mu q^i (q-q^{-1})\over 1+ \mu q^{2i-1}} A^{(0,i)}\, -\,\mu q^i B^{(0,i)}\Bigr)\, .\ea
To continue, we need a following auxiliary result:

\begin{lem}\lb{lemma6.3}
{}For $1\leq i\leq n$, one has
\be\lb{A-tri} {(-1)^{i-1} A^{(0,i)}\over 1+\mu q^{2i-1}}\, =\,\sum_{j=0}^{i-1}
{(-q)^j\over 1+\mu q^{2j+1}}\Bigl( M^{\overline{i-j}} a_j\, +\, {\mu q^j(q-q^{-1})\over 1+\mu q^{2j-1}}
(A^{(i-j,j)}-q^{-1} B^{(i-j,j)})\Bigr) .\ee\end{lem}

\noindent {\bf Proof.}~ Use the recursion (\ref{rek1}) for $A^{(i-j-1,j+1)}$ to calculate
$$ {A^{(i-j-1,j+1)}\over 1+\mu q^{2j+1}} +{A^{(i-j,j)}\over 1+\mu q^{2j-1}} =
{q^{j}\over 1+\mu q^{2j+1}}\Bigl( M^{\overline{m+i-j}}a_j +{\mu q^j (q-q^{-1})\over 1+\mu q^{2j-1}}
(A^{(i-j,j)}-q^{-1} B^{(i-j,j)}) \Bigr)\, .$$
Compose an alternating sum of the above relations for $0\leq j \leq i-1$
and take into account the condition $A^{(i,0)}=0$. \hfill$\blacksquare$

\medskip
Using the relation (\ref{A-tri}), we finish the calculation
\be\lb{H-H} H^{(i+2)} - H^{(i)} g \, =\, (1+\mu q^{2i+1})^{-1}\Bigl(
(-q)^i I a_i\, +\, (-1)^{i+1}(A^{(0,i)}+\mu q^{2i} B^{(0,i)})\Bigr)\,\;\; \forall\; i=0,1,\dots ,n-2.\ee

Now it is straightforward to get
\be\lb{H-otvet} H^{(i)}\, =\,\sum_{j=1}^{[i/2]}{(-1)^{i-1}\over
1+\mu q^{2(i-2j)+1}}\Bigl( A^{(0,i-2j)} + \mu q^{2(i-2j)} B^{(0,i-2j)} - q^{i-2j} I a_{i-2j}
\Bigr)\, g^{j-1}\, , \;\; i=2,3,\dots ,n.\ee

{}Finally, substituting the expression (\ref{H-otvet}) back into eq.(\ref{H-1}), we obtain a formula
\be\lb{J-otvet} J^{(i)}\, =\, (-1)^{i-1} A^{(0,i)}\, +\,
\sum_{j=1}^{[i/2]}{(-1)^{i-1}(1-q^2)\over 1+\mu q^{2(i-2j)+1}}\Bigl(
A^{(0,i-2j)} + \mu q^{2(i-2j)} B^{(0,i-2j)} - q^{i-2j} I a_{i-2j}\Bigr) g^j\, , \ee
which is valid for $0\leq i\leq n$ (the cases $i=0,1$ can be easily  checked).

\smallskip
Taking the R-trace of eq.(\ref{J-otvet}), we obtain the Newton relations (\ref{Newton-a}). Here, in
the calculation of the R-trace of $B^{(0,i-2j)}$, we took into account the formulas (\ref{phi-xi-traces}).

\medskip
The formulas (\ref{Newton-s}) can be deduced from the relations (\ref{Newton-a}) by a substitution
$q\rightarrow -q^{-1}$,~ $a_j\rightarrow s_j$. This is justified by the existence of the BMW algebras
homomorphism (\ref{homS-A}) ${\cal W}_n(q,\mu)\rightarrow {\cal W}_n(-q^{-1},\mu)$ and a fact that one and
the same R-matrix $R$ generates representations of both algebras ${\cal W}_n(q,\mu)$ and
${\cal W}_n(-q^{-1},\mu)$.

\medskip
The relation (\ref{Wronski}) is proved by induction on $n$. The cases $n=0,1,2$ are easily checked with the use of
eqs.(\ref{Newton-a}) and (\ref{Newton-s}). Then, making an induction assumption, we derive the Wronski relations
for arbitrary $n>2$. To this end, we take a difference of  eqs.(\ref{Newton-s}) and (\ref{Newton-a})
$${\displaystyle \sum_{i=0}^{n-1}} \Bigl( q^{-i} s_i\, p_{n-i}\, -\, (-q)^i a_i\, p_{n-i}
\Bigr)\, =\, n_q (s_n\, +\, (-1)^n a_n )\, +\,\mbox{terms proportional to $g$}$$
and substitute for $p_{n-i}$ in the first/second term of the left hand side its expression from the Newton
relation (\ref{Newton-a})/(\ref{Newton-s}) (with  $n$ replaced by $n-i$). As a result, all terms,
containing the power sums, cancel and, after rearranging the summations, we get
$$n_q \sum_{i=0}^n (-1)^i a_i s_{n-i}\, =\, -\sum_{i=1}^{[n/2]}
(q^{1-n+2i}+q^{n-1-2i})g^i\sum_{j=0}^{n-2i}(-1)^j a_j s_{n-2i-j}\, .$$
By the induction assumption, the double sum in the right hand side of this relation vanishes identically
(when $n$ is odd, the second sum vanishes for all values of the index $i$; when $n$ is even, the second sum
is different from zero only for two values of the index $i$, $i=n/2$ and $i=n/2-1$, and
these two summands cancel). \hfill$\blacksquare$

\begin{cor}\lb{corollary6.4}
Let  ${\cal M}(R,F)$ belong to one of the four specific series  $Sp(2k)$, $O^\pm(2\ell)$ or
$O^{^\bullet}(2\ell -1)$ of the BMW-type  quantum matrix algebras. Assume that the algebra parameter $q$
fulfills the conditions $i_q\neq 0$, $i=2,3,\dots n$  for some $n$~\footnote{For up to, respectively, $n=k$,
$n=2\ell$ and $n=(2\ell -1)$, these conditions enter the initial settings for the algebras of the types
$Sp(2k)$, $O^\pm(2\ell)$ and $O^{^\bullet}(2\ell -1)$ (see the definitions \ref{definition3.11} and
\ref{definition4.1}). For the algebras of the orthogonal $O^\pm(2\ell)$ and $O^{^\bullet}(2\ell -1)$-types,
the additional restrictions, $q^{2\ell +2}\neq -1$ and, respectively, $q^{2\ell +1}\neq -1$, were imposed
in the theorems \ref{theorem5.6} and \ref{theorem5.8}. Formally, they are not used in the proof of the
corollary. However, they are necessary if one wants to establish a relation between the spectral
variables $\nu_i$ and the eigenvalues of the matrix $M$.} and define the elements $a_n$ and  $s_n$ of
the characteristic subalgebra (respectively, ${\cal C}(R,F)$, ${\cal C}^\pm(R,F)$ or
${\cal C}^{^\bullet}(R,F)\,$) with the help of the Newton recurrent relations (\ref{Newton-a}) and (\ref{Newton-s}) (see
the remark \ref{remark6.2}). Then the elements $s_n$ and $p_n$ have following images under the respective
homomorphism $\pi_{Sp(2k)}$, $\pi_{O^\pm(2\ell)}$ or $\pi_{O^{^\bullet}(2\ell -1)}$ of the
corresponding characteristic
algebra to the algebra of the spectral variables (for the definitions of these homomorphisms, see the
corollaries \ref{corollary5.5}, \ref{corollary5.7} and \ref{corollary5.9}):

\vspace{5mm}
\noindent
\underline{the type $Sp(2k)$}
\ba\lb{para-s-Sp}\pi_{Sp(2k)}:&& s_n \mapsto h_n(\nu_1,\nu_2,\dots ,\nu_{2k}) ,\qquad
p_n \mapsto q^{n-1}\sum_{i=1}^{2k} d_i \nu_i^n\, ,\ea
where
\ba\lb{para-p-Sp} &&
\quad d_i\, := {\nu_i - q^{-4}\nu_{2k+1-i}\over \nu_i -\nu_{2k+1-i}}
\prod _{j=1\atop j\neq i,\, 2k +1-i}^{2k} {\nu_i - q^{-2}\nu_j\over \nu_i - \nu_j}\, ;\ea
\underline{the type $O^+(2\ell)$}
\ba\lb{para-s-O+} \pi_{O^+(2\ell)}:&& s_n \mapsto h_n(\nu_1,\dots ,\nu_{2\ell})
- \nu_0^2 h_{n-2}(\nu_1,\dots ,\nu_{2\ell}) ,\qquad p_n \mapsto q^{n-1}\sum_{i=1}^{2\ell} d_i \nu_i^n\, ,
\ea
where
\ba\lb{para-p-O+}&&\quad d_i\, := \prod _{j=1\atop j\neq i,\, 2\ell +1-i}^{2\ell}
{\nu_i - q^{-2}\nu_j\over \nu_i - \nu_j}\, ;\ea
\underline{the type $O^-(2\ell)$}
\ba\lb{para-s-O-} \pi_{O^-(2\ell)}:&& s_n \mapsto h_n(\nu_1,\dots ,\nu_{2\ell-2}) , \qquad
p_n \mapsto q^{n-1}\sum_{i=1}^{2\ell-2} d_i \nu_i^n + \bigl(1 + (-1)^n\bigr) q^{1-2\ell} \nu_0^n\, ,
\ea
where
\ba\lb{para-p-O-} && \hspace{2.5mm} \quad d_i\, :=
{\nu_i - q^{-4}\nu_{2\ell -1-i}\over \nu_i - \nu_{2\ell -1-i}}
\prod _{j=1\atop j\neq i,\, 2\ell -1-i}^{2\ell -2} {\nu_i - q^{-2}\nu_j\over \nu_i - \nu_j}\, ;\ea
\underline{the type $O^{^{\bullet\!}}(2\ell-1)$}
\ba\lb{para-s-O-odd}\pi_{O^{^\bullet}(2\ell-1)}:&& s_n \mapsto h_n(\nu_1,\dots ,\nu_{2\ell-2})
+\nu_0 h_{n-1}(\nu_1,\dots ,\nu_{2\ell-2}) ,\quad p_n \mapsto q^{n-1}\sum_{i=1}^{2\ell-2} d_i \nu_i^n +
q^{2-2\ell}\nu_0^n\, ,\hspace{10mm}\ea
where
\ba\lb{para-p-O-odd}&&\hspace{16.5mm}
\quad d_i\, := {\nu_i - q^{-2}\nu_o\over \nu_i - \nu_o} \prod _{j=1\atop j\neq i,\, 2\ell -1-i}^{2\ell -2}
{\nu_i - q^{-2}\nu_j\over \nu_i - \nu_j}\, .\ea
Here, in eqs.(\ref{para-s-Sp}), (\ref{para-s-O+}), (\ref{para-s-O-}) and (\ref{para-s-O-odd}),
$h_n(\dots)$ denotes the complete symmetric function in its arguments.

\medskip
The power sums contain the rational functions $d_i$ in the spectral variables and are themselves
rational functions in $\{ \nu_i\}$. However, as it follows from the
Newton recursion (\ref{Newton-a}), the power sums simplify, modulo the
relations (\ref{specSp}), to polynomials in the spectral variables.\end{cor}

\noindent {\bf Proof.}~ For the proof, we use a following auxiliary statement:

\begin{lem}\lb{lemma6.5}
In the assumptions of the remark \ref{remark6.2} consider iterations
\ba\lb{mod-a}a'_0=a_0\, , && a'_1=a_1\, , \qquad a'_i = a_i + a'_{i-2} g\, ;\\[1em] \lb{mod-s}
s'_0=s_0\, , && s'_1=s_1\, , \qquad s'_i = s_i + s'_{i-2} g\, ;\\[1em] \lb{mod-p1}
p'_0=(1-\mu^2 q^{-2})/(q-q^{-1})\, ,  && p'_1=p_1\, , \qquad p'_i = p_i + (q^2p'_{i-2}-p_{i-2}) g\, ;
\\[1em] \lb{mod-p2} p''_0=(1-\mu^2 q^{2})/(q-q^{-1})\, ,  &&
p''_1=p_1\, , \qquad p''_i = p_i + (q^{-2}p''_{i-2}-p_{i-2}) g\,\quad \mbox{for}\;\; i\geq 2.\ea
The modified sequences $\{a'_i\}_{i=0}^n$ and $\{s'_i\}_{i=0}^n$ and the two modified sequences of the
power sums, $\{p'_i\}_{i=0}^n$ and $\{p''_i\}_{i=0}^n$, satisfy following versions of the Newton and Wronski
relations
\ba\lb{mod-N} \sum_{i=0}^{n-1} (-q)^i a_i p'_{n-i}\, =\,  (-1)^{n-1} n_q a_n\, ,
\qquad \sum_{i=0}^{n-1} q^{-i} s_i p''_{n-i}\, =\,  n_q s_n\, && \forall\; n\geq 1\, ;\\[1em]
\lb{mod-W}\sum_{i=0}^{n} (-1)^i a'_i s_{n-i}\, =\,\sum_{i=0}^{n} (-1)^i a_i s'_{n-i}\, =\,\delta_{n,0}\,
&& \forall\; n\geq 0\, .\ea\end{lem}

\noindent {\bf Proof.}~ For $n<2$, the equalities (\ref{mod-N})--(\ref{mod-W}) are clearly satisfied. For
$n\geq 2$, one can check them inductively, applying the iterative formulas (\ref{mod-a})--(\ref{mod-p2}).
E.g., the relations (\ref{mod-W}) follow  from  a calculation
$$\sum_{i=0}^n (-1)^i a'_i s_{n-i}\, =\,\sum_{i=0}^n (-1)^i (a_i + a'_{i-2} g) s_{n-i}\, =\,
\delta_{n,0} -\delta_{n,2}\, g +\sum_{i=0}^{n-2} (-1)^i a'_i s_{n-2-i}\, g
\, =\, \delta_{n,0}\, ,$$
where the relations (\ref{mod-a}) and (\ref{Wronski}) and the induction assumption were successively used.

\smallskip
{}For the modified Newton relations (\ref{mod-N}), arguing in the same way, one finally
finds conditions on the zero-th modified power sums
$$p'_0 = q^{-2}(p_0-\mu+q)\, , \qquad p''_0 = q^{2}(p_0-\mu -q^{-1})\, ,$$
which fix the initial settings for the iterations (\ref{mod-p1}) and (\ref{mod-p2}).\hfill$\blacksquare$

\medskip
We now notice that in all four cases, considered in the corollary, the images of the elements
$a_i$, $i=1,\dots ,n$, are given by the elementary symmetric functions (see eqs.(\ref{rep-charSp}),
(\ref{rep-charO-}), (\ref{rep-charO+}) and (\ref{rep-charOodd})$\,$). Hence, by the Wronski relations,
(\ref{mod-W}) the images of the modified elements $s'_n$, $i=1,\dots ,n$, are the complete symmetric
functions in the same arguments. Using then  eq.(\ref{mod-s}) and taking into account the relation
$h_n(\nu_0,\nu_1,\dots)= \sum_{i=0}^n \nu_0^i h_{n-i}(\nu_1,\dots)$, it is easy to check the formulas
for the images of the elements $s_n$, which are given in eqs.(\ref{para-s-Sp})--(\ref{para-s-O-odd}).

\smallskip
To check the formulas for the power sums, we use a  following statement, which was proved in \cite{GS}: if
the elements $a_i$ (respectively, $s_i$) for $i=0,1,\dots ,n\geq 1$ are realized as the elementary
symmetric polynomials $e_i$ (respectively, the complete symmetric polynomials $h_i$)
in some set of variables $\{\nu_i\}_{i=1}^p$, then the elements $p'_n$ (respectively, $p''_n$), defined
by eqs.(\ref{mod-N}), have following expressions in terms of the variables $\nu_i$
\be\lb{formula-GS} p'_n  (\mbox{respectively,~} p''_n) = q^{n-1} \sum_{i-1}^p \widehat d_i \nu_i^n ,
\qquad \mbox{where}\quad \widehat d_i :=
\prod_{j=1\atop j\neq i}^{p}{\nu_i -q^{-2} \nu_j\over \nu_i - \nu_j}.\ee
The latter formulas for $p'_i$ (respectively, for $p''_i$) suit better for the calculation of the images
of $p_i$ in the cases $O^+(2\ell)$ and $O^{^\bullet}(2\ell -1)$ (respectively, $Sp(2k)$ and
$O^-(2\ell$)). In all cases, the calculation goes essentially in the same way. Below we shall describe
it in the $Sp(2k)$ case.

\medskip
Assuming that the relation (\ref{formula-GS}) stays valid for $p''_0$ (note, $p''_0$ is not fixed by the
recursion (\ref{mod-N})$\,$) and taking the ansatz (\ref{para-s-Sp}) for the power sums $p_i$ for
$i=0,1,\dots ,n$, we make use of the recursion (\ref{mod-p2}). Upon substitutions, we find that
the relations (\ref{mod-p2}) hold valid provided that
\be\lb{d-hat-d}d_i \, =\, {\nu_i^2 - q^{-4}\nu_0^2\over\nu_i^2 -q^{-2}\nu_0^2}\,\,\widehat d_i\, .\ee
Taking into account the relations (\ref{specSp}), which are imposed on $\nu_i\in {\cal E}_{2k}$ in the
$Sp(2k)$-case, we observe that the conditions (\ref{d-hat-d}) dictate the choice (\ref{para-p-Sp})
for $d_i$.

\smallskip
It remains to check the initial settings for the recursion (\ref{mod-p2}). They are:
\ba\lb{init-1}
&& p''_0\, =\,  q^{-1} \displaystyle \sum_{i=1}^{2k} \widehat d_i\,
=\, {1-\mu^2 q^2\over q-q^{-1}}|_{\mu=-q^{-1-2k}}\, =\, q^{-2k} (2k)_q\, \ea
from eq.(\ref{mod-N});
\ba\lb{init-2}
&& p_0\, =\,  q^{-1} \displaystyle \sum_{i=1}^{2k} d_i\,
=\, {\rm Tr}_R I|_{\mu=-q^{-1-2k}}\, =\, q^{-1-2k}\bigl( (2k+1)_q -1\bigr)\, \ea
from eq.(\ref{P-01}); and
\ba\lb{init-3}&& p_1=p''_1\quad \Longrightarrow\quad \displaystyle
\sum_{i=1}^{2k}\nu_i (d_i - \widehat d_i)\, =\, 0\, \ea
from eq.(\ref{mod-N}).

To verify them, we use expansions of following rational functions
$$w_1(z) := \prod_{i=1}^{2k} {z-q^{-2} \nu_i\over z-\nu_i}\, , \qquad
w_2(z) := {\nu_0^2 w_1(z)\over z^2 - q^{-2} \nu_0^2}\, , \qquad w_3(z) := z w_2(z)$$
in simple ratios. Expanding $w_1(z)$ and evaluating the result at $z=0$, we prove immediately the
condition (\ref{init-1}). We make one more comment on a less trivial check of
the condition (\ref{init-2}). Expanding $w_2(z)$, we obtain
$$w_2(z)\, =\, \sum_{i=1}^{2k} q^2 (d_i-\widehat d_i) {\nu_i\over z-\nu_i}\, +\,  {q\nu_0\over 2}\Bigl(
{w_1(q^{-1}\nu_0)\over z-q^{-1}\nu_0} -{w_1(-q^{-1}\nu_0)\over z+q^{-1}\nu_0}\Bigr)\, .$$
Here, for the transformation of the first term in the right hand side, we used the formulas
(\ref{formula-GS}) and (\ref{d-hat-d}) and applied the relations (\ref{specSp}), which confine variables
$\nu_i\in {\cal E}_{2k}$ in the $Sp(2k)$ case. The  relations (\ref{specSp}) also allow us to calculate
$w_1(\pm q^{-1}\nu_0)=q^{-2k}$. Thus, evaluating $w_2(z)$ at $z=0$, we obtain
$$w_2(0)\, =\, -q^{2-4k}\, =\,  -q^3 (p_0-p''_0)\, -\, q^{2-2k}\, ,$$
wherefrom the condition (\ref{init-2}) follows. A check of the condition (\ref{init-3}), by the expansion
and evaluation of $w_3(z)$ at $z=0$, is a similar calculation. \hfill$\blacksquare$

\appendix
\section{On twists in quasitriangular Hopf algebras}\label{append}

Here we shall  discuss universal (i.e., quasi-triangular Hopf algebraic) counterparts
of relations from the subsection \ref{subsec3.2}, especially from the lemma \ref{lemma3.6}:
we shall see, in the part {\bf 7} of the appendix, that these relations have a quite transparent
meaning, they reflect the properties of the twisted universal R-matrix.

\medskip
We don't give an introduction to the theory of quasitriangular Hopf algebras assuming that the reader
has some basic knowledge on the subject (see, e.g., \cite{CP}, the chapter 4).

\paragraph{{\bf 1.}} Let ${\cal A}$ be a Hopf algebra; $m, \D ,\e$ and $S$ denote the multiplication, comultiplication,
counit and antipode, respectively.

\smallskip
Assume that ${\cal A}$ is quasitriangular with a universal R-matrix $\cR =a\ot b$ (this is a symbolic
notation, instead of $\sum_i a_i\ot b_i$). One has $(S\ot S)\cR =\cR$. The universal R-matrix
$\cR$ is invertible, its inverse is related to $\cR$ by formulas~ $\cR^{-1}=S(a)\ot b$ ~or~
$(\id\ot S)(\cR^{-1})=\cR$.

\smallskip
{}For elements in ${\cal A}\ot {\cal A}$, the "skew" product $\odot$ is defined as the product in
${\cal A}^{\mathrm{op}}\ot {\cal A}$, where ${\cal A}^{\mathrm{op}}$ denotes the algebra with the opposite multiplication. In
other words, the skew product of two elements, $x\ot y$ and $\ti{x}\ot\ti{y}$ is
$(x\ot y)\odot (\ti{x}\ot\ti{y})=\ti{x}x\ot y\ti{y}$.
{}For a skew invertible element ${\cal X}\in {\cal A}\ot {\cal A}$, we shall denote its skew inverse by $\psi_{\cal X}$.
The universal R-matrix $\cR$ has a skew inverse, $\psi_\cR =a\ot S(b)$. The element $\psi_\cR$
is invertible, $(\psi_\cR)^{-1}=a\ot S^2(b)$. The element $\cR^{-1}$ is skew invertible as well,
its skew inverse is $\psi_{(\cR^{-1})} =S^2(a)\ot b$. All these formulas are present in \cite{D3}.
We shall see below that there are similar formulas for the twisting element $\cF$. However,
the properties of the twisting element $\cF$ and of the universal R-matrix $\cR$ are different,
for instance, the square of the antipode is given by $S^2(x)=u_{_\cR}\, x\, (u_{_\cR} )^{-1}$, where
$u_{_\cR} =S(b)a$, but there is no analogue of such formula for $\cF$. Because of this difference, we
felt obliged to give some proofs of the relations for $\cF$.

\medskip
Let $\rho$ be a representation of the algebra ${\cal A}$ in a vector space $V$. For an element ${\cal X}\in {\cal A}\ot {\cal A}$,
denote by $\hat{\rho}({\cal X})\in {\rm End}(V^{\otimes 2})$ an operator $\hat{\rho}({\cal X})=P\cdot (\rho\ot\rho)({\cal X})$
(recall that $P$ is the permutation operator). The skew product $\odot$ translates into a following
product $\hat{\odot}$ for the elements of ${\rm End}(V^{\otimes 2})$:
$(X\hat{\odot}Y)_{13}=\tr_{(2)}(X_{12}Y_{23})$. In other words, if ${\cal X}\odot {\cal Y}={\cal Z}$ then
$\hat{\rho}({\cal X})\;\hat{\odot}\;\hat{\rho}({\cal Y})=\hat{\rho}({\cal Z})$. For an operator
$X\in {\rm End}(V^{\otimes 2})$, its skew inverse $\Psi_X$, in the sense explained in
the subsection \ref{subsec3.1}, is presicely the inverse with respect to the product $\hat{\odot}$.

\paragraph{{\bf 2.}} A following lemma is well known (see, e.g., \cite{CP},
the chapter 4, and references therein).

\begin{lem}\lb{lemma3.2.0.1}
Consider an invertible element $\cF =\a\ot\b\in {\cal A}\ot {\cal A}$ (we use the symbolic notation,
$\a\ot\b =\sum_i \a_i\ot\b_i$, like for the universal R-matrix) and
let $\cF^{-1}=\g\ot\d$. Assume that the element $\cF$ satisfies
\be \cF_{12}\; (\D\ot\id)(\cF )=\cF_{23}\; (\id\ot\D )(\cF )\ .\lb{uu1}\ee
Assume also that
\be \e (\a )\,\b =\a\,\e(\b )=1\ .\lb{uu2'}\ee
Then an element $\vf =\a\, S(\b )$ is invertible, its inverse is
\be (\vf )^{-1}=S(\g )\,\d\ .\lb{vfin}\ee
One also has
\be S(\a )\, (\vf )^{-1}\,\b =1\ \ \ {\mathrm{and}}\ \ \ \g\,\vf\, S(\d )=1\ .\lb{uu14}\ee
Twisting the coproduct by the element $\cF$,
\be \D_\cF (a)=\cF\;\D (a)\;\cF^{-1}\ ,\ee
one obtains another quasitriangular structure on ${\cal A}$ with
\be \cR_\cF =\cF_{21}\;\cR\;\cF^{-1}\ \lb{anor}\ee
and
\be\ S_\cF (a)=\vf\; S(a)\; (\vf )^{-1}\ \lb{uu17}\ee
(the counit does not change).
\end{lem}

An element $\cF$, satisfying conditions (\ref{uu1}) and (\ref{uu2'}) is called {\em twisting} element.
We shall denote by ${\cal A}_\cF$ the resulting "twisted" quasitriangular Hopf algebra.

\begin{rem}{\rm On the representation level, the formula (\ref{anor}) transforms (compare with
eq.(\ref{R_f})$\,$) into $\hat{\rho}(\cR_\cF )=
\hat{\rho}(\cF )_{21}\hat{\rho}(\cR )_{21}\hat{\rho}(\cF )_{21}^{-1}$.
Below, when we talk about matrix counterparts of universal formulas, one should keep in mind this
difference in conventions.\lb{rrr1}}\end{rem}

\paragraph{{\bf 3.}} Assume, in addition to eq.(\ref{uu1}), that
\be (\D \ot\id )\; (\cF )=\cF_{13}\;\cF_{23}\ \lb{uu18}\ee
and
\be (\id\ot\D )\; (\cF )=\cF_{13}\;\cF_{12}\ .\lb{uu19}\ee

\smallskip
Now the conditions (\ref{uu2'}) follow from eqs.(\ref{uu18}) and (\ref{uu19}) and the invertibility of
the twisting element $\cF$: applying $\e\ot\id\ot\id$ to eq.(\ref{uu18}), we find $(\e\ot\id ) (\cF)=1$; applying $\id\ot\id\ot\e$ to eq.(\ref{uu19}), we find $(\id\ot\e ) (\cF)=1$.

\medskip
Since $\D^{{\mathrm{op}}}(x)\cR =\cR\D (x)$ for any element $x\in {\cal A}$ (where $\D^{{\mathrm{op}}}$ is the
opposite comultiplication), it follows from eq.(\ref{uu18}) that
\be \cR_{12}\;\cF_{13}\;\cF_{23}=\cF_{23}\;\cF_{13}\;\cR_{12}\ .\lb{uu20}\ee

\smallskip
Similarly, eq.(\ref{uu19}) implies
\be \cR_{23}\;\cF_{13}\;\cF_{12}=\cF_{12}\;\cF_{13}\;\cR_{23}\ .\lb{uu22}\ee

\smallskip
When both relations (\ref{uu18}) and (\ref{uu19}) are satisfied, the relation (\ref{uu1}) is equivalent to
the Yang--Baxter equation for the twisting element $\cF$:
\be \cF_{12}\;\cF_{13}\;\cF_{23}=\cF_{23}\;\cF_{13}\;\cF_{12}\ .\lb{uu24}\ee

\begin{rem}{\rm One also has
\be (\D_\cF\ot\id ) (\cF_{21})=\cF_{31}\;\cF_{32}\ \ {\mathrm{and}}\ \
(\id\ot\D_\cF )(\cF_{21})=\cF_{31}\;\cF_{21}\ .\ee
Therefore, one can twist $\D_\cF$ again, now by the element $\cF_{21}$.

\smallskip
On the matrix level, this corresponds to the second conjugation of $\hat{\rho}(\cR )$ by
$\hat{\rho}(\cF )$,~~
$$\hat{\rho}\bigl(\, (\cR_\cF )_{\cF_{21}}\,\bigr)\,
=\,\hat{\rho}(\cF )^2\; \hat{\rho}(\cR )\;\hat{\rho}(\cF )^{-2}\ .$$
\lb{rrr2}}\end{rem}

\begin{rem}{\rm The element $\cF_{21}^{-1}$ satisfies the conditions (\ref{uu1}), (\ref{uu18}) and
(\ref{uu19}) if the element $\cF$ does. Thus, one can twist the coproduct $\D$ by the element $\cF_{21}^{-1}$ as well.\lb{rrr3}}\end{rem}

\paragraph{{\bf 4.}} The conditions (\ref{uu2'}), (\ref{uu18}), (\ref{uu19})
imply the invertibility and skew-invertibility of the element $\cF$. The formulas
for its inverse and skew inverse are similar to the corresponding formulas for the universal R-matrix
$\cR$ (in particular, we reproduce the standard formulas for $\cR$ since we can take $\cF =\cR$).

\begin{lem}\lb{lemma3.2.0.2}
Assume that the conditions (\ref{uu2'}) and (\ref{uu18})
are satisfied. Then the element $\cF$ is invertible, its inverse is
\be \cF^{-1}=S(\a )\ot\b \ .\lb{uu26}\ee
Assume that the conditions (\ref{uu2'}) and (\ref{uu19})
are satisfied. Then the element $\cF$ is skew invertible, with the skew inverse
\be \psi_\cF =\a\ot S(\b )\ .\lb{uu34}\ee
Assume that the conditions (\ref{uu2'}), (\ref{uu18}) and (\ref{uu19})
are satisfied. Then
\be (S\ot S)(\cF )=\cF\ .\lb{uu29}\ee
Moreover, the element $\psi_\cF $ is invertible, its inverse is
\be (\psi_\cF )^{-1}=\a\ot S^2(\b )\ \lb{uu35}\ee
and the element $\cF^{-1}$ is skew-invertible, its skew inverse reads
\be \psi_{(\cF^{-1})} =S^2(\a)\ot \b \ .\lb{uu34'}\ee
\end{lem}

\nin {\bf Proof.~} The calculations are similar to those, from textbooks, for the universal R-matrix.
We include this proof for a completness only.

\smallskip
Applications of $m_{12}\circ S_1$ and $m_{12}\circ S_2$ to eq.(\ref{uu18}) imply eq.(\ref{uu26})
(here $m_{12}$ is the multiplication of the first and the second tensor arguments; $S_1$ is an
operation of taking the antipode of the first tensor argument, {\it etc.}).

\smallskip
Applications of $m_{23}\circ S_2$ and $m_{23}\circ S_3$ to eq.(\ref{uu19}) establish
eq.(\ref{uu34}).

\smallskip
Given eq.(\ref{uu34}), the statement, that the element $\psi_\cF$ is a  left skew inverse of
the element $\cF$, reads in components:
\be \a\a'\ot S(\b')\b =1\ ,\lb{uupsi}\ee
where primes are used to distinguish different
summations terms, the expression $\a\a'\ot S(\b')\b$ stands for $\sum_{i,j}\a_i\a_j\ot S(\b_j)\b_i$.
Applying $S_1$ to this equation, we find $(S(\a' )\ot S(\b'))\cdot (S(\a )\ot\b)=1$ which means
that the element $S(\a' )\ot S(\b')$ is a left inverse of the element $S(\a )\ot\b$.
However, the latter element is, by eq.(\ref{uu26}), the inverse of $\cF$. Therefore,
eq.(\ref{uu29}) follows.

\smallskip
Applying $S_2$ to eq.(\ref{uupsi}), we find that $\a\ot S^2(\b )$ is a right inverse of the
element $\psi_\cF$.

\smallskip
Applying $S_1^2$ to eq.(\ref{uupsi}) and using eq.(\ref{uu29}), we find that
$S^2(\a )\ot \b $ is a right skew inverse of the element $\cF^{-1}$.

\smallskip
We shall not repeat details for the left inverse of the element $\psi_\cF$ and the left
skew inverse of the element $\cF^{-1}$, calculations are analogous.\hfill$\blacksquare$

\begin{rem}{\rm There is a further generalization of the formulas from the lemma
\ref{lemma3.2.0.2}. Start with the element $\cF$ and alternate operations "take an inverse" and
"take a skew inverse". Then the next operation is always possible, the result is always invertible
and skew invertible. One arrives, after $n$ steps, at $S^n(\a )\ot\b$ if the first operation was
"take an inverse"; if the first operation was "take a skew inverse" then one arrives
at $\a\ot S^n(\b )$ (see \cite{D3}, the section 8). \lb{rrr4}}\end{rem}

\smallskip
{\it From now on, we shall assume that the twisting element $\cF$ is invertible
and satisfies the conditions (\ref{uu1}), (\ref{uu18}) and (\ref{uu19}).}

\paragraph{{\bf 5.}} We turn now to the Hopf algebraic meaning of relations from the subsection \ref{subsec3.2}.

\smallskip
The square of the antipode in an almost cocommutative Hopf algebra, with a universal R-matrix
$\cR =a\ot b$, satisfies the property $S^2(x)=u_{_\cR} x(u_{_\cR})^{-1}$, where $u_{_\cR} =S(b)a$,
for any element $x\in {\cal A}$.
In a matrix representation of an algebra ${\cal A}$, the element $u_{_\cR}$ maps to the matrix
$D_{\hat{\rho}({\cal R})}$
(and the element $S(u_{_\cR} )$ maps to the matrix $C_{\hat{\rho}({\cal R})}$), so an identity (which follows from
eq.(\ref{uu29})$\,$)
$$\begin{array}{rcl}(1\ot u_{_\cR})\;\cF^{-1}\; (1\ot (u_{_\cR})^{-1})&\equiv& (1\ot u_{_\cR} )
(S(\a )\ot\b )(1\ot (u_{_\cR})^{-1}=S(\a )\ot S^2(\b )
\\[1em]
&=&\a\ot S(\b )\equiv\psi_\cF\end{array}$$
becomes one of the relations from the lemma \ref{lemma3.3}. In a similar manner, one can interpret other
relations from the lemma \ref{lemma3.3}.

\medskip
Such an interpretation is not, however, unique. For instance, applying $m_{12}\circ S_2$ to eq.(\ref{uu24})
and using eq.(\ref{uu26}), one finds
$$ \vf\ot 1=\a'\vf S(\a )\ot\b\b'\ ,$$
which, after an application of $S_2$, becomes, due to eqs.(\ref{uu34}) and (\ref{uu29}),
\be \vf\ot 1=\psi_\cF\; (\vf\ot 1)\;\cF\ .\lb{uu39}\ee
Similarly, applying $(\id\ot S)\circ m_{23}\circ\tau_{23}\circ S_3$ (where $\tau$ is the flip,
$\tau (x\ot y)=y\ot x$) to eq.(\ref{uu24}) and using eqs.(\ref{uu29}) and (\ref{uu35}), one finds
$$ 1\ot\vf =\a\a'\ot S(\b')\vf S^2(\b )\ ,$$
which, after an application of $S_1$, becomes, with the help of eq.(\ref{uu29}),
\be 1\ot\vf =\cF\; (1\ot\vf )\;\psi_\cF\ .\lb{uu36}\ee
In the matrix picture, eqs.(\ref{uu39}) and (\ref{uu36}) are also equivalent to particular cases of the
relations from the lemma \ref{lemma3.3} -- but this time we did not use the fact that the square of the
antipode is given by the conjugation by the element $u_{_\cR}$.

\smallskip
Below we shall make use of another version of  eqs.(\ref{uu39}) and (\ref{uu36}).

\smallskip
Writing eqs.(\ref{uu39}) and (\ref{uu36}) as $(\vf\ot 1)\cF^{-1}=\psi_\cF (\vf\ot 1)$ and
$ \cF^{-1}(1\ot\vf) =(1\ot\vf )\psi_\cF$, respectively, and using the expressions
for $\psi_\cF$, $(\psi_\cF )^{-1}$
and $ \cF^{-1}$ from the lemma \ref{lemma3.2.0.2}, we find, in components:
\be\vf S(\a ) \ot\b =\a\vf\ot S(\b ) \ \lb{uu41}\ee
and, respectively,
\be S(\a )\ot\b\vf  =\a\ot\vf S(\b )\ .\lb{uu37}\ee

Applying $S_1$ or $S_2$ to eqs.(\ref{uu41}) and (\ref{uu37}), we obtain corresponding formulas with
$\vf$ replaced by $S(\vf )$. These formulas, together with eqs.(\ref{uu41}) and (\ref{uu37}), imply
\be\begin{array}{ccc} \cF\cdot (\vf S(\vf )\ot 1)&=&(\vf S(\vf )\ot 1)\cdot \cF\ ,\\[1em]
\cF\cdot (1\ot \vf S(\vf ))&=&(1\ot \vf S(\vf ))\cdot \cF\ .\end{array}\lb{fvsv}\ee

It follows, from a compatibility of eqs.(\ref{uu39}) and (\ref{uu36}) (express the element $\psi_\cF$
in terms of $\cF$ and $\vf$ in two ways), that
\be \cF_{12}\cdot (\vf\ot\vf ) =(\vf\ot\vf )\cdot \cF_{12}\ .\lb{fvvvvf}\ee

The relations (\ref{fvsv}) and (\ref{fvvvvf}) are universal analogues of the matrix equalities
(\ref{FCD}) and (\ref{FCC}) (for certain choices of the compatible pairs of the R-matrices)
from the corollary \ref{corollary3.4}.

\paragraph{{\bf 6.}} We shall now establish a Hopf algebraic counterpart of the relation (\ref{XffD}).

\smallskip
We need some information about the element $\vf$. The inverse to the element $\vf$ is given by
eq.(\ref{vfin}); it follows from the lemma \ref{lemma3.2.0.2} that $(\vf )^{-1}=S^2(\a )\b$.

\smallskip
By eq.(\ref{uu29}), one has $S(\vf )=S(\b )\a$ and, then,
$S^2(\vf )=\vf $. Since $S^2(x)=u_{_\cR} x(u_{_\cR})^{-1}$ for any element $x\in {\cal A}$, we conclude that the element
$u_{_\cR}$ commutes with the element $\vf$ and, similarly, with the element $S(\vf )$.

\smallskip
Making the flip in the relations (\ref{uu41}) and (\ref{uu37}), multiplying them out
and comparing, we find that the elements $\vf$ and $S(\vf )$ commute.

\begin{rem}{\rm In fact, more is true. Applying $\id\ot S^j$ to eq.(\ref{uu41}), we obtain
$\vf\a\ot S^{j-1}(\b )=\a\vf\ot S^{j+1}(\b )$ (we used eq.(\ref{uu29}) to rearrange the powers of the antipode).
In a similar way, applying $S^{-j}\ot\id$ to
eq.(\ref{uu37}), we obtain $\a\ot S^{j-1}(\b )\vf =\a\ot\vf S^{j+1}(\b )$. Multiplying out and comparing
the right hand sides, we find that the element $\vf$ commutes with the elements $S^k(\a )\b$
$\forall$ $k\in {\Bbb Z}$.

\smallskip
The same procedure, applied to the flipped versions of eqs.(\ref{uu41}) and (\ref{uu37}) shows
that the element $\vf$ commutes with the elements $S^k(\b )\a$ $\forall$ $k\in {\Bbb Z}$.

\smallskip
Applying the antipode to these commutativity relations, we find that the element $S(\vf )$ commutes with
the elements $S^k(\a )\b$ and $S^k(\b )\a$ $\forall$ $k\in {\Bbb Z}$ as well.}\lb{rrrq}\end{rem}

There is a closed formula for the coproduct of the element $\vf$, again similar to the
standard formula for the coproduct of the element $u_{_\cR}$.

\begin{lem}\lb{lemma3.2.0.3} One has
\be \D (\vf )=\cF_{12}^{-1}\;\cF_{21}^{-1}\cdot (\vf\ot\vf )\ .\lb{deltavf}\ee
\end{lem}

\nin {\bf Proof.~} Together, eqs.(\ref{uu18}) and (\ref{uu19}) imply
\be (\D\ot\D )(\cF )=\cF_{14}\cF_{13}\cF_{24}\cF_{23}\ .\lb{ddf}\ee
Therefore, the coproduct of $\vf$ can be written in a form
\be \D (\vf )=\a_{(1)} S(\b_{(2)})\ot\a_{(2)} S(\b_{(1)} )=\a\a' S(\b\b'')\ot\a'' \vf S(\b')\ \lb{vvs1}\ee
(we use the Sweedler notation for the coproduct, $\D (x)=x_{(1)}\ot
x_{(2)}$ for an element $x\in {\cal A}$).

Using eq.(\ref{uu37}), we continue to rewrite the expression (\ref{vvs1}):
\be \D (\vf )=\a S(\a') S(\b\b'')\ot\a''\b' \vf \ .\lb{vvs2}\ee
The relation (\ref{uu24}), in a form $\cF_{13}\cF_{23}\cF_{12}^{-1}=\cF_{12}^{-1}\cF_{23}\cF_{13}$,
reads, in components,
\be \a S(\a')\ot\a''\b'\ot\b\b''=S(\a )\a''\ot\b\a'\ot\b'\b''\ .\lb{cfmff}\ee
Using eq.(\ref{cfmff}), we transform the right hand side of eq.(\ref{vvs2}) to a form
\be \D (\vf )=S(\a )\a''S(\b'')S(\b')\ot\b\a' \vf =S(\a )\vf S(\b')\ot\b\a' \vf\ .\lb{vvs3}\ee
Using again eq.(\ref{uu37}), we obtain
\be \D (\vf )=S(\a )\b'\vf\ot\b S(\a' ) \vf\ ,\lb{vvs4}\ee
which, by eq.(\ref{uu26}), is a component form of eq.(\ref{deltavf}).
\hfill$\blacksquare$

\medskip
Applying the flip to eq.(\ref{deltavf}), we find
$\D^{{\mathrm{op}}} (\vf )=\cF_{21}^{-1}\;\cF_{12}^{-1}\cdot (\vf\ot\vf )$. Since
$\D^{{\mathrm{op}}} (\vf )\;\cR =\cR\;\D (\vf )$, we conclude
\be (\cR_\cF )_{\cF_{21}}\; (\vf\ot\vf )=(\vf\ot\vf )\;\cR\ .\lb{twitwi}\ee
The translation of the equality (\ref{twitwi}) into the matrix language is equivalent to the
relation (\ref{XffD}) (see the remarks \ref{rrr1} and \ref{rrr2}).

\begin{rem}{\rm It follows from eq.(\ref{deltavf}) that
\be\D (S(\vf ))=(S(\vf )\ot S(\vf ))\cdot\cF_{12}^{-1}\;\cF_{21}^{-1}\ .\lb{svf1}\ee
Eq. (\ref{fvvvvf}), together with eqs.(\ref{deltavf}) and (\ref{svf1}), implies that an
element
\be\varphi :=\vf S(\vf )^{-1}\lb{deffi}\ee
is group-like, $\D (\varphi )=\varphi\ot\varphi$. Therefore,
$S(\varphi )=\varphi^{-1}=S(\vf )(\vf )^{-1}$ but $S(\varphi )=S(\vf S(\vf )^{-1})=(\vf )^{-1}S(\vf )$, which shows
again that $\vf$ commutes with $S(\vf )$.\lb{rrr6}}\end{rem}

\paragraph{{\bf 7.}} The twisted Hopf algebra ${\cal A}_\cF$ is quasitriangular, so we can write the usual identities for its universal
R-matrix $\cR_\cF =\cF_{21}\cR\cF^{-1}$. The relations from the lemma \ref{lemma3.6} are the matrix
counterparts of some of these identities.

\smallskip
{}For the twisted Hopf algebra ${\cal A}_\cF$, one finds, with the help of the first of eqs.(\ref{uu14}),
that $u_{(\cR_\cF)}=\varphi\, u_{_\cR}$, where the element $\varphi$ is defined in eq.(\ref{deffi}) (on the matrix level, this becomes one of the relations (\ref{CDtwist})$\,$). In particular,
\be (S_\cF)^2\, (x)=\varphi\; S^2(x)\;\varphi^{-1}\ .\lb{sfsf}\ee

\medskip
({\it i}) The relation (\ref{R_f-fin}) is a consequence of, for example, the identity
\be (\id\ot S_\cF )((\cR_\cF )^{-1})=\cR_\cF\ .\lb{ee1drf}\ee
We have
\be\begin{array}{rcl} \cR_\cF &=&(\id\ot S_\cF )((\cR_\cF )^{-1})=(\id\ot S_\cF )(\cF\cR^{-1}\cF_{21}^{-1})
=(\id\ot S_\cF )(\a S(a)\b'\ot\b b S(\a'))\\[1em]
&=& \a S(a)\b'\ot \vf S^2(\a')S(b)S(\b )(\vf )^{-1}=\a a\b'\ot \vf S^2(\a')b S(\b )(\vf )^{-1}\ .
\end{array}\lb{le36-1}\ee
Here we used eq.(\ref{uu17}) and the identities from the lemma \ref{lemma3.2.0.2} for $\cF$ and $\cR$.
Applying $S^2\ot S$ to eq.(\ref{uu41}), we find
\be \vf S^2(\a )\ot\b =\a\vf\ot\b\ ,\lb{vsvs1}\ee
since $S^2(\vf )=\vf$.
Using the relation (\ref{vsvs1}) and the relation (\ref{uu37}) in a form
$S(\a )\ot (\vf)^{-1}\b =\a\ot S(\b )(\vf)^{-1}$, we rewrite the last expression in eq.(\ref{le36-1}):
$$ \cR_\cF =S(\a )a\b'\ot \a'\vf b (\vf )^{-1}\b \ $$
or
\be \cR_\cF =\cF_{21}\odot \Bigl( (1\ot\vf )\cR (1\ot\vf )^{-1}\Bigr)\odot\cF^{-1}\ ,\ee
which, on the matrix level, is equivalent to the relation (\ref{R_f-fin}).

\medskip
({\it ii}) Next,
$$\begin{array}{rcl} \psi_{(\cR_\cF )}&=&(\id\ot S_\cF)(\cR_\cF)
=(\id\ot S_\cF)(\cF_{21}\cR\cF^{-1})=(\id\ot S_\cF)(\b a S(\a')\ot\a b\b')\\[1em]
&=&\b a S(\a')\ot \vf S(\b')S(b)S(\a )(\vf )^{-1}
=\b a \a'\ot \vf \b' S(b)S(\a )(\vf )^{-1}\ \end{array}$$
or
\be (1\ot\vf )^{-1}\; \psi_{(\cR_\cF )}\;(1\ot\vf )=\cF\odot\psi_\cR\odot\cF_{21}^{-1}\ ,\ee
which, on the matrix level, is equivalent to the relation (\ref{Psi_R_f}).

\medskip
({\it iii}) To obtain another formula for $\psi_{(\cR_\cF )}$, we start with the identity
$\psi_{(\cR_\cF )}=(\id\ot (S_\cF )^2)((\cR_\cF )^{-1})$, which is a direct consequence of the identities from
the lemma \ref{lemma3.2.0.2}:
\be\begin{array}{rcl}
\psi_{(\cR_\cF )}&=&(\id\ot (S_\cF )^2 )( \cF\cR^{-1}\cF_{21}^{-1})
=(\id\ot (S_\cF )^2 )\Bigl( \a S(a)\b'\ot\b b S(\a')\Bigr) \\[1em]
&=&\a S(a)\b'\ot\varphi S^2(\b )S^2(b)S^3(\a')\varphi^{-1}
=\a a\b'\ot\varphi S^2(\b )S(b)S^3(\a')\varphi^{-1}\\[1em]
&=&\a a\b'\ot S(\vf )^{-1}\b\vf S(b)(\vf)^{-1}S(\a')S(\vf )\ .\end{array}\lb{drfof}\ee
Here we used the identities from the lemma \ref{lemma3.2.0.2}, relations
$\a\ot\vf S^2(\b )=\a\ot\b\vf$ and
$S^3(\a )(\vf )^{-1}\ot\b =(\vf )^{-1}S(\a )\ot\b$, which follow from eqs.(\ref{uu41}) and
(\ref{uu37}), and the formula (\ref{sfsf}) for the square of the twisted antipode.

Eq.(\ref{drfof}) can be rewritten as
\be \bigl( 1\ot S(\vf )\bigr)\; \psi_{(\cR_\cF )}\;\bigl( 1\ot S(\vf )^{-1}\bigr)
=\cF\; (1\ot \vf)\;\psi_\cR\; (1\ot(\vf)^{-1})\;\cF_{21}^{-1}\ ,\ee
which, in the matrix picture, is equivalent to eq.(\ref{Psi_R_f-another}).

\medskip
({\it iv}) The property $(S_\cF\ot S_\cF )(\cR_\cF )=\cR_\cF$ leads to
\be (\vf\ot\vf )\;\cF^{-1}\;\cR\;\cF_{21}= \cF_{21}\;\cR\;\cF^{-1}\; (\vf\ot\vf )\ .\lb{anom}\ee
Since the twisting element $\cF$ commutes with $\vf\ot\vf$, the formula (\ref{anom}) is another
manifestation of the relation (\ref{XffD}).

\begin{rem}{\rm We conclude the appendix with several more properties of the group-like
element $\varphi$ defined in eq.(\ref{deffi}).

\medskip
We have $\cR\cdot (\varphi\ot\varphi )=(\varphi\ot\varphi )\cdot\cR$. To see this, apply the antipode to
the relation (\ref{twitwi}) and compare with the same relation (\ref{twitwi}).

\medskip
Recall that a quasitriangular Hopf algebra ${\cal A}$ is called a ribbon Hopf algebra if it contains
a ribbon element $r$, that is, a central element such that $r^2=u_{_\cR} S(u_{_\cR} )$ and
$\D (r)=\cR_{12}^{-1}\cR_{21}^{-1}\cdot (r\ot r)$ (see \cite{Resh}, or \cite{CP}, the chapter 4).
The twisted algebra ${\cal A}_\cF$ is a ribbon Hopf algebra if the algebra ${\cal A}$ is; for the ribbon element
of the algebra ${\cal A}_\cF$, one can choose $r_{_\cF} = \varphi r$.
}\end{rem}

\addtocontents{toc}{\contentsline {section}{\numberline {} References}{\pageref{refer}}}

\end{document}